%% file: contact.tex
\documentclass[11pt,a4paper]{article}
\usepackage{amsmath,amssymb,amsthm} 
\input rlepsf

\numberwithin{equation}{section}

\def\co{\colon\thinspace}
\newcommand{\R}{{\mathbb R}}
\newcommand{\C}{{\mathbb C}}
\newcommand{\Z}{{\mathbb Z}}
\newcommand{\SL}{\mbox{\rm SL}}
\newcommand{\SO}{\mbox{\rm SO}}

\newcommand{\U}{\mbox{\rm U}}

\newcommand{\Sp}{\mbox{\rm Sp}}

\newcommand{\ozeta}{\overline{\zeta}}
\newcommand{\oz}{\overline{z}}
\newcommand{\ow}{\overline{w}}

\newtheorem{thm}{Theorem}[section]
\newtheorem{lem}[thm]{Lemma}
\newtheorem{prop}[thm]{Proposition}
\newtheorem{cor}[thm]{Corollary}
\newtheorem{defn}[thm]{Definition}
\newtheorem{exam}[thm]{Example}
\newtheorem{rem}[thm]{Remark}

\begin{document}

\title{\bf Contact Geometry
\\{\small\rm To appear in {\it Handbook
of Differential Geometry}, vol.~2}\\{\small\rm (F.J.E.\ Dillen
and L.C.A.\ Verstraelen, eds.)}
}
\author{Hansj\"org Geiges}
\date{\small\it Mathematisches Institut, Universit\"at zu K\"oln,\\
Weyertal 86--90, 50931 K\"oln, Germany\\
E-mail: geiges@math.uni-koeln.de}

\maketitle

\tableofcontents

\newpage


\section{Introduction}
Over the past two decades, contact geometry has undergone a veritable
metamorphosis: once the ugly duckling known as `the odd-dimensional
analogue of symplectic geometry', it has now evolved into a proud field of
study in its own right. As is typical for a period of rapid development
in an area of mathematics, there are a fair number of folklore results
that every mathematician working in the area knows, but no references
that make these results accessible to the novice. I therefore take
the present article as an opportunity to take stock of some of that folklore.

There are many excellent surveys covering specific aspects of contact
geometry (e.g.\ classification questions in dimension~$3$, dynamics
of the Reeb vector field, various notions of symplectic fillability,
transverse and Legendrian knots and links). All these topics deserve
to be included in a comprehensive survey, but an attempt to do so
here would have left this article in the `to appear' limbo for much
too long.

Thus, instead of adding yet another survey, my plan here is to
cover in detail some of the more fundamental differential topological
aspects of contact geometry. In doing so, I have not tried to hide
my own idiosyncrasies and preoccupations. Owing to a relatively leisurely
pace and constraints of the present format,
I have not been able to cover quite as much
material as I should have wished. None the less, I hope that the reader
of the present handbook chapter will be better prepared to study some of
the surveys I alluded to -- a guide to these surveys will be
provided -- and from there
to move on to the original literature.

A book chapter with comparable aims is Chapter~8 in~\cite{abklr94}.
It seemed opportune to be brief on topics that are covered extensively
there, even if it is done at the cost of leaving out
some essential issues. I hope to return to the material of the present
chapter in a yet to be written more comprehensive monograph.

\vspace{2mm}

\noindent {\bf Acknowledgements.} I am grateful to Fan Ding and
Jes\'us Gonzalo for
their attentive reading of the original manuscript.
I also thank John Etnyre and Stephan Sch\"onenberger
for allowing me to use a couple of their figures
(viz., Figures~\ref{figure:charfolonS2}
and~\ref{figure:ex1} of the present text, respectively).

\input{section2}

\input{section3}

\input{section4}



\end{document}

%% file: section2.tex
\section{Contact manifolds}
\label{section:conman}
Let $M$ be a differential manifold and $\xi\subset TM$ a field of
hyperplanes on $M$. Locally such a hyperplane field can always be
written as the kernel of a non-vanishing $1$--form~$\alpha$. One way
to see this is to choose an auxiliary Riemannian metric $g$ on $M$ and
then to define $\alpha =g (X,.)$, where $X$ is a local non-zero
section of the line bundle $\xi^{\perp}$ (the orthogonal complement
of~$\xi$ in $TM$). We see that the existence of a globally defined
$1$--form $\alpha$ with $\xi =\ker\alpha$ is equivalent to the
orientability (hence triviality) of $\xi^{\perp}$, i.e.\ the
coorientability of~$\xi$. Except for an example below, I shall
always assume this condition.

If $\alpha$ satisfies the Frobenius integrability condition
\[ \alpha\wedge d\alpha =0,\]
then $\xi$ is an integrable hyperplane field (and vice versa), and its
integral submanifolds form a codimension~$1$ foliation of~$M$. Equivalently,
this integrability condition can be written as
\[ X,Y\in\xi\Longrightarrow [X,Y]\in\xi.\]
An integrable hyperplane field is locally of the form $dz=0$, where $z$
is a coordinate function on~$M$. Much is known, too, about the global
topology of foliations, cf.~\cite{tamu92}.

Contact structures are in a certain sense the exact opposite of
integrable hyperplane fields.

\begin{defn}
Let $M$ be a manifold of odd dimension $2n+1$. A {\bf contact structure}
is a maximally non-integrable hyperplane field $\xi =\ker\alpha
\subset TM$, that is, the defining $1$--form $\alpha$ is required
to satisfy
\[ \alpha\wedge (d\alpha )^n\neq 0\]
(meaning that it vanishes nowhere). Such a $1$--form $\alpha$ is
called a {\bf contact form}. The pair $(M,\xi )$ is called a
{\bf contact manifold}.
\end{defn}

\begin{rem}
\label{rem:posneg}
{\rm
Observe that in this case $\alpha\wedge (d\alpha)^n$ is a volume form
on~$M$; in particular, $M$ needs to be orientable.
The condition $\alpha\wedge (d\alpha )^n\neq 0$ is independent
of the specific choice of $\alpha$ and thus is indeed a property of
$\xi=\ker\alpha$: Any other $1$--form defining the same hyperplane
field must be of the form $\lambda\alpha$ for some smooth function
$\lambda\co M\rightarrow {\mathbb R}\setminus\{ 0\}$, and we have
\[ (\lambda\alpha )\wedge (d(\lambda\alpha ))^n=\lambda\alpha\wedge
(\lambda\, d\alpha +d\lambda\wedge\alpha )^n=\lambda^{n+1}
\alpha\wedge (d\alpha )^n\neq 0.\]

We see that if $n$ is odd,
the sign of this volume form depends only on~$\xi$, not the choice
of~$\alpha$. This makes it possible, given an orientation of~$M$,
to speak of {\it positive} and {\it negative} contact structures.}
\end{rem}

\begin{rem}
\label{rem:symplectic}
{\rm An equivalent formulation of the contact condition is that
we have $(d\alpha )^n|_{\xi}\neq 0$. In particular, for every point $p\in M$,
the $2n$--dimensional subspace $\xi_p\subset T_pM$ is a vector space on which
$d\alpha$ defines a skew-symmetric form of maximal rank, that is,
$(\xi_p, d\alpha |_{\xi_p})$ is a {\it symplectic} vector space.
A consequence of this fact is that there exists a complex bundle
structure $J\co\xi\rightarrow\xi$ compatible with $d\alpha$
(see~\cite[Prop.~2.63]{mcsa98}), i.e.\ a bundle endomorphism satisfying
\begin{itemize}
\item $J^2=-\mbox{\rm id}_{\xi}$,
\item $d\alpha (JX,JY)=d\alpha (X,Y)$ for all $X,Y\in\xi$,
\item $d\alpha (X,JX)>0$ for $0\neq X\in\xi$.
\end{itemize}
}
\end{rem}

\begin{rem}
{\rm
The name `contact structure' has its origins in the fact that one
of the first historical sources of contact manifolds are the so-called
spaces of contact elements (which in fact have to do with `contact'
in the differential geometric sense), see \cite{arno78} and~\cite{geig01}.
}
\end{rem}

In the $3$--dimensional case the contact condition can also be
formulated as
\[ X,Y\in\xi \; \mbox{\rm linearly independent}\Longrightarrow
[X,Y]\not\in\xi;\]
this follows immediately from the equation
\[ d\alpha (X,Y)=X(\alpha (Y))-Y(\alpha (X))-\alpha ([X,Y])\]
and the fact that the contact condition (in dim.~$3$) may be
written as $d\alpha |_{\xi}\neq 0$.

In the present article I shall take it for granted that contact
structures are worthwhile objects of study. As I hope to illustrate,
this is fully justified by the beautiful mathematics to which
they have given rise. For an apology of contact structures
in terms of their origin (with hindsight) in physics and the multifarious
connections with other areas of mathematics I refer the reader to
the historical surveys \cite{lutz88} and~\cite{geig01}. Contact structures
may also be justified on the grounds that they are generic objects:
A generic $1$--form $\alpha$ on an odd-dimensional manifold satisfies
the contact condition outside a smooth hypersurface, see~\cite{mart70}.
Similarly, a generic $1$--form $\alpha$ on a $2n$--dimensional
manifold satisfies the condition $\alpha\wedge (d\alpha )^{n-1}\neq 0$
outside a submanifold of codimension~$3$; such `even-contact manifolds'
have been studied in \cite{ginz92}, for instance, but on the whole
their theory is not as rich or well-motivated as that of contact
structures.

\begin{defn}
Associated with a contact form $\alpha$ one has the so-called
{\bf Reeb vector field}~$R_{\alpha}$, defined by the
equations
\begin{itemize}
\item[(i)] $d\alpha (R_{\alpha},.)\equiv 0$,
\item[(ii)] $\alpha (R_{\alpha})\equiv 1$.
\end{itemize}
\end{defn}

As a skew-symmetric form of maximal rank~$2n$, the form $d\alpha|_{T_pM}$
has a $1$--dimensional kernel for each $p\in M^{2n+1}$. Hence equation~(i)
defines a unique line field $\langle R_{\alpha}\rangle$ on~$M$. The
contact condition $\alpha\wedge (d\alpha )^n\neq 0$ implies that
$\alpha$ is non-trivial on that line field, so a global vector field
is defined by the additional normalisation condition~(ii).
\subsection{Contact manifolds and their submanifolds}
We begin with some examples of contact manifolds; the simple verification
that the listed $1$--forms are contact forms is left to
the reader.

\begin{exam}
\label{exam:standardR}
{\rm
On ${\mathbb R}^{2n+1}$ with cartesian coordinates $(x_1,y_1,
\ldots ,x_n,y_n,z)$, the $1$--form
\[ \alpha_1 =dz+\sum_{i=1}^nx_i\, dy_i\]
is a contact form.  }
\end{exam}

\begin{figure}[h]
\centerline{\relabelbox\small
\epsfysize 2truein \epsfbox{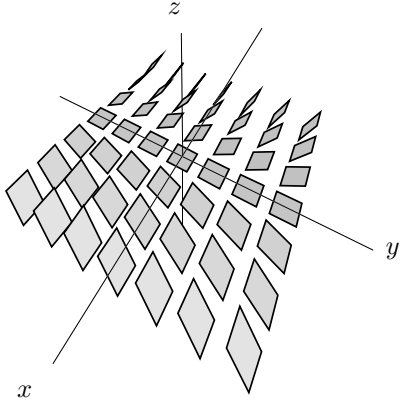}
\extralabel <-4.8cm,.1cm> {$x$}
\extralabel <.1cm,2.0cm> {$y$}
\extralabel <-2.8cm,5.2cm> {$z$}
\endrelabelbox}
\caption{The contact structure $\ker (dz+x\, dy)$.}
\label{figure:ex1}
\end{figure}

\begin{exam}
\label{exam:standard2}
{\rm On ${\mathbb R}^{2n+1}$ with polar coordinates $(r_j,\varphi_j)$
for the $(x_j,y_j)$--plane, $j=1,\ldots ,n$, the $1$--form
\[\alpha_2=dz+\sum_{j=1}^nr_j^2\, d\varphi_j=
dz+\sum_{j=1}^n(x_j\, dy_j-y_j\, dx_j)\]
is a contact form.  }
\end{exam}

\begin{defn}
Two contact manifolds $(M_1,\xi_1)$ and $(M_2,\xi_2)$ are called
{\bf contactomorphic} if there is a diffeomorphism $f\co M_1\rightarrow
M_2$ with $Tf(\xi_1)=\xi_2$, where $Tf\co TM_1\rightarrow TM_2$
denotes the differential of~$f$. If $\xi_i=\ker\alpha_i$, $i=1,2$,
this is equivalent to the existence of a nowhere zero function
$\lambda\co M_1\rightarrow {\mathbb R}$ such that $f^*\alpha_2=
\lambda\alpha_1$.
\end{defn}

\begin{exam}
{\rm The contact manifolds $({\mathbb R}^{2n+1},\xi_i=\ker\alpha_i)$, $i=1,2$,
from the preceding examples are contactomorphic. An explicit contactomorphism
$f$ with $f^*\alpha_2=\alpha_1$ is given by
\[ f(x,y,z)=\bigl( (x+y)/2,(y-x)/2,z+xy/2\bigr) ,\]
where $x$ and $y$ stand for $(x_1,\ldots ,x_n)$ and $(y_1,\ldots ,
y_n)$, respectively, and $xy$ stands for $\sum_jx_jy_j$.
Similarly, both these contact structures are
contactomorphic to $\ker (dz-\sum_j y_j\, dx_j)$. Any of these
contact structures is called the {\bf standard contact structure
on} ${\mathbb R}^{2n+1}$.
}
\end{exam}

\begin{exam}
\label{exam:standardS}
{\rm The {\bf standard contact structure on the unit sphere}
$S^{2n+1}$ in ${\mathbb R}^{2n+2}$ (with cartesian coordinates
$(x_1,y_1,\ldots ,x_{n+1},y_{n+1})$) is defined by the contact form
\[ \alpha_0 =\sum_{j=1}^{n+1}(x_j\, dy_j-y_j\, dx_j).\]
With $r$ denoting the radial coordinate on ${\mathbb R}^{2n+1}$
(that is, $r^2=\sum_j(x_j^2+y_j^2)$) one checks easily that $\alpha_0\wedge
(d\alpha_0)^n\wedge r\, dr\neq 0$ for $r\neq 0$. Since $S^{2n+1}$
is a level surface of $r$ (or~$r^2$), this verifies the contact
condition.

Alternatively, one may regard $S^{2n+1}$ as the unit sphere in
${\mathbb C}^{n+1}$
with complex structure~$J$ (corresponding to complex coordinates
$z_j=x_j+iy_j$, $j=1,\ldots ,n+1$).
Then $\xi_0=\ker\alpha_0$ defines at each point
$p\in S^{2n+1}$ the complex (i.e.\ $J$--invariant) subspace of
$T_pS^{2n+1}$, that is,
\[ \xi_0=TS^{2n+1}\cap J(TS^{2n+1}).\]
This follows from the observation that $\alpha =-r\, dr\circ J$. The
hermitian form $d\alpha (.,J.)$ on $\xi_0$ is called the {\it Levi form}
of the hypersurface $S^{2n+1}\subset {\mathbb C}^{n+1}$. The contact condition
for $\xi$ corresponds to the positive definiteness of that Levi form,
or what in complex analysis is called the {\it strict pseudoconvexity}
of the hypersurface. For more on the question of pseudoconvexity
from the contact geometric viewpoint see~\cite[Section~8.2]{abklr94}. Beware
that the `complex structure' in their Proposition~8.14 is not required to
be integrable, i.e.\ constitutes what is more commonly referred to
as an `almost complex structure'.
}
\end{exam}

\begin{defn}
Let $(V,\omega )$ be a {\bf symplectic manifold} of dimension $2n+2$, that is,
$\omega$ is a closed ($d\omega =0$) and non-degenerate ($\omega^{n+1}
\neq 0$) $2$--form on~$V$. A vector field $X$ is called a {\bf Liouville
vector field} if ${\mathcal L}_X\omega =\omega$, where ${\mathcal L}$
denotes the Lie derivative.
\end{defn}

With the help of Cartan's formula
${\mathcal L}_X=d\circ i_X+i_X\circ d$ this may be rewritten
as $d(i_X\omega )=\omega$. Then the $1$--form $\alpha=i_X\omega$ defines
a contact form on any hypersurface $M$ in $V$ transverse to~$X$. Indeed,
\[ \alpha\wedge (d\alpha )^n=i_X\omega\wedge (d(i_X\omega ))^n
=i_X\omega\wedge \omega^n =\frac{1}{n+1}i_X(\omega^{n+1}),\]
which is a volume form on $M\subset V$ provided $M$ is transverse
to~$X$.

\begin{exam}
{\rm With $V={\mathbb R}^{2n+2}$, symplectic form $\omega =\sum_j dx_j\wedge
dy_j$, and Liouville vector field $X=\sum_j(x_j\partial_{x_j}+
y_j\partial_{y_j})/2=r\partial_r/2$, we recover the standard contact structure
on~$S^{2n+1}$.}
\end{exam}

For finer issues relating to hypersurfaces in symplectic manifolds
transverse to a Liouville vector field I refer the reader
to~\cite[Section~8.2]{abklr94}.

Here is a further useful example of contactomorphic manifolds.

\begin{prop}
\label{prop:standard-S-R}
For any point $p\in S^{2n+1}$, the manifold $(S^{2n+1}\setminus\{ p\},
\xi_0)$ is contactomorphic to $({\mathbb R}^{2n+1},\xi_2)$.
\end{prop}

\begin{proof}
The contact manifold $(S^{2n+1},\xi_0)$ is a homogeneous space under
the natural $\U (n+1)$--action, so we are free to choose $p=(0,\ldots ,
0,-1)$. Stereographic projection from $p$ does almost, but not quite
yield the desired contactomorphism. Instead, we use a map that is
well-known in the theory of Siegel domains (cf.~\cite[Chapter~8]{akhi90})
and that looks a bit like a complex analogue of stereographic
projection; this was suggested in \cite[Exercise 3.64]{mcsa98}.

Regard $S^{2n+1}$ as the unit sphere in ${\mathbb C}^{n+1}=
{\mathbb C}^n\times {\mathbb C}$ with cartesian coordinates
$(z_1,\ldots ,z_n,w)=(z,w)$. We identify ${\mathbb R}^{2n+1}$
with ${\mathbb C}^n\times{\mathbb R}\subset {\mathbb C}^n
\times {\mathbb C}$ with coordinates
$(\zeta_1,\ldots ,\zeta_n,s)=(\zeta ,s)=(\zeta ,\mbox{\rm Re}\,\sigma )$,
where $\zeta_j=x_j+iy_j$. Then
\begin{eqnarray*}
\alpha_2 & = & ds+\sum_{j=1}^n(x_j\, dy_j-y_j\, dx_j)\\
         & = & ds+\frac{i}{2}(\zeta\, d\ozeta -\ozeta\, d\zeta ).
\end{eqnarray*}
and
\[ \alpha_0=\frac{i}{2}(z\, d\oz -\oz\, dz+w\, d\ow -\ow\, dw).\]
Now define a smooth map $f\co S^{2n+1}\setminus \{ (0,-1)\}$ by
\[ (\zeta ,s)=f(z,w)=\left( \frac{z}{1+w},
-\frac{i(w-\ow)}{2|1+w|^2}\right) .\]
Then
\begin{eqnarray*}
f^*ds & = & -\frac{i\, dw}{2|1+w|^2}+
            \frac{i\, d\ow}{2|1+w|^2}\\
      &   & \mbox{}+\frac{i(w-\ow)}{2(1+w)}\frac{dw}{|1+w|^2}
            +\frac{i(w-\ow)}{2(1+\ow)}\frac{d\ow}{|1+w|^2}\\
      & = & \frac{i}{2|1+w|^2}\left(-dw+d\ow+
            \frac{w-\ow}{1+w}dw+\frac{w-\ow}{1+\ow}d\ow\right)
\end{eqnarray*}
and
\begin{eqnarray*}
f^*(\zeta\, d\ozeta -\ozeta\, d\zeta ) & = & \frac{z}{1+w}\left(
            \frac{d\oz}{1+\ow}-\frac{\oz}{(1+\ow )^2}d\ow\right)\\
      &   & \mbox{}-\frac{\oz}{1+\ow}\left( \frac{dz}{1+w}-
            \frac{z}{(1+w)^2}dw\right)\\
      & = & \frac{1}{|1+w|^2}\left( z\, d\oz-\oz dz+|z|^2
             \left( \frac{dw}{1+w}-\frac{d\ow}{1+\ow}\right)\right).
\end{eqnarray*}
Along $S^{2n+1}$ we have
\begin{eqnarray*}
|z|^2=1-|w|^2 & = & (1-w)(1+\ow )+(w-\ow )\\
              & = & (1-\ow )(1+w)-(w-\ow ),
\end{eqnarray*}
whence
\begin{eqnarray*}
|z|^2\left(\frac{dw}{1+w}-\frac{d\ow}{1+\ow}\right) & = &
       (1-\ow)\, dw-\frac{w-\ow}{1+w}dw\\
  &  & \mbox{}-(1-w)\, d\ow-\frac{w-\ow}{1+\ow}d\ow.
\end{eqnarray*}
From these calculations we conclude $f^*\alpha_2=\alpha_0/|1+w|^2$. So it
only remains to show that $f$ is actually a diffeomorphism of
$S^{2n+1}\setminus\{ (0,-1)\}$ onto ${\mathbb R}^{2n+1}$. To that end,
consider the map
\[ \widetilde{f}\co ({\mathbb C}^n\times {\mathbb C})\setminus
({\mathbb C}^n\times\{ -1\})\longrightarrow ({\mathbb C}^n\times {\mathbb C})
\setminus ({\mathbb C}^n\times\{ -i/2\})\]
defined by
\[ (\zeta ,\sigma )=\widetilde{f}(z,w)=\left(\frac{z}{1+w},
-\frac{i}{2}\frac{w-1}{w+1}\right).\]
This is a biholomorphic map with inverse map
\[ (\zeta ,\sigma )\longmapsto \left(\frac{2\zeta }{1-2i\sigma },\frac{1+2i
\sigma}{1-2i\sigma}\right) .\]
We compute
\begin{eqnarray*}
\mbox{\rm Im}\,\sigma & = & -\frac{w-1}{4(w+1)}-
                          \frac{\ow-1}{4(\ow+1)}\\
  & = & -\frac{(w-1)(\ow+1)+(\ow-1)(w+1)}{4|1+w|^2}\\
  & = & \frac{1-|w|^2}{2|1+w|^2}.
\end{eqnarray*}
Hence for $(z,w)\in S^{2n+1}\setminus\{ (0,-1)\}$ we have
\[ \mbox{\rm Im}\,\sigma =\frac{|z|^2}{2|1+w|^2}=\frac{1}{2}
|\zeta |^2;\]
conversely, any point $(\zeta ,\sigma )$ with $\mbox{\rm Im}\,\sigma =
|\zeta |^2/2$ lies in the image of $\widetilde{f}
|_{S^{2n+1}\setminus\{ (0,-1)\}}$, that is, $\widetilde{f}$
restricted to $S^{2n+1}\setminus\{ (0,-1)\}$ is a diffeomorphism onto
$\{ \mbox{\rm Im}\,\sigma =|\zeta |^2/2\}$.
Finally, we compute
\begin{eqnarray*}
\mbox{\rm Re}\,\sigma & = & -\frac{i(w-1)}{4(w+1)}+
                          \frac{i(\ow-1)}{4(\ow+1)}\\
   & = & -i\frac{(w-1)(\ow+1)-(\ow-1)(w+1)}{4|1+w|^2}\\
   & = & -\frac{i(w-\ow)}{2|1+w|^2},
\end{eqnarray*}
from which we see that for $(z,w)\in S^{2n+1}\setminus\{ (0,-1)\}$
and with $(\zeta ,\sigma )=\widetilde{f}(z,w)$ we have
$f(z,w)=(\zeta ,\mbox{\rm Re}\,\sigma )$. This concludes the proof.
\end{proof}

At the beginning of this section I mentioned that one may allow
contact structures that are not coorientable, and hence not defined
by a global contact form.

\begin{exam}
{\rm Let $M={\mathbb R}^{n+1}\times {\mathbb R}P^n$ with
cartesian coordinates $(x_0,\ldots ,x_n)$ on the ${\mathbb R}^{n+1}$--factor
and
homogeneous coordinates $[y_0:\ldots :y_n]$ on the ${\mathbb R}P^n$--factor.
Then
\[ \xi =\ker \bigl( \sum_{j=0}^ny_j\, dx_j\bigr) \]
is a well-defined hyperplane field on $M$, because the $1$--form
on the right-hand side is well-defined up to scaling by a
non-zero real constant.
On the open submanifold $U_k=\{ y_k\neq 0\}\cong {\mathbb R}^{n+1}\times
{\mathbb R}^n$ of $M$ we have $\xi =\ker \alpha_k$ with
\[ \alpha_k =dx_k+\sum_{j\neq k}\left(\frac{y_j}{y_k}\right)\, dx_j\]
an honest $1$--form on~$U_k$. This is the standard contact form
of Example~\ref{exam:standardR}, which proves that $\xi$ is a contact
structure on~$M$.

If $n$ is even, then $M$ is not orientable, so there can be no global
contact form defining~$\xi$ (cf.\ Remark~\ref{rem:posneg}), i.e.\
$\xi$ is {\it not coorientable}. Notice,
however, that a contact structure on a manifold of dimension $2n+1$
with $n$ even is always {\it orientable}: the sign of $(d\alpha )^n|_{\xi}$
does not depend on the choice of local $1$--form defining~$\xi$.

If $n$ is odd, then $M$ is orientable, so it would be possible that
$\xi$ is the kernel of a globally defined $1$--form. However, since the
sign of $\alpha\wedge (d\alpha )^n$, for $n$ odd,
is independent of the choice of local $1$--form defining~$\xi$, it is also
conceivable that no global contact form exists. (In fact, this
consideration shows that any manifold of dimension $2n+1$, with $n$ odd,
admitting a contact structure (coorientable or not) needs to
be orientable.) This is indeed what happens,
as we shall prove now.
}
\end{exam}

\begin{prop}
Let $(M,\xi )$ be the contact manifold of the preceding example. Then
$TM/\xi$ can be identified with the canonical line bundle on~${\mathbb R}
P^n$ (pulled back to~$M$). In particular, $TM/\xi$ is a non-trivial line
bundle, so $\xi$ is not coorientable.
\end{prop}

\begin{proof}
For given $y=[y_0:\ldots :y_n]\in {\mathbb R}P^n$,
the vector $y_0\partial_{x_0}+\cdots +y_n\partial_{x_n}\in
T_x\R^{n+1}$ is well-defined
up to a non-zero real factor (and independent of $x\in\R^{n+1}$),
and hence defines a line
$\ell_y$ in~$T_x\R^{n+1}\cong {\mathbb R}^{n+1}$. The set
\begin{eqnarray*}
E & = & \{ (t,x,y)\co x\in {\mathbb R}^{n+1},\; y\in
{\mathbb R}P^n,\; t\in\ell_y \}\\
  & \subset &  T\R^{n+1}\times\R P^n
\subset T(\R^{n+1}\times\R P^n)=TM
\end{eqnarray*}
with projection $(t,x,y)\mapsto (x,y)$ defines a line sub-bundle
of $TM$ that restricts to the canonical line bundle
over $\{ x\}\times {\mathbb R}P^n\equiv {\mathbb R}P^n$ for each
$x\in {\mathbb R}^{n+1}$. The canonical line bundle over
${\mathbb R}P^n$ is well-known to be
non-trivial~\cite[p.~16]{mist74}, so the same holds for~$E$.

Moreover, $E$ is clearly complementary to~$\xi$, i.e.\
$TM/\xi\cong E$, since
\[ \sum_{j=0}^n y_j\, dx_j (\sum_{k=0}^ny_k\partial_{x_k})=
\sum_{j=0}^n y_j^2\neq 0.\]
This proves that that $\xi$ is not coorientable.
\end{proof}

To sum up, in the example above we have one of the following two
situations:
\begin{itemize}
\item If $n$ is odd, then $M$ is orientable; $\xi$ is neither orientable
nor coorientable.
\item If $n$ is even, then $M$ is not orientable; $\xi$ is not coorientable,
but it is orientable.
\end{itemize}

We close this section with the definition of the most important
types of submanifolds.

\begin{defn}
Let $(M,\xi )$ be a contact manifold.

(i) A submanifold $L$ of $(M, \xi )$ is called an {\bf isotropic} submanifold
if $T_xL\subset\xi_x$ for all $x\in  L$.

(ii) A submanifold $M'$ of $M$ with contact structure $\xi'$ is called
a {\bf contact submanifold} if $TM'\cap\xi |_{M'} =\xi'$.
\end{defn}

Observe that if $\xi =\ker\alpha$ and $i\co M'\rightarrow M$ denotes the
inclusion map, then the condition for $(M',\xi')$ to be a contact
submanifold of $(M,\xi )$ is that $\xi'=\ker (i^*\alpha )$. In particular,
$\xi'\subset \xi|_{M'}$ is a symplectic sub-bundle with respect to
the symplectic bundle structure on $\xi$ given by~$d\alpha$.

The following is a manifestation of the maximal non-integrability
of contact structures.

\begin{prop}
Let $(M,\xi )$ be a contact manifold of dimension $2n+1$ and $L$
an isotropic submanifold. Then $\dim L\leq n$.
\end{prop}

\begin{proof}
Write $i$ for the inclusion of $L$ in $M$ and let $\alpha$ be an
(at least locally defined)
contact form defining~$\xi$. Then the condition for $L$ to be
isotropic becomes $i^*\alpha\equiv 0$. It follows that $i^*d\alpha \equiv 0$.
In particular, $T_pL\subset\xi_p$ is an isotropic subspace of
the symplectic vector space $(\xi_p ,d\alpha |_{\xi_p})$, i.e.\
a subspace on which the symplectic form restricts to zero. From Linear
Algebra we know that this implies $\dim T_pL\leq (\dim\xi_p)/2=n$.
\end{proof}

\begin{defn}
An isotropic submanifold $L\subset (M^{2n+1},\xi )$ of maximal
possible dimension $n$ is called a {\bf Legendrian submanifold}.
\end{defn}

In particular, in a $3$--dimensional contact manifold there are two
distinguished types of knots: {\bf Legendrian knots} on the one hand,
{\bf transverse\footnote{Some people like to call them
`transversal knots', but I adhere to J.H.C. Whitehead's
dictum, as quoted in~\cite{hirs76}:
``{\it Transversal} is a noun; the adjective is
{\it transverse}.''} knots}
on the other, i.e.\ knots that are everywhere
transverse to the contact structure. If $\xi$ is cooriented by a contact
form $\alpha$ and $\gamma\co S^1\rightarrow (M,\xi =\ker\alpha )$ is
oriented, one can speak of a {\it positively} or {\it negatively}
transverse knot, depending on whether $\alpha (\dot{\gamma})>0$
or $\alpha (\dot{\gamma})<0$.

\subsection{Gray stability and the Moser trick}
The Gray stability theorem that we are going to prove in this section says
that there are no non-trivial deformations of contact structures on
closed manifolds. In fancy language, this means that contact structures
on closed manifolds have discrete moduli. First a preparatory lemma.

\begin{lem}
\label{lem:Lie}
Let $\omega_t$, $t\in [0,1]$, be a smooth family of differential $k$--forms
on a manifold $M$ and $(\psi_t)_{t\in [0,1]}$ an isotopy of~$M$.
Define a time-dependent vector field $X_t$ on $M$ by
$X_t\circ\psi_t=\dot{\psi}_t$, where the dot denotes derivative
with respect to~$t$ (so that $\psi_t$ is the flow of~$X_t$).
Then
\[ \frac{d}{dt}\bigl( \psi_t^*\omega_t\bigr)
=\psi_t^*\bigl( \dot{\omega}_t+
{\mathcal L}_{X_t}\omega_t\bigr) .\]
\end{lem}

\begin{proof}
For a time-independent $k$--form $\omega$ we have
\[ \frac{d}{dt}\bigl( \psi_t^*\omega \bigr) 
=\psi_t^*\bigl( {\mathcal L}_{X_t}\omega
\bigr) .\]
This follows by observing that
\begin{itemize}
\item[(i)] the formula holds for functions,
\item[(ii)] if it holds for differential forms $\omega$ and~$\omega'$,
then also for $\omega\wedge\omega'$,
\item[(iii)] if it holds for $\omega$, then also for~$d\omega$,
\item[(iv)] locally functions and differentials of functions generate
the algebra of differential forms.
\end{itemize}

We then compute
\begin{eqnarray*}
\frac{d}{dt}(\psi_t^*\omega_t) & = & \lim_{h\rightarrow 0}
\frac{\psi_{t+h}^*\omega_{t+h}-\psi_t^*\omega_t}{h}\\
  & = & \lim_{h\rightarrow 0}\frac{\psi_{t+h}^*\omega_{t+h}-
\psi_{t+h}^*\omega_t + \psi_{t+h}^*\omega_t -\psi_t^*\omega_t}{h}\\
  & = & \lim_{h\rightarrow 0}\psi_{t+h}^*\left( \frac{\omega_{t+h}-\omega_t}{h}
\right) +\lim_{h\rightarrow 0}\frac{\psi_{t+h}^*\omega_t-\psi_t^*\omega_t}{h}\\
  & = & \psi_t^*\bigl( \dot{\omega}_t+{\mathcal L}_{X_t}\omega_t\bigr) .
\qed
\end{eqnarray*}
\renewcommand{\qed}{}
\end{proof}

\begin{thm}[Gray stability]
\label{thm:gray}
Let $\xi_t$, $t\in [0,1]$, be a smooth family of contact structures on
a closed manifold~$M$. Then there is an isotopy $(\psi_t)_{t\in [0,1]}$ of $M$
such that
\[ T\psi_t(\xi_0)=\xi_t\;\;\mbox{\rm for each}\; t\in [0,1].\]
\end{thm}

\begin{proof}
The simplest proof of this result rests on what is known as the
{\bf Moser trick}, introduced by J.~Moser~\cite{mose65} in the context
of stability results for (equicohomologous) volume
and symplectic forms. J.~Gray's original proof~\cite{gray59} was based on
deformation theory \`a la Kodaira-Spencer. The idea of the Moser trick
is to assume that $\psi_t$ is the flow of a time-dependent vector
field~$X_t$. The desired equation for $\psi_t$ then translates
into an equation for $X_t$. If that equation can be solved, the
isotopy $\psi_t$ is found by integrating~$X_t$; on a closed
manifold the flow of $X_t$ will be globally defined.

Let $\alpha_t$ be a smooth family of $1$--forms with $\ker\alpha_t
=\xi_t$. The equation in the theorem then translates into
\[ \psi_t^*\alpha_t=\lambda_t\alpha_0,\]
where $\lambda_t\co M\rightarrow {\mathbb R}^+$ is a suitable
smooth family of smooth functions. Differentiation of this equation
with respect to $t$ yields, with the help of the preceding lemma,
\[ \psi_t^*\bigl( \dot{\alpha}_t+{\mathcal L}_{X_t}\alpha_t
\bigr) = \dot{\lambda}_t\alpha_0=\frac{\dot{\lambda}_t}{\lambda_t}
\psi_t^*\alpha_t,\]
or, with the help of Cartan's formula ${\mathcal L}_X=d\circ i_X+
i_X\circ d$ and with $\mu_t=\frac{d}{dt}(\log\lambda_t)\circ\psi_t^{-1}$,
\[ \psi_t^*\bigl( \dot{\alpha}_t+d(\alpha_t(X_t))+
i_{X_t}d\alpha_t\bigr) =\psi_t^*(\mu_t\alpha_t).\]
If we choose $X_t\in\xi_t$, this equation will be satisfied if
\begin{equation}
\label{eqn:gray1}
\dot{\alpha}_t+i_{X_t}d\alpha_t=\mu_t\alpha_t.
\end{equation}
Plugging in the Reeb vector field $R_{\alpha_t}$ gives
\begin{equation}
\label{eqn:gray2}
\dot{\alpha}_t(R_{\alpha_t})=\mu_t.
\end{equation}
So we can use (\ref{eqn:gray2}) to define $\mu_t$, and then the
non-degeneracy of $d\alpha_t|_{\xi_t}$ and the fact that $R_{\alpha_t}
\in\ker (\mu_t\alpha_t-\dot{\alpha}_t)$ allow us to find a unique
solution $X_t\in\xi_t$ of~(\ref{eqn:gray1}).
\end{proof}

\begin{rem}
\label{rem:gray}
{\rm
(1) Contact {\it forms} do {\it not} satisfy stability, that is, in general
one cannot find an isotopy $\psi_t$ such that $\psi_t^*\alpha_t
=\alpha_0$. For instance, consider the following family of contact forms on
$S^3\subset {\mathbb R}^4$:
\[ \alpha_t=(x_1\, dy_1-y_1\, dx_1)+(1+t)(x_2\, dy_2-y_2\, dx_2),\]
where $t\geq 0$ is a real parameter. The Reeb vector field of $\alpha_t$ is
\[ R_{\alpha_t}=(x_1\,\partial_{y_1}-y_1\,\partial_{x_1})+\frac{1}{1+t}
(x_2\,\partial_{y_2}-y_2\,\partial_{x_2}).\]
The flow of $R_{\alpha_0}$ defines the Hopf fibration, in particular all
orbits of $R_{\alpha_0}$ are closed. For $t\in {\mathbb R}^+\setminus
{\mathbb Q}$, on the other hand, $R_{\alpha_t}$ has only two periodic orbits.
So there can be no isotopy with $\psi_t^*\alpha_t=\alpha_0$, because
such a $\psi_t$ would also map $R_{\alpha_0}$ to~$R_{\alpha_t}$.

(2) Y.~Eliashberg~\cite{elia93} has shown that on the open manifold
${\mathbb R}^3$ there are likewise no non-trivial deformations of
contact structures, but on $S^1\times {\mathbb R}^2$ there does exist a
continuum of non-equivalent contact structures.

(3) For further applications of this theorem it is useful to observe
that at points $p\in M$ with $\dot{\alpha}_{t,p}$ identically zero in~$t$
we have $X_t(p)\equiv 0$, so such points remain stationary under
the isotopy~$\psi_t$.
}
\end{rem}
\subsection{Contact Hamiltonians}
A vector field $X$ on the contact manifold $(M,\xi =\ker\alpha )$ is called an
{\bf infinitesimal automorphism} of the contact structure
if the local flow of $X$ preserves~$\xi$ (The study of such automorphisms
was initiated by P.~Libermann, cf.~\cite{lima87}). By slight abuse of notation,
we denote this flow by $\psi_t$; if $M$ is not closed, $\psi_t$ (for
a fixed $t\neq 0$) will not in general be defined on all of~$M$. The
condition for $X$ to be an infinitesimal automorphism can be written as
$T\psi_t(\xi )=\xi$, which is equivalent to ${\mathcal L}_X\alpha =
\lambda \alpha$ for some function $\lambda\co M\rightarrow {\mathbb R}$
(notice that this condition is independent of the choice of $1$--form
$\alpha$ defining~$\xi$). The local flow of $X$ preserves $\alpha$ if
and only if ${\mathcal L}_X\alpha =0$.

\begin{thm}
With a fixed choice of contact form $\alpha$ there is a one-to-one
correspondence between infinitesimal automorphisms $X$ of $\xi=\ker\alpha$
and smooth functions $H\co M\rightarrow {\mathbb R}^+$. The correspondence
is given by
\begin{itemize}
\item $X\longmapsto H_X=\alpha (X)$;
\item $H\longmapsto X_H$, defined uniqely by $\alpha (X_H)=H$ and
$i_{X_H}d\alpha =dH(R_{\alpha})\alpha -dH$.
\end{itemize}
\end{thm}

The fact that $X_H$ is uniquely defined by the equations in the theorem
follows as in the preceding section from the fact that $d\alpha$ is
non-degenerate on $\xi$ and $R_{\alpha}\in\ker (dH(R_{\alpha})\alpha -dH)$.

\begin{proof}
Let $X$ be an infinitesimal automorphism of $\xi$. Set $H_X=\alpha (X)$
and write $dH_X+i_Xd\alpha =
{\mathcal L}_X\alpha =\lambda\alpha$ with $\lambda\co M\rightarrow
{\mathbb R}$. Applying this last equation to $R_{\alpha}$ yields
$dH_X (R_{\alpha})=\lambda$. So $X$ satisfies the
equations $\alpha (X)=H_X$ and $i_Xd\alpha =dH_X (R_{\alpha})\alpha -dH_X$.
This means that $X_{H_X}=X$.

Conversely, given $H\co M\rightarrow {\mathbb R}$ and with $X_H$ as defined
in the theorem, we have
\[ {\mathcal L}_{X_H}\alpha =i_{X_H}d\alpha +d(\alpha (X_H))=
dH(R_{\alpha})\alpha,\]
so $X_H$ is an infinitesimal automorphism of~$\xi$. Moreover, it
is immediate from the definitions that $H_{X_H}=\alpha (X_H)=H$.
\end{proof}

\begin{cor}
Let $(M,\xi =\ker\alpha )$ be a closed contact
manifold and $H_t\co M\rightarrow\R$, $t\in [0,1]$,
a smooth family of functions. Let $X_t=X_{H_t}$ be the corresponding
family of infinitesimal automorphisms of $\xi$ (defined via
the correspondence described in the preceding theorem). Then
the globally defined flow $\psi_t$ of the time-dependent vector
field $X_t$ is a contact isotopy of $(M,\xi )$, that is,
$\psi_t^*\alpha =\lambda_t\alpha$ for some smooth family of
functions $\lambda_t\co M\rightarrow\R^+$.
\end{cor}

\begin{proof}
With Lemma~\ref{lem:Lie} and the preceding proof we have
\[ \frac{d}{dt}\bigl( \psi_t^*\alpha \bigr) =
\psi_t^*\bigl( {\mathcal L}_{X_t}\alpha\bigr)=
\psi_t^*\bigl( dH_t(R_{\alpha})\alpha\bigr)
=\mu_t\psi_t^*\alpha\]
with $\mu_t=dH_t(R_{\alpha})\circ\psi_t$. Since $\psi_0=\mbox{\rm id}_M$
(whence $\psi_0^*\alpha =\alpha$) this implies that,
with
\[ \lambda_t=\exp \bigl( \int_0^t\mu_s\, ds\bigr) ,\]
we have $\psi_t^*\alpha =\lambda_t\alpha$. 
\end{proof}

This corollary will be used in Section~\ref{section:i-extension}
to prove various isotopy extension theorems from isotopies of
special submanifolds to isotopies of the ambient contact
manifold. In a similar vein, contact Hamiltonians can be used to show
that standard general position arguments from differential topology
continue to hold in the contact geometric setting. Another application
of contact Hamiltonians is a proof of the fact that the contactomorphism
group of a connected contact manifold acts transitively on
that manifold~\cite{boot69}. (See~\cite{bany97} for more on the general
structure of contactomorphism groups.)
\subsection{Darboux's theorem and neighbourhood theorems}
The flexibility of contact structures inherent in the Gray stability theorem
and the possibility to construct contact isotopies via contact Hamiltonians
results in a variety of theorems that can be summed up as saying that
there are no local invariants in contact geometry. Such theorems form the
theme of the present section.

In contrast with Riemannian geometry, for instance, where the local
structure coming from the curvature gives rise to a rich theory,
the interesting questions in contact geometry thus appear only at the global
level. However, it is actually that local flexibility that allows us
to prove strong global theorems, such as the existence of contact
structures on certain closed manifolds.
\subsubsection{Darboux's theorem}
\begin{thm}[Darboux's theorem]
\label{thm:darboux}
Let $\alpha$ be a contact form on the $(2n+1)$--dimensional
manifold $M$ and $p$ a point on~$M$.  Then there are coordinates
$x_1,\ldots ,x_n,y_1,\ldots ,y_n,z$ on a neighbourhood $U\subset M$ of $p$
such that
\[ \alpha|_U=dz+\sum_{j=1}^nx_j\, dy_j.\]
\end{thm}

\begin{proof}
We may assume without loss of generality that $M={\mathbb R}^{2n+1}$
and $p=0$ is the origin of ${\mathbb R}^{2n+1}$.
Choose linear coordinates $x_1,\ldots ,x_n,y_1,
\ldots y_n,z$ on ${\mathbb R}^{2n+1}$ such that
\[ \mbox{\rm on}\; T_0{\mathbb R}^{2n+1}:
\left\{ \begin{array}{l}
\alpha(\partial_z)=1,\;\; i_{\partial_z}d\alpha =0,\\
\partial_{x_j},\partial_{y_j}\in\ker\alpha\; (j=1,\ldots ,n),\;\;
d\alpha=\sum_{j=1}^ndx_j\wedge dy_j.
\end{array}\right. \]
This is simply a matter of linear algebra (the normal form theorem
for skew-symmetric forms on a vector space).

Now set $\alpha_0=dz+\sum_jx_j\, dy_j$ and consider the family of $1$--forms
\[ \alpha_t=(1-t)\alpha_0+t\alpha ,\; t\in [0,1],\]
on ${\mathbb R}^{2n+1}$. Our choice of coordinates ensures that
\[ \alpha_t=\alpha,\;\; d\alpha_t=d\alpha\;\; \mbox{\rm at the origin}.\]
Hence, on a sufficiently small neighbourhood of the
origin, $\alpha_t$ is a contact form for all $t\in [0,1]$.

We now want to use the Moser trick to find an isotopy $\psi_t$ of a
neighbourhood of the origin such that $\psi_t^*\alpha_t=\alpha_0$.
This aim seems to be in conflict with our earlier remark that
contact forms are not stable, but as we shall see presently, locally
this equation can always be solved.

Indeed, differentiating $\psi_t^*\alpha_t=\alpha_0$ (and assuming that
$\psi_t$ is the flow of some time-dependent vector field~$X_t$) we find
\[ \psi_t^*\bigl( \dot{\alpha}_t+{\mathcal L}_{X_t}\alpha_t\bigr) =0,\]
so $X_t$ needs to satisfy
\begin{equation}
\label{eqn:darboux1}
\dot{\alpha}_t+d(\alpha_t(X_t))+i_{X_t}d\alpha_t=0.
\end{equation}
Write $X_t=H_tR_{\alpha_t}+Y_t$ with $Y_t\in\ker\alpha_t$. Inserting
$R_{\alpha_t}$ in~(\ref{eqn:darboux1}) gives
\begin{equation}
\label{eqn:darboux2}
\dot{\alpha}_t(R_{\alpha_t})+dH_t(R_{\alpha_t})=0.
\end{equation}
On a neighbourhood of the origin, a smooth family of functions $H_t$ satisfying
(\ref{eqn:darboux2}) can always be found by integration, provided only that
this neighbourhood has been chosen so small that none of the $R_{\alpha_t}$
has any closed orbits there. Since $\dot{\alpha_t}$ is zero at the
origin, we may require that $H_t(0)=0$ and $dH_t|_0=0$ for
all $t\in [0,1]$. Once $H_t$ has been chosen, $Y_t$ is defined uniquely 
by~(\ref{eqn:darboux1}), i.e.\ by
\[ \dot{\alpha}_t+dH_t+ i_{Y_t}d\alpha_t=0.\]
Notice that with our assumptions on $H_t$ we have $X_t(0)=0$ for all~$t$.

Now define $\psi_t$ to be the local flow of $X_t$. This local flow
fixes the origin, so there it is defined for all~$t\in [0,1]$.
Since the domain of
definition in ${\mathbb R}\times M$ of a local flow on a manifold $M$
is always open (cf.~\cite[8.11]{brja73}),
we can infer\footnote{To be absolutely precise, one ought to work with
a family $\alpha_t$, $t\in \R$, where $\alpha_t\equiv\alpha_0$ for
$t\leq\varepsilon$ and $\alpha_t\equiv\alpha_1$ for $t\geq 1-
\varepsilon$, i.e.\ a {\it technical homotopy} in the sense
of~\cite{brja73}. Then $X_t$ will be defined for all $t\in\R$, and the
reasoning of \cite{brja73} can be applied.}
that $\psi_t$ is actually defined for
all $t\in [0,1]$ on a sufficiently small neighbourhood of the origin
in~${\mathbb R}^{2n+1}$. This concludes the proof of the theorem
(strictly speaking, the local coordinates in the statement of the
theorem are the coordinates $x_i\circ\psi_1^{-1}$ etc.).
\end{proof}

\begin{rem}
{\rm
The proof of this result given in~\cite{abklr94} is incomplete: It is not
possible, as is suggested there, to prove the Darboux theorem for
contact {\it forms} if one requires $X_t\in\ker\alpha_t$.
}
\end{rem}
\subsubsection{Isotropic submanifolds}
Let $L\subset (M,\xi =\ker\alpha )$ be an isotropic submanifold
in a contact manifold with cooriented contact structure. Write
$(TL)^{\perp}\subset\xi|_L$ for the subbundle of $\xi|_L$ that is
symplectically orthogonal to $TL$ with respect to the symplectic
bundle structure $d\alpha|_{\xi}$. The conformal class of this symplectic
bundle structure depends only on the contact structure~$\xi$,
not on the choice of contact form $\alpha$ defining~$\xi$: If
$\alpha$ is replaced by $f\alpha$ for some smooth function
$f\co M\rightarrow {\mathbb R}^+$, then $d(f\alpha)|_{\xi}=
f\, d\alpha|_{\xi}$. So the bundle $(TL)^{\perp}$ is determined by~$\xi$.

The fact that $L$ is isotropic implies $TL\subset (TL)^{\perp}$.
Following Weinstein~\cite{wein91}, we call the quotient bundle
$(TL)^{\perp}/TL$ with the conformal symplectic structure induced by
$d\alpha$ the {\bf conformal symplectic normal bundle} of $L$ in $M$ and
write
\[ \mbox{\rm CSN}(M,L)=(TL)^{\perp}/TL.\]
So the normal bundle $NL=(TM|_L)/TM$ of $L$ in $M$ can be split as
\[ NL\cong (TM|_L)/(\xi|_L)\oplus (\xi|_L)/(TL)^{\perp}
\oplus\mbox{\rm CSN}(M,L).\]
Observe that if $\dim M=2n+1$ and $\dim L=k\leq n$, then the ranks of
the three summands in this splitting are $1$, $k$ and $2(n-k)$,
respectively.
Our aim in this section is to show that a neighbourhood of $L$ in $M$ is
determined, up to contactomorphism, by the isomorphism type (as a conformal
symplectic bundle) of $\mbox{\rm CSN}(M,L)$.

The bundle $(TM|_L)/(\xi|_L)$ is a trivial line bundle because $\xi$ is
cooriented. The bundle $(\xi|_L)/(TL)^{\perp}$ can be identified with
the cotangent bundle $T^*L$ via the well-defined bundle
isomorphism
\[ \begin{array}{rrcl}
\Psi\co & (\xi|_L)/(TL)^{\perp} & \longrightarrow & T^*L\\
        &  Y                    &\longmapsto      &i_Yd\alpha.
\end{array} \]
($\Psi$ is obviously injective and well-defined by the definition of
$(TL)^{\perp}$, and the ranks of the two bundles are equal.)

Although $\Psi$ is well-defined on the quotient $(\xi|_L)/(TL)^{\perp}$,
to proceed further we need to choose an isotropic complement of
$(TL)^{\perp}$ in $\xi|_L$. Restricted to each fibre $\xi_p$, $p\in L$,
such an isotropic complement of $(T_pL)^{\perp}$ exists. There are two
ways to obtain a smooth bundle of such isotropic complements. The first
would be to carry over Arnold's corresponding discussion of Lagrangian
subbundles of symplectic bundles~\cite{arno67} to the isotropic case
in order to show that the space of isotropic complements of $U^{\perp}
\subset V$, where $U$ is an isotropic subspace in a symplectic vector
space~$V$, is convex. (This argument uses generating functions for isotropic
subspaces.) Then by a partition of unity argument the desired complement
can be constructed on the bundle level.

A slightly more pedestrian approach is to define this isotropic complement
with the help of a complex bundle structure $J$ on $\xi$ compatible
with~$d\alpha$ (cf.\ Remark~\ref{rem:symplectic}).
The condition $d\alpha (X,JX)>0$ for $0\neq X\in
\xi$ implies that $(T_pL)^{\perp}\cap J(T_pL)=\{ 0\}$ for all $p\in L$,
and so a dimension count shows that $J(TL)$ is indeed a complement of
$(TL)^{\perp}$ in~$\xi|_L$. (In a similar vein, $\mbox{\rm CSN}(M,L)$ can
be identified as a sub-bundle of~$\xi$, viz., the orthogonal complement
of $TL\oplus J(TL)\subset\xi$ with respect to the bundle metric
$d\alpha (.,J.)$ on~$\xi$.)

On the Whitney sum $TL\oplus T^*L$ (for any manifold~$L$) there is
a canonical symplectic bundle structure $\Omega_L$ defined by
\[ \Omega_{L,p}(X+Y,X'+Y')=Y(X')-Y'(X)\;\;\mbox{\rm for}\; X,X'\in T_pL;\,
Y,Y'\in T_p^*L.\]

\begin{lem}
\label{lem:i-nbhd}
The bundle map
\[ \mbox{\rm id}_{TL}\oplus\Psi\co (TL\oplus J(TL),d\alpha )\longrightarrow
(TL\oplus T^*L,\Omega_L)\]
is an isomorphism of symplectic vector bundles.
\end{lem}

\begin{proof}
We only need to check that $\mbox{\rm id}_{TL}\oplus\Psi$ is a symplectic
bundle map. Let $X,X',Y,Y'$ be as above. We can write $Y=J_pZ, Y'=J_pZ'$
with $Z,Z'\in T_pL$. It follows that
\[ d\alpha (Y,Y')=d\alpha (JZ,JZ')=d\alpha (Z,Z')=0,\]
since $L$ is an isotropic submanifold. For the same reason $d\alpha (X,X')
=0$. Hence
\begin{eqnarray*}
d\alpha (X+Y,X'+Y') & = & d\alpha (Y,X')-d\alpha (Y',X)\\
   & = & \Psi (Y)(X')-\Psi (Y')(X)\\
   & = & \Omega_L(X+\Psi (Y),X'+\Psi (Y')).
\qed
\end{eqnarray*}
\renewcommand{\qed}{}
\end{proof}

\begin{thm}
\label{thm:i-nbhd}
Let $(M_i,\xi_i)$, $i=0,1$, be contact manifolds with closed isotropic
submanifolds~$L_i$. Suppose there is an isomorphism of conformal
symplectic normal bundles $\Phi\co
\mbox{\rm CSN}(M_0,L_0)\rightarrow \mbox{\rm CSN}(M_1,L_1)$ that
covers a diffeomorphism $\phi\co L_1\rightarrow L_2$. Then $\phi$
extends to a contactomorphism $\psi\co{\mathcal N}(L_0)\rightarrow
{\mathcal N}(L_1)$ of suitable neighbourhoods ${\mathcal N}(L_i)$
of $L_i$ such that $T\psi|_{\mbox{\rm\scriptsize CSN}(M_0,L_0)}$ and
$\Phi$ are bundle homotopic (as symplectic bundle isomorphisms).
\end{thm}

\begin{cor}
Diffeomorphic (closed)
Legendrian submanifolds have contactomorphic neighbourhoods.
\end{cor}

\begin{proof}
If $L_i\subset M_i$ is Legendrian, then $\mbox{\rm CSN}(M_i,L_i)$ has rank~$0$,
so the conditions in the theorem, apart from the existence of a diffeomorphism
$\phi\co L_1\rightarrow L_2$, are void.
\end{proof}

\begin{exam}
\label{exam:Legendrian}
{\rm Let $S^1\subset (M^3,\xi )$ be a Legendrian knot in a contact
$3$--manifold. Then with a coordinate $\theta\in [0,2\pi ]$ along
$S^1$ and coordinates $x,y$ in slices transverse to~$S^1$, the
contact structure
\[ \cos\theta\, dx-\sin\theta\, dy =0\]
provides a model for a neighbourhood of~$S^1$.
}
\end{exam}

\begin{proof}[Proof of Theorem \ref{thm:i-nbhd}]
Choose contact forms $\alpha_i$ for $\xi_i$, $i=0,1$, scaled in such
a way that $\Phi$ is actually an isomorphism of symplectic vector bundles
with respect to the symplectic bundle structures on $\mbox{\rm CSN}(M_i,L_i)$
given by $d\alpha_i$. Here we think of $\mbox{\rm CSN}(M_i,L_i)$ as
a sub-bundle of $TM_i|_{L_i}$ (rather than as a quotient bundle).

We identify $(TM_i|_{L_i})/(\xi_i|_{L_i})$ with the trivial line bundle
spanned by the Reeb vector field $R_{\alpha_i}$. In total, this identifies
\[ NL_i=\langle R_{\alpha_i}\rangle\oplus J_i(TL_i)\oplus
\mbox{\rm CSN}(M_i,L_i)\]
as a sub-bundle of $TM_i|_{L_i}$.

Let $\Phi_R\co\langle R_{\alpha_0}\rangle\rightarrow
\langle R_{\alpha_1}\rangle$ be the obvious bundle isomorphism defined
by requiring that $R_{\alpha_0}(p)$ map to $R_{\alpha_1}(\phi (p))$.

Let $\Psi_i\co J_i(TL_i)\rightarrow T^*L_i$ be the isomorphism
defined by taking the interior product with~$d\alpha_i$. Notice that
\[ T\phi\oplus (\phi^*)^{-1}\co (TL_0\oplus T^*L_0,\Omega_{L_0})
\rightarrow (TL_1\oplus T^*L_1,\Omega_{L_1})\]
is an isomorphism of symplectic vector bundles. With Lemma~\ref{lem:i-nbhd}
it follows that
\[ T\phi\oplus\Psi_1^{-1}\circ (\phi^*)^{-1}\circ\Psi_0\co
(TL_0\oplus J_0(TL_0),d\alpha_0)\rightarrow
(TL_1\oplus J_1(TL_1),d\alpha_1) \]
is an isomorphism of symplectic vector bundles.

Now let
\[ \widetilde{\Phi}\co NL_0\longrightarrow NL_1\]
be the bundle isomorphism (covering $\phi$) defined by
\[ \widetilde{\Phi}=\Phi_R\oplus\Psi_1^{-1}\circ (\phi^*)^{-1}\circ\Psi_0
\oplus\Phi .\]
Let $\tau_i\co NL_i\rightarrow M_i$ be tubular maps, that is, the
$\tau$ (I suppress the index $i$ for better readability) are embeddings such
that $\tau|_L$ -- where $L$ is identified with the zero section of~$NL$ --
is the inclusion $L\subset M$, and $T\tau$ induces the identity on
$NL$ along~$L$ (with respect to the splittings $T(NL)|_L=TL\oplus NL
=TM|_L$).

Then $\tau_1\circ\widetilde{\Phi}\circ\tau_0^{-1}\co
{\mathcal N}(L_0)\rightarrow {\mathcal N}(L_1)$ is a diffeomorphism of
suitable neighbourhoods ${\mathcal N}(L_i)$ of $L_i$ that induces
the bundle map
\[ T\phi\oplus\widetilde{\Phi}\co TM_0|_{L_0}\longrightarrow
TM_1|_{L_1}.\]
By construction, this bundle map pulls $\alpha_1$ back to $\alpha_0$
and $d\alpha_1$ to $d\alpha_0$. Hence, $\alpha_0$ and
$(\tau_1\circ\widetilde{\Phi}\circ\tau_0^{-1})^*\alpha_1$
are contact forms on ${\mathcal N}(L_0)$ that coincide on $TM_0|_{L_0}$,
and so do their differentials.

Now consider the family of $1$--forms
\[ \beta_t=(1-t)\alpha_0+
t(\tau_1\circ\widetilde{\Phi}\circ\tau_0^{-1})^*\alpha_1,\;\; t\in [0,1].\]
On $TM_0|_{L_0}$ we have $\beta_t\equiv \alpha_0$ and $d\beta_t\equiv
d\alpha_0$. Since the contact condition $\alpha\wedge (d\alpha )^n\neq 0$
is an open condition, we may assume -- shrinking ${\mathcal N}(L_0)$
if necessary -- that $\beta_t$ is a contact form on ${\mathcal N}(L_0)$
for all $t\in [0,1]$. By the Gray stability theorem (Thm.~\ref{thm:gray})
and Remark~\ref{rem:gray}~(3)
following its proof, we find an isotopy $\psi_t$ of
${\mathcal N}(L_0)$, fixing~$L_0$, such that $\psi_t^*\beta_t=
\lambda_t\alpha_0$ for some smooth family of smooth
functions $\lambda_t\co {\mathcal N}(L_0)
\rightarrow {\mathbb R}^+$.

(Since ${\mathcal N}(L_0)$ is not a closed manifold, $\psi_t$ is
{\it a priori} only a local flow. But on $L_0$ it is stationary and hence
defined for all~$t$. As in the proof of the Darboux theorem
(Thm.~\ref{thm:darboux}) we conclude that $\psi_t$ is defined for all
$t\in [0,1]$ in a sufficiently small neighbourhood of~$L_0$, so
shrinking ${\mathcal N}(L_0)$ once again, if necessary, will ensure
that $\psi_t$ is a global flow on~${\mathcal N}(L_0)$.)

We conclude that $\psi =\tau_1\circ\widetilde{\Phi}
\circ\tau_0^{-1}\circ\psi_1$ is the desired contactomorphism.
\end{proof}

\begin{rem}
{\rm
With a little more care one can actually achieve $T\psi_1=\mbox{\rm id}$
on $TM_0|_{L_0}$, which implies in particular that
$T\psi|_{\mbox{\rm\scriptsize CSN}(M_0,L_0)}=\Phi$, cf.~\cite{wein91}.
(Remember that there
is a certain freedom in constructing an isotopy via the Moser trick
if the condition $X_t\in\xi_t$ is dropped.) The key point is
the generalised Poincar\'e lemma, cf.~\cite[p.~361]{lima87},
which allows to write
a closed differential form $\gamma$ given in a neighbourhood of the zero
section of a bundle and vanishing along that zero section as
an exact form $\gamma =d\eta$ with $\eta$ and its partial derivatives
with respect to all coordinates (in any chart) vanishing along the zero
section. This lemma is applied first to $\gamma =d(\beta_1-\beta_0)$,
in order to find (with the symplectic Moser trick) a diffeomorphism
$\sigma$ of a neighbourhood of $L_0\subset M_0$ with
$T\sigma =\mbox{\rm id}$ on $TM_0|_{L_0}$ and such that
$d\beta_0 =d(\sigma^*\beta_1)$. It is then applied once again
to $\gamma = \beta_0-\sigma^*\beta_1$.

(The proof of the symplectic neighbourhood theorem in~\cite{mcsa98}
appears to be incomplete in this respect.)
}
\end{rem}

\begin{exam}
{\rm
Let $M_0=M_1=\R^3$ with contact forms $\alpha_0=dz+x\, dy$
and $\alpha_1=dz+(x+y)\, dy$ and $L_0=L_1=0$ the origin in~$\R^3$.
Thus
\[ \mbox{\rm CSN}(M_0,L_0)=\mbox{\rm CSN}(M_1,L_1)=\mbox{\rm span}
\{\partial_x,\partial_y\}\subset T_0\R^3.\]
We take $\Phi=\mbox{\rm id}_{\mbox{\scriptsize\rm CSN}}$.

Set $\alpha_t=dz+(x+ty)\, dy$. The Moser trick with $X_t\in\ker
\alpha_t$ yields $X_t=-y\partial_x$, and hence $\psi_t(x,y,z)=
(x-ty,y,z)$. Then
\[ T\psi_1=\left(\begin{array}{ccc}1&-1&0\\
0&1&0\\0&0&1\end{array}\right),\]
which does not restrict to $\Phi$ on $\mbox{\rm CSN}$.

However, a different solution for $\psi_t^*\alpha_t=\alpha_0$ is
$\psi_t(x,y,z)=(x,y,z-ty^2/2)$, found by integrating $X_t=-y^2\partial_z
/2$ (a multiple of the Reeb vector field of~$\alpha_t$). Here we get
\[ T\psi_1=\left(\begin{array}{ccc}
1&0&0\\
0&1&0\\
0&-y&1\end{array}\right) ,\]
hence $T\psi_1|_{T_0\R^3}=\mbox{\rm id}$, so in particular
$T\psi_1|_{\mbox{\scriptsize\rm CSN}}=\Phi$.
}
\end{exam}
\subsubsection{Contact submanifolds}
Let $(M',\xi'=\ker\alpha')\subset (M,\xi=\ker\alpha )$ be a contact
submanifold, that is, $TM'\cap \xi|_{M'}=\xi'$. As before we write
$(\xi')^{\perp}\subset\xi|_{M'}$ for the symplectically orthogonal
complement of $\xi'$ in $\xi|_{M'}$. Since $M'$ is a contact
submanifold (so $\xi'$ is a symplectic sub-bundle of
$(\xi|_{M'},d\alpha )$), we have
\[ TM'\oplus (\xi')^{\perp}=TM|_{M'},\]
i.e.\ we can identify $(\xi')^{\perp}$ with the normal bundle $NM'$.
Moreover, $d\alpha$ induces a conformal symplectic structure on
$(\xi')^{\perp}$, so we call $(\xi')^{\perp}$ the {\bf conformal
symplectic normal bundle} of $M'$ in $M$ and write
\[ \mbox{\rm CSN}(M,M')=(\xi')^{\perp}.\]

\begin{thm}
\label{thm:c-nbhd}
Let $(M_i,\xi_i)$, $i=0,1$, be contact manifolds with compact
contact submanifolds $(M_i',\xi_i')$. Suppose there is an isomorphism
of conformal symplectic normal bundles $\Phi\co\mbox{\rm CSN}(M_0,M_0')
\rightarrow\mbox{\rm CSN}(M_1,M_1')$ that covers a contactomorphism
$\phi\co (M_0',\xi_0')\rightarrow (M_1',\xi_1')$. Then $\phi$ extends
to a contactomorphism $\psi$
of suitable neighbourhoods ${\mathcal N}(M_i')$ of
$M_i'$ such that $T\psi|_{\mbox{\rm\scriptsize CSN}(M_0,M_0')}$ and
$\Phi$ are bundle homotopic (as symplectic bundle isomorphisms)
up to a conformality.
\end{thm}

\begin{exam}
\label{exam:c-nbhd}
{\rm A particular instance of this theorem is the case of a transverse knot
in a contact manifold $(M,\xi )$, i.e.\ an embedding $S^1\hookrightarrow
(M,\xi )$ transverse to~$\xi$. Since the symplectic group
$\Sp (2n)$ of linear transformations of ${\mathbb R}^{2n}$ preserving
the standard symplectic structure $\omega_0=\sum_{i=1}^ndx_i\wedge dy_i$
is connected, there is only one conformal symplectic
${\mathbb R}^{2n}$--bundle over $S^1$ up to conformal equivalence.
A model for the neighbourhood of a transverse knot is given by
\[ \bigl( S^1\times{\mathbb R}^{2n},\xi =\ker \bigl( d\theta +\sum_{i=1}^n
(x_i\, dy_i-y_i\, dx_i)\bigr)\bigr) ,\]
where $\theta$ denotes the $S^1$--coordinate; the theorem says that
in suitable local coordinates the neighbourhood of any
transverse knot looks like this model.}
\end{exam}

\begin{proof}[Proof of Theorem~\ref{thm:c-nbhd}]
As in the proof of Theorem~\ref{thm:i-nbhd} it is sufficient to find
contact forms $\alpha_i$ on $M_i$ and a bundle map $TM_0|_{M_0'}
\rightarrow TM_1|_{M_1'}$, covering $\phi$ and inducing $\Phi$,
that pulls back $\alpha_1$ to $\alpha_0$ and $d\alpha_1$ to~$d\alpha_0$;
the proof then concludes as there with a stability argument.

For this we need to make a judicious choice of~$\alpha_i$. The essential
choice is made separately on each~$M_i$, so I suppress the subscript $i$
for the time being. Choose a contact form $\alpha'$ for $\xi'$ on~$M'$.
Write $R'$ for the Reeb vector field of~$\alpha'$. Given any contact form
$\alpha$ for $\xi$ on $M$ we may first scale it such that $\alpha (R')
\equiv 1$ along~$M'$. Then $\alpha|_{TM'}=\alpha'$, and hence
$d\alpha|_{TM'}=d\alpha'$. We now want to scale $\alpha$ further such that
its Reeb vector field $R$ coincides with $R'$ along~$M'$. To this end
it is sufficient to find a smooth function $f\co M\rightarrow
{\mathbb R}^+$ with $f|_{M'}\equiv 1$ and $i_{R'}d(f\alpha )\equiv 0$
on $TM|_{M'}$. This last equation becomes
\[ 0=i_{R'}d(f\alpha )=i_{R'}(df\wedge\alpha +f\, d\alpha )=
-df+i_{R'}d\alpha\;\;\mbox{\rm on}\;\; TM|_{M'}.\]
Since $i_{R'}d\alpha|_{TM'}=i_{R'}d\alpha'\equiv 0$, such an $f$
can be found.

The choices of $\alpha_0'$ and $\alpha_1'$ cannot be made independently
of each other; we may first choose $\alpha_1'$, say, and then define
$\alpha_0'=\phi^*\alpha_1'$. Then define $\alpha_0,\alpha_1$ as
described and scale $\Phi$ such that it is a symplectic bundle
isomorphism of
\[ ((\xi_0')^{\perp},d\alpha_0)\longrightarrow ((\xi_1')^{\perp},
d\alpha_1).\]
Then
\[ T\phi\oplus \Phi\co TM_0|_{M_0'}\longrightarrow TM_1|_{M_1'}\]
is the desired bundle map that pulls back $\alpha_1$ to $\alpha_0$
and $d\alpha_1$ to $d\alpha_0$.
\end{proof}

\begin{rem}
{\rm The condition that $R_i\equiv R_i'$ along $M'$ is necessary for
ensuring that $(T\phi\oplus\Phi ) (R_0)=R_1$, which guarantees (with
the other stated conditions) that $(T\phi\oplus\Phi )^*(d\alpha_1)=
d\alpha_0$. The condition $d\alpha_i|_{TM_i'}=d\alpha_i'$
and the described choice of $\Phi$ alone would only give
$(T\phi\oplus\Phi )^*(d\alpha_1|_{\xi_1})=d\alpha_0|_{\xi_0}$.
}
\end{rem}
\subsubsection{Hypersurfaces}
Let $S$ be an oriented hypersurface in a contact manifold $(M,\xi =
\ker\alpha )$ of dimension~$2n+1$. In a neighbourhood of $S$ in~$M$,
which we can identify with $S\times {\mathbb R}$ (and $S$ with
$S\times\{ 0\}$), the contact form $\alpha$ can be written as
\[ \alpha =\beta_r+u_r\, dr, \]
where $\beta_r$, $r\in {\mathbb R}$, is a smooth family of $1$--forms
on $S$ and $u_r\co S\rightarrow {\mathbb R}$ a smooth family of
functions. The contact condition $\alpha\wedge (d\alpha )^n\neq 0$ then
becomes
\begin{eqnarray}
\label{eqn:hyper1}
0\neq \alpha\wedge (d\alpha )^n & = & (\beta_r+u_r\, dr)\wedge
        (d\beta_r-\dot{\beta}_r\wedge dr+du_r\wedge dr)^n\nonumber\\
  & = & (-n\beta_r\wedge\dot{\beta}_r+n\beta_r\wedge du_r+u_r\, d\beta_r)
        \wedge (d\beta_r)^{n-1}\wedge dr,
\end{eqnarray}
where the dot denotes derivative with respect to~$r$.
The intersection $TS\cap (\xi|_S)$ determines a distribution (of
non-constant rank) of subspaces of~$TS$. If $\alpha$ is written as above,
this distribution is given by the kernel of~$\beta_0$, and hence, at a
given $p\in S$, defines either the full tangent space $T_pS$
(if~$\beta_{0,p}=0$) or a $1$--codimensional subspace both of $T_pS$
and $\xi_p$ (if~$\beta_{0,p}\neq 0$). In the former case, the
symplectically orthogonal complement $(T_pS\cap\xi_p)^{\perp}$ (with
respect to the conformal symplectic structure $d\alpha$ on~$\xi_p$)
is~$\{ 0\}$; in the latter case, $(T_pS\cap\xi_p)^{\perp}$ is a
$1$--dimensional subspace of $\xi_p$ contained in $T_pS\cap\xi_p$.

From that it is intuitively clear what one should mean by a
`singular $1$--dimensional foliation', and we make the following
somewhat provisional definition.

\begin{defn}
The {\bf characteristic foliation} $S_{\xi}$ of a hypersurface
$S$ in $(M,\xi )$ is the singular $1$--dimensional foliation
of $S$ defined by $(TS\cap\xi|_S)^{\perp}$.
\end{defn}

\begin{exam}
{\rm If $\dim M=3$ and $\dim S=2$, then $(T_pS\cap\xi_p)^{\perp}=
T_pS\cap\xi_p$ at the points $p\in S$ where $T_pS\cap\xi_p$ is
$1$--dimensional.
Figure~\ref{figure:charfolonS2} shows the characteristic foliation
of the unit $2$--sphere in $(\R^3,\xi_2)$, where $\xi_2$ denotes the
standard contact structure of Example~\ref{exam:standard2}: The only
singular points are $(0,0,\pm 1)$; away from these points the
characteristic foliation is spanned by
\[ (y-xz)\partial_x-(x+yz)\partial_y+(x^2+y^2)\partial_z.\]

\begin{figure}[ht]
        {\epsfysize=1.2truein\centerline{\epsfbox[0 719 130 841]{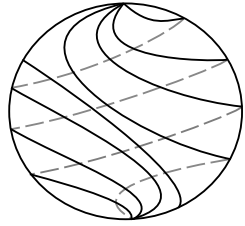}}}
        \caption{The characteristic foliation on  $S^2\subset (\R^3,\xi_2)$.}
        \label{figure:charfolonS2}
\end{figure}
}
\end{exam}

The following lemma helps to clarify the notion of singular $1$--dimensional
foliation.

\begin{lem}
\label{lem:singfol}
Let $\beta_0$ be the $1$--form induced on $S$ by a contact form
$\alpha$ defining~$\xi$, and let $\Omega$ be a volume form on~$S$.
Then $S_{\xi}$ is defined by the vector field $X$ satisfying
\[ i_X\Omega =\beta_0\wedge (d\beta_0)^{n-1}.\]
\end{lem}

\begin{proof}
First of all, we observe that $\beta_0\wedge (d\beta_0)^{n-1}\neq 0$
outside the zeros of~$\beta_0$: Arguing by contradiction, assume
$\beta_{0,p}\neq 0$ and $\beta_0\wedge (d\beta_0)^{n-1}|_p=0$ at some
$p\in S$. Then $(d\beta_0)^n|_p\neq 0$ by~(\ref{eqn:hyper1}). On
the codimension~$1$ subspace $\ker\beta_{0,p}$ of $T_pS$ the symplectic
form $d\beta_{0,p}$ has maximal rank $n-1$. It follows that
$\beta_0\wedge (d\beta_0)^{n-1}|_p\neq 0$ after all, a contradiction.

Next we want to show that $X\in\ker\beta_0$. We observe
\begin{equation}
\label{eqn:hyper2}
0=i_X(i_X\Omega )=\beta_0(X)(d\beta_0)^{n-1}-(n-1)\beta_0\wedge
i_Xd\beta_0\wedge (d\beta_0)^{n-2}.
\end{equation}

Taking the exterior product of this equation with $\beta_0$ we get
\[ \beta_0(X)\beta_0\wedge (d\beta_0)^{n-1}=0.\]
By our previous consideration this implies $\beta_0(X)=0$.

It remains to show that for $\beta_{0,p}\neq 0$ we have
\[ d\beta_0(X(p),v)=0\;\;\mbox{\rm for all}\; v\in T_pS\cap\xi_p.\]
For $n=1$ this is trivially satisfied, because in that case $v$ is
a multiple of~$X(p)$. I suppress the point~$p$ in the following calculation,
where we assume $n\geq 2$.
From~(\ref{eqn:hyper2}) and with $\beta_0(X)=0$ we have
\begin{equation}
\label{eqn:hyper3}
\beta_0\wedge i_Xd\beta_0\wedge (d\beta_0)^{n-2}=0.
\end{equation}
Taking the interior product with $v\in TS\cap\xi$ yields
\[ -d\beta_0(X,v)\beta_0\wedge (d\beta_0)^{n-2}+
(n-2)\beta_0\wedge i_Xd\beta_0\wedge i_vd\beta_0\wedge (d\beta_0)^{n-3}=0.\]
(Thanks to the coefficient $n-2$ the term $(d\beta_0)^{n-3}$ is not
a problem for $n=2$.)
Taking the exterior product of that last equation with $d\beta_0$,
and using~(\ref{eqn:hyper3}), we find
\[ d\beta_0(X,v)\beta_0\wedge (d\beta_0)^{n-1}=0,\]
and thus $d\beta_0(X,v)=0$.
\end{proof}

\begin{rem}
{\rm (1) We can now give a more formal definition of `singular $1$--dimensional
foliation' as an equivalence class of vector fields $[X]$, where $X$
is allowed to have zeros and $[X]=[X']$ if there is a nowhere zero
function on all of $S$ such that $X'=fX$. Notice that the non-integrability
of contact structures and the reasoning at the beginning
of the proof of the lemma imply that the zero set of $X$ does not
contain any open subsets of~$S$.

(2) If the contact structure $\xi$ is cooriented rather than just coorientable,
so that $\alpha$ is well-defined up to multiplication with a {\it positive}
function, then this lemma allows to give an orientation to the
characteristic foliation: Changing $\alpha$ to $\lambda\alpha$ with
$\lambda\co M\rightarrow {\mathbb R}^+$ will change $\beta_0\wedge
(d\beta_0)^{n-1}$ by a factor $\lambda^n$.
}
\end{rem}

We now restrict attention to surfaces in contact $3$--manifolds,
where the notion of characteristic foliation has proved to be
particularly useful.

The following theorem is due to E.~Giroux~\cite{giro91}.

\begin{thm}[Giroux]
\label{thm:Giroux}
Let $S_i$ be closed surfaces in
contact $3$--manifolds $(M_i,\xi_i)$,
$i=0,1$ (with $\xi_i$ coorientable),
and $\phi\co S_0\rightarrow S_1$ a diffeomorphism with
$\phi (S_{0,\xi_0})=S_{1,\xi_1}$ as oriented characteristic
foliations. Then there is a contactomorphism
$\psi\co {\mathcal N}(S_0)\rightarrow {\mathcal N}(S_1)$ of suitable
neighbourhoods ${\mathcal N}(S_i)$ of $S_i$ with $\psi (S_0)=S_1$
and such that $\psi|_{S_0}$ is isotopic to $\phi$ via an isotopy
preserving the characteristic foliation.
\end{thm}

\begin{proof}
By passing to a double cover, if necessary, we may assume that the
$S_i$ are orientable hypersurfaces. Let $\alpha_i$ be contact
forms defining~$\xi_i$. Extend $\phi$ to a diffeomorphism (still
denoted~$\phi$) of neighbourhoods of $S_i$ and consider the contact
forms $\alpha_0$ and $\phi^*\alpha_1$ on a neighbourhood of~$S_0$,
which we may identify with $S_0\times{\mathbb R}$.

By rescaling $\alpha_1$
we may assume that $\alpha_0$ and $\phi^*\alpha_1$ induce the same
form $\beta_0$ on $S_0\equiv S_0\times\{ 0\}$,
and hence also the same form $d\beta_0$.

Observe that the expression on the right-hand side of
equation~(\ref{eqn:hyper1}) is linear in $\dot{\beta}_r$ and~$u_r$. This
implies that convex linear combinations of solutions
of~(\ref{eqn:hyper1}) (for $n=1$) with the same $\beta_0$
(and $d\beta_0$) will
again be solutions of (\ref{eqn:hyper1}) for sufficiently small~$r$. This
reasoning applies to
\[ \alpha_t:= (1-t)\alpha_0+t\phi^*\alpha_1,\; t\in [0,1].\]
(I hope the reader will forgive the slight abuse of notation, with $\alpha_1$
denoting both a form on $M_1$ and its pull-back $\phi^*\alpha_1$ to~$M_0$.)
As is to be expected, we now use the Moser trick to find an isotopy
$\psi_t$ with $\psi_t^*\alpha_t=\lambda_t\alpha_0$, just
as in the proof of Gray stability (Theorem~\ref{thm:gray}). In particular,
we require as there that the vector field $X_t$ that we want to integrate
to the flow $\psi_t$ lie in the kernel of~$\alpha_t$.

On $TS_0$ we have $\dot{\alpha}_t\equiv 0$ (thanks to the assumption
that $\alpha_0$ and $\phi^*\alpha_1$ induce the same form $\beta_0$
on~$S_0$). In particular,
if $v$ is a vector in $S_{0,\xi_0}$, then
by equation~(\ref{eqn:gray1}) we have $d\alpha_t(X_t,v)=0$,
which implies that $X_t$ is a multiple of~$v$, hence tangent to~$S_{0,\xi_0}$.
This shows that the flow of $X_t$ preserves $S_0$ and
its characteristic foliation. More formally, we have
\[ {\mathcal L}_{X_t}\alpha_t=d(\alpha_t(X_t))+i_{X_t}d\alpha_t=
i_{X_t}d\alpha_t,\]
so with $v$ as above we have ${\mathcal L}_{X_t}\alpha_t(v)=0$, which shows
that ${\mathcal L}_{X_t}\alpha_t|_{TS_0}$ is a multiple of $\alpha_0|_{TS_0}
=\beta_0$. This implies that the (local) flow of $X_t$ changes $\beta_0$
by a conformal factor.

Since $S_0$ is closed, the local flow of $X_t$ restricted
to $S_0$ integrates up to $t=1$, and
so the same is true\footnote{Cf.\ the proof (and the footnote therein)
of Darboux's theorem (Thm.~\ref{thm:darboux}).}
in a neighbourhood of~$S_0$. Then $\psi =\phi\circ
\psi_1$ will be the desired diffeomorphism
${\mathcal N}(S_0)\rightarrow {\mathcal N}(S_1)$.
\end{proof}

As observed previously in the proof of Darboux's theorem for contact
{\it forms}, the Moser trick allows more flexibility if we drop the
condition $\alpha_t(X_t)=0$. We are now going to exploit this extra freedom
to strengthen Giroux's theorem slightly. This will be important later
on when we want to extend isotopies of hypersurfaces.

\begin{thm}
\label{thm:ch-nbhd}
Under the assumptions of the preceding theorem we can find
$\psi\co {\mathcal N}(S_0)\rightarrow {\mathcal N}(S_1)$ satisfying
the stronger condition that $\psi|_{S_0}=\phi$.
\end{thm}

\begin{proof}
We want to show that the isotopy $\psi_t$ of the preceding proof may be
assumed to fix $S_0$ pointwise. As there, we may assume
$\dot{\alpha}_t|_{TS_0}\equiv 0$.

If the condition that $X_t$ be tangent to
$\ker\alpha_t$ is dropped, the condition $X_t$ needs to satisfy so that
its flow will pull back $\alpha_t$ to $\lambda_t\alpha_0$ is
\begin{equation}
\label{eqn:hyper1'}
\dot{\alpha}_t+d(\alpha_t(X_t))+i_{X_t}d\alpha_t=\mu_t\alpha_t,
\end{equation}
where $\mu_t$ and $\lambda_t$ are related by $\mu_t=\frac{d}{dt}(\log
\lambda_t)\circ\psi_t^{-1}$, cf.\ the proof of the Gray stability theorem
(Theorem~\ref{thm:gray}).

Write $X_t=H_tR_t+Y_t$ with $R_t$ the Reeb vector field of $\alpha_t$
and $Y_t\in\xi_t=\ker\alpha_t$.
Then condition~(\ref{eqn:hyper1'}) translates into
\begin{equation}
\label{eqn:hyper3'}
\dot{\alpha}_t+dH_t+i_{Y_t}d\alpha_t=\mu_t\alpha_t.
\end{equation}
For given $H_t$ one determines $\mu_t$ from this equation by inserting
the Reeb vector field~$R_t$; the equation then
admits a unique solution $Y_t\in\ker\alpha_t$
because of the non-degeneracy of $d\alpha_t|_{\xi_t}$.

Our aim now is to ensure that $H_t\equiv 0$ on $S_0$ and $Y_t\equiv 0$
along~$S_0$. The latter we achieve by imposing the condition
\begin{equation}
\label{eqn:hyper4'}
\dot{\alpha}_t+dH_t=0\;\;\mbox{\rm along}\; S_0
\end{equation}
(which entails with (\ref{eqn:hyper3'}) that $\mu_t|_{S_0}\equiv 0$). 
The conditions $H_t\equiv 0$ on $S_0$ and (\ref{eqn:hyper4'}) can be
simultaneously satisfied thanks to $\dot{\alpha}_t|_{TS_0}\equiv 0$.

We can therefore find a smooth family of smooth
functions $H_t$ satisfying these conditions,
and then define $Y_t$
by~(\ref{eqn:hyper3'}). The flow of the vector field $X_t=
H_tR_t+Y_t$ then defines an isotopy $\psi_t$ that fixes
$S_0$ pointwise (and thus is defined for all $t\in [0,1]$ in a neighbourhood
of~$S_0$). Then $\psi=\phi\circ\psi_1$ will be the diffeomorphism we
wanted to construct.
\end{proof}
\subsubsection{Applications}
\label{section:nbhd-applications}
Perhaps the most important consequence of the neighbourhood theorems
proved above is that they allow us to perform differential topological
constructions such as surgery or similar cutting and pasting operations
in the presence of a contact structure, that is, these constructions can be
carried out on a contact manifold in such a way that the resulting manifold
again carries a contact structure.

One such construction that I shall explain in detail in
Section~\ref{section:3-manifolds} is the surgery of contact
$3$--manifolds along transverse knots, which enables us to construct
a contact structure on every closed, orientable $3$--manifold.

The concept of {\it characteristic foliation} on surfaces in contact
$3$--manifolds has proved seminal for the classification
of contact structures on $3$--manifolds, although it has recently
been superseded by the notion of {\it dividing curves}. It is clear
that Theorem~\ref{thm:Giroux} can be used to cut and paste contact
manifolds along hypersurfaces with the same characteristic foliation.
What actually makes this useful in dimension~$3$ is that there are ways
to manipulate the characteristic foliation of a surface by isotoping
that surface inside the contact $3$--manifold.

The most important result in that direction is the {\it Elimination
Lemma} proved by Giroux~\cite{giro91}; an improved version is due to
D.~Fuchs, see~\cite{elia93a}. This lemma says that under suitable
assumptions it is possible to cancel singular points of the characteristic
foliation in pairs by a $C^0$--small isotopy of the surface (specifically:
an elliptic and a hyperbolic point of the same sign -- the sign
being determined by the matching or non-matching of the orientation of
the surface $S$ and the contact structure $\xi$ at the singular point
of~$S_{\xi}$). For instance, Eliashberg~\cite{elia92} has shown that
if a contact $3$--manifold $(M,\xi )$ contains an embedded disc $D'$ such
that $D'_{\xi}$ has a limit cycle, then one can actually find a
so-called {\it overtwisted disc}: an embedded disc $D$ with boundary
$\partial D$ tangent to~$\xi$ (but $D$ transverse to $\xi$
along~$\partial D$, i.e.\ no singular points of $D_{\xi}$
on~$\partial D$) and with $D_{\xi}$ having exactly one singular
point (of elliptic type); see Section~\ref{subsection3.5}.

Moreover, in the generic situation it is possible, given surfaces
$S\subset (M,\xi )$ and $S'\subset (M',\xi')$ with $S_{\xi}$
{\it homeomorphic} to~$S'_{\xi'}$, to perturb one of the surfaces so
as to get {\it diffeomorphic} characteristic foliations.

Chapter~8 of~\cite{abklr94} contains a section on surfaces in contact
$3$--manifolds, and in particular a proof of the Elimination
Lemma. Further introductory reading on the matter
can be found in the lectures of J.~Etnyre~\cite{etny-lectures};
of the sources cited above I recommend~\cite{elia93a} as a starting
point.

In~\cite{giro91} Giroux initiated the study of {\it convex surfaces} in
contact $3$--manifolds. These are surfaces $S$ with an infinitesimal
automorphism $X$ of the contact structure $\xi$ with $X$
transverse to~$S$. For such surfaces, it turns out, much less
information than the characteristic foliation $S_{\xi}$ is needed
to determine $\xi$ in a neighbourhood of~$S$, viz., only the
so-called {\it dividing set} of~$S_{\xi}$. This notion lies at
the centre of most of the recent advances in the classification
of contact structures on $3$--manifolds \cite{giro00}, \cite{hond00},
\cite{hond00a}. A brief introduction to convex surface theory
can be found in~\cite{etny-lectures}.

\subsection{Isotopy extension theorems}
\label{section:i-extension}
In this section we show that the isotopy extension theorem of differential
topology -- an isotopy of a closed submanifold extends to an isotopy
of the ambient manifold -- remains valid for the various distinguished
submanifolds of contact manifolds. The neighbourhood theorems proved
above provide the key to the corresponding isotopy
extension theorems. For simplicity, I assume throughout
that the ambient contact manifold $M$ is closed; all isotopy extension theorems
remain valid if $M$ has non-empty boundary~$\partial M$, provided the isotopy
stays away from the boundary. In that case, the isotopy of $M$ found
by extension keeps a neighbourhood of $\partial M$ fixed. A further
convention throughout is that our ambient isotopies $\psi_t$ are
understood to start at $\psi_0=\mbox{\rm id}_M$.

\subsubsection{Isotropic submanifolds}

An embedding $j\co L\rightarrow (M,\xi =\ker\alpha )$ is called
{\bf isotropic} if $j(L)$ is an isotropic submanifold of
$(M,\xi )$, i.e.\ everywhere
tangent to the contact structure~$\xi$. Equivalently, one needs to
require $j^*\alpha\equiv 0$.

\begin{thm}
\label{thm:i-isotopy}
Let $j_t\co L\rightarrow (M,\xi =\ker\alpha )$, $t\in [0,1]$, be an
isotopy of isotropic embeddings of a closed manifold $L$ in a contact
manifold $(M,\xi )$.
Then there is a compactly supported
contact isotopy $\psi_t\co M\rightarrow M$
with $\psi_t(j_0(L))=j_t(L)$.
\end{thm}

\begin{proof}
Define a time-dependent vector field $X_t$ along $j_t(L)$ by
\[ X_t\circ j_t=\frac{d}{dt}j_t.\]
To simplify notation later on, we assume that $L$ is a submanifold of
$M$ and $j_0$ the inclusion $L\subset M$. Our aim is to find a (smooth)
family of compactly supported, smooth functions $\widetilde{H}_t\co
M\rightarrow {\mathbb R}$ whose Hamiltonian vector field $\widetilde{X}_t$
equals $X_t$ along $j_t(L)$. Recall that $\widetilde{X}_t$ is defined
in terms of $\widetilde{H}_t$ by
\[ \alpha (\widetilde{X}_t)=\widetilde{H}_t,\;\; i_{\widetilde{X}_t}d\alpha
=d\widetilde{H}_t (R_{\alpha})\alpha-d\widetilde{H}_t,\]
where, as usual, $R_{\alpha}$ denotes the Reeb vector field of~$\alpha$.

Hence, we need
\begin{equation}
\label{eqn:i-isotopy1}
\alpha (X_t)=\widetilde{H}_t,\;\; i_{X_t}d\alpha
=d\widetilde{H}_t (R_{\alpha})\alpha-d\widetilde{H}_t
\;\;\mbox{\rm along}\; j_t(L).
\end{equation}
Write $X_t=H_tR_{\alpha}+Y_t$ with $H_t\co j_t(L)\rightarrow
{\mathbb R}$ and $Y_t\in\ker\alpha$. To satisfy (\ref{eqn:i-isotopy1})
we need
\begin{equation}
\label{eqn:i-isotopy2}
\widetilde{H}_t=H_t\;\;\mbox{\rm along}\; j_t(L).
\end{equation}
This implies
\[ d\widetilde{H}_t(v)=dH_t(v)\;\;\mbox{\rm for}\; v\in T(j_t(L)).\]
Since $j_t$ is an isotopy of isotropic embeddings, we have
$T(j_t(L))\subset\ker\alpha$. So a prerequisite for
(\ref{eqn:i-isotopy1}) is that
\begin{equation}
\label{eqn:i-isotopy3}
d\alpha (X_t,v)=-dH_t(v)\;\;\mbox{\rm for}\; v\in T(j_t(L)).
\end{equation}
We have
\begin{eqnarray*}
d\alpha (X_t,v)+dH_t(v) & = & d\alpha (X_t,v)+d(\alpha (X_t))(v)\\
  & = & i_v(i_{X_t}d\alpha+d(i_{X_t}\alpha ))\\
  & = & i_v({\mathcal L}_{X_t}\alpha ),
\end{eqnarray*}
so equation~(\ref{eqn:i-isotopy3}) is equivalent to
\[ i_v({\mathcal L}_{X_t}\alpha )=0\;\;\mbox{\rm for}\; v\in T(j_t(L)).\]
But this is indeed tautologically satisfied: The fact that $j_t$ is an isotopy
of isotropic embeddings can be written as $j_t^*\alpha\equiv 0$; this
in turn implies $0=\frac{d}{dt}(j_t^*\alpha )=j_t^*
({\mathcal L}_{X_t}\alpha )$,
as desired.

This means that we can define $\widetilde{H}_t$ by prescribing
the value of $\widetilde{H}_t$ along $j_t(L)$ (with~(\ref{eqn:i-isotopy2}))
and the differential of $\widetilde{H}_t$ along $j_t(L)$
(with~(\ref{eqn:i-isotopy1})), where we are free to impose
$d\widetilde{H}_t(R_{\alpha})=0$, for instance. The calculation we just
performed shows that these two requirements are consistent with each
other. Any function satisfying these requirements along $j_t(L)$ can
be smoothed out to zero outside a tubular neighbourhood of $j_t(L)$,
and the Hamiltonian flow of this $\widetilde{H}_t$ will be the desired
contact isotopy extending~$j_t$.

One small technical point is to ensure that the resulting family of
functions $\widetilde{H}_t$ will be smooth in~$t$. To achieve this, we proceed
as follows. Set $\hat{M}=M\times [0,1]$ and
\[ \hat{L}=\bigcup_{q\in L,t\in [0,1]}(j_t(q),t),\]
so that $\hat{L}$ is a submanifold of~$\hat{M}$.
Let $g$ be an auxiliary Riemannian
metric on $M$ with respect to which $R_{\alpha}$
is orthogonal to $\ker\alpha$. Identify the normal bundle $N\hat{L}$ of
$\hat{L}$ in $\hat{M}$ with a sub-bundle of $T\hat{M}$ by requiring its fibre
at a point $(p,t)\in\hat{L}$ to be the $g$--orthogonal subspace of
$T_p(j_t(L))$ in $T_pM$. Let $\tau\co N\hat{L}\rightarrow \hat{M}$
be a tubular map.

Now define a smooth function $\hat{H}\co N\hat{L}\rightarrow {\mathbb R}$
as follows, where $(p,t)$ always denotes a point of $\hat{L}
\subset N\hat{L}$.
\begin{itemize}
\item $\hat{H}(p,t)=\alpha (X_t)$,
\item $d\hat{H}_{(p,t)}(R_{\alpha})=0$,
\item $d\hat{H}_{(p,t)}(v)=-d\alpha (X_t,v)$ for $v\in\ker\alpha_p
\subset T_pM\subset T_{(p,t)}\hat{M}$,
\item $\hat{H}$ is linear on the fibres of $N\hat{L}\rightarrow\hat{L}$.
\end{itemize}
Let $\chi\co\hat{M}\rightarrow [0,1]$ be a smooth function with
$\chi\equiv 0$ outside a small neighbourhood ${\mathcal N}_0
\subset\tau (N\hat{L})$
of $\hat{L}$ and $\chi\equiv 1$ in a smaller neighbourhood
${\mathcal N}_1\subset
{\mathcal N}_0$ of~$\hat{L}$. For $(p,t)\in\hat{M}$, set
\[ \widetilde{H}_t(p)=
\left\{ \begin{array}{ll}
\chi (p,t)\hat{H}(\tau^{-1}(p,t))\;\; & \mbox{\rm for}\;
                                           (p,t)\in\tau(N\hat{L})\\
0 & \mbox{\rm for}\; (p,t)\not\in\tau (N\hat{L}).
\end{array}\right.
\]
This is smooth in $p$ and $t$, and the Hamiltonian flow $\psi_t$ of
$\widetilde{H}_t$ (defined globally since $\widetilde{H}_t$ is
compactly supported) is the desired contact isotopy.
\end{proof}
\subsubsection{Contact submanifolds}
An embedding $j\co (M',\xi ')\rightarrow (M,\xi)$ is called a
{\bf contact embedding} if
\[ (j(M'), Tj(\xi ')) \]
is a contact submanifold
of $(M,\xi )$, i.e.\
\[ T(j(M))\cap\xi|_{j(M)}=Tj(\xi '). \]
If $\xi =\ker\alpha$, this can be reformulated as $\ker j^*\alpha =\xi '$.

\begin{thm}
\label{thm:c-isotopy}
Let $j_t\co (M',\xi ')\rightarrow (M,\xi )$, $t\in [0,1]$, be an isotopy
of contact embeddings of the closed contact manifold $(M',\xi ')$
in the contact manifold $(M,\xi )$. Then there is a compactly supported
contact isotopy $\psi_t\co M\rightarrow M$ with $\psi_t(j_0(M'))=j_t(M')$.
\end{thm}

\begin{proof}
In the proof of this theorem we follow a slightly different strategy
from the one in the isotropic case. Instead of directly finding an
extension of the Hamiltonian $H_t\co j_t(M')\rightarrow\R$, we first use the
neighbourhood theorem for contact submanifolds to extend $j_t$ to an
isotopy of contact embeddings of tubular neighbourhoods.

Again we assume that $M'$ is a submanifold of $M$ and $j_0$ the
inclusion $M'\subset M$. As earlier, $NM'$ denotes the normal bundle of $M'$
in~$M$. We also identify $M'$ with the zero section of $NM'$, and we use
the canonical identification
\[ T(NM')|_{M'}=TM'\oplus NM'.\]
By the usual isotopy extension theorem from differential topology we find
an isotopy
\[ \phi_t\co NM'\longrightarrow M\]
with $\phi_t|_{M'}=j_t$.

Choose contact forms $\alpha ,\alpha'$ defining $\xi$ and $\xi'$,
respectively. Define $\alpha_t=\phi_t^*\alpha$. Then $TM'\cap\ker
\alpha_t=\xi'$. Let $R'$ denote the Reeb vector field of~$\alpha'$.
Analogous to the proof of Theorem~\ref{thm:c-nbhd}, we first find
a smooth family of smooth functions $g_t\co M'\rightarrow\R^+$ such that
$g_t\alpha_t|_{TM'}=\alpha'$, and then a family $f_t\co NM'\rightarrow
\R^+$ with $f_t|_{M'}\equiv 1$ and
\[ df_t=i_{R'}d(g_t\alpha_t)\;\; \mbox{\rm on}\;\; T(NM')|_{M'}.\]
Then $\beta_t=f_tg_t\alpha_t$ is a family of contact forms on $NM'$
representing the contact structure $\ker (\phi_t^*\alpha )$ and with
the properties
\begin{eqnarray*}
\beta_t|_{TM'} & = & \alpha',\\
d\beta_t|_{TM'} & = & d\alpha',\\
\ker (d\beta_t) & = & R'\;\;\mbox{\rm along}\;\; M'.
\end{eqnarray*}
The family $(NM',d\beta_t)$ of symplectic vector bundles may be thought of as
a symplectic vector bundle over $M'\times [0,1]$, which is necessarily
isomorphic to a bundle pulled back from $M'\times\{ 0\}$
(cf.~\cite[Cor.~3.4.4]{huse94}). In other words,
there is a smooth family of symplectic bundle isomorphisms
\[ \Phi_t\co (NM',d\beta_0)\longrightarrow (NM',d\beta_t).\]
Then
\[ \mbox{\rm id}_{TM'}\oplus\Phi_t\co T(NM')|_{M'}
\longrightarrow T(NM')|_{M'}\]
is a bundle map that pulls back $\beta_t$ to $\beta_0$ and
$d\beta_t$ to~$d\beta_0$.

By the now familiar stability argument we find a smooth family of
embeddings
\[ \varphi_t\co {\mathcal N}(M')\longrightarrow NM'\]
for some neighbourhood ${\mathcal N}(M')$ of the
zero section $M'$ in $NM'$ with
$\varphi_0=\mbox{\rm inclusion}$, $\varphi_t|_{M'}=\mbox{\rm id}_{M'}$
and $\varphi_t^*\beta_t=\lambda_t\beta_0$, where $\lambda_t\co
{\mathcal N}(M')\rightarrow\R^+$. This means that
\[ \phi_t\circ\varphi_t\co {\mathcal N}(M')\longrightarrow M\]
is a smooth family of contact embeddings of $({\mathcal N}(M'),\ker\beta_0)$
in $(M,\xi )$.

Define a time-dependent vector field $X_t$ along $\phi_t\circ\varphi_t
({\mathcal N}(M'))$ by
\[ X_t\circ \phi_t\circ\varphi_t=\frac{d}{dt}(\phi_t\circ\varphi_t).\]
This $X_t$ is clearly an infinitesimal automorphism of~$\xi$: By
differentiating the equation $\varphi_t^*\phi_t^*\alpha=\mu_t
\phi_0^*\alpha$ (where $\mu_t\co {\mathcal N}(M')\rightarrow\R^+$)
with respect to $t$ we get
\[ \varphi_t^*\phi_t^*({\mathcal L}_{X_t}\alpha )=\dot{\mu}_t\phi_0^*\alpha
=\frac{\dot{\mu}_t}{\mu_t} \varphi_t^*\phi_t^*\alpha ,\]
so ${\mathcal L}_{X_t}\alpha$ is a multiple of~$\alpha$ (since
$\phi_t\circ\varphi_t$ is a diffeomorphism onto its image).

By the theory of contact Hamiltonians, $X_t$ is the Hamiltonian
vector field of a Hamiltonian function $\hat{H}_t$ defined on
$\phi_t\circ\varphi_t ({\mathcal N}(M'))$. Cut off this
function with a bump function so as to obtain $H_t\co M\rightarrow {\mathbb R}$
with $H_t\equiv \hat{H}_t$ near $\phi_t\circ\varphi_t(M')$ and
$H_t\equiv 0$ outside a slightly larger neighbourhood of
$\phi_t\circ\varphi_t(M')$.
Then the Hamiltonian flow $\psi_t$ of $H_t$ satisfies
our requirements.
\end{proof}
\subsubsection{Surfaces in $3$--manifolds}
\begin{thm}
Let $j_t\co S\rightarrow (M,\xi =\ker\alpha )$, $t\in [0,1]$, be an
isotopy of embeddings of a closed surface $S$
in a $3$--dimensional contact manifold $(M,\xi )$.
If all $j_t$ induce
the same characteristic foliation on~$S$, then there is a compactly supported
isotopy $\psi_t\co M\rightarrow M$
with $\psi_t(j_0(S))= j_t(S)$.
\end{thm}

\begin{proof}
Extend $j_t$ to a smooth family of embeddings $\phi_t\co S\times{\mathbb R}
\rightarrow M$, and identify $S$ with $S\times\{ 0\}$. The assumptions
say that all $\phi_t^*\alpha$ induce the same characteristic foliation
on~$S$. By the proof of Theorem~\ref{thm:ch-nbhd} and in analogy with
the proof of Theorem~\ref{thm:c-isotopy} we find a smooth
family of embeddings
\[ \varphi_t\co S\times (-\varepsilon ,\varepsilon )\longrightarrow
S\times {\mathbb R}\]
for some $\varepsilon >0$ with $\varphi_0=\mbox{\rm inclusion}$,
$\varphi_t|_{S\times\{ 0\} }=\mbox{\rm id}_S$ and
 $\varphi_t^*\phi_t^*\alpha
=\lambda_t\phi_0^*\alpha$, where $\lambda_t\co S\times
(-\varepsilon ,\varepsilon )\rightarrow {\mathbb R}^+$. In other words,
$\phi_t\circ\varphi_t$ is a smooth family of contact embeddings of
$(S\times (-\varepsilon ,\varepsilon ),\ker\phi_0^*\alpha)$ in $(M,\xi )$.

The proof now concludes exactly as the proof of
Theorem~\ref{thm:c-isotopy}.
\end{proof}
\subsection{Approximation theorems}
A further manifestation of the (local) flexibility of contact structures
is the fact that a given submanifold can, under fairly weak (and
usually obvious) topological conditions, be approximated (typically
$C^0$--closely) by a contact submanifold or an isotropic submanifold,
respectively. The most general results in this direction are best
phrased in M.~Gromov's language of {\it $h$-principles}. For a recent
text on $h$-principles that puts particular emphasis on
symplectic and contact geometry see~\cite{elmi02}; a brief and perhaps more
gentle introduction to $h$-principles can be found in~\cite{geig02}.

In the present section I restrict attention to the $3$--dimensional
situation, where the relevant approximation theorems can be proved
by elementary geometric {\it ad hoc} techniques.

\begin{thm}
\label{thm:approximation}
Let $\gamma\co S^1\rightarrow (M,\xi )$ be a knot, i.e.\
an embedding of $S^1$, in a contact $3$--manifold. Then $\gamma$ can be
$C^0$--approximated by a Legendrian knot isotopic
to~$\gamma$. Alternatively, it can be
$C^0$--approximated by a positively as well as a negatively
transverse knot.
\end{thm}

In order to prove this theorem, we first consider embeddings
$\gamma\co (a,b)\rightarrow (\R^3,\xi )$ of an open interval
in $\R^3$ with its standard contact structure $\xi =\ker\alpha$, where
$\alpha =dz+x\, dy$. Write $\gamma (t)=(x(t),y(t),z(t))$. Then
\[ \alpha (\dot{\gamma})=\dot{z}+x\dot{y},\]
so the condition for a Legendrian curve reads $\dot{z}+x\dot{y}
\equiv 0$; for a positively (resp.\ negatively) transverse curve,
$\dot{z}+x\dot{y}>0$ (resp.~$<0$).

There are two ways to visualise such curves. The first is via its
{\bf front projection}
\[ \gamma_F(t)=(y(t),z(t)),\]
the second via its {\bf Lagrangian projection}
\[ \gamma_L(t)=(x(t),y(t)).\]
\subsubsection{Legendrian knots}
If $\gamma (t)=(x(t),y(t),z(t))$ is a Legendrian curve in~$\R^3$, then
$\dot{y}=0$ implies $\dot{z}=0$, so there the front projection has a singular
point. Indeed, the curve $t\mapsto (t,0,0)$ is an example of a Legendrian
curve whose front projection is a single point. We call a Legendrian
curve {\it generic} if $\dot{y}=0$ only holds at isolated points
(which we call {\bf cusp points}), and there $\ddot{y}\neq 0$.

\begin{figure}[ht]
        {\epsfysize=3cm \centerline{\epsfbox{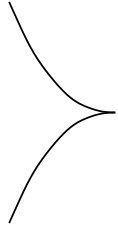}}}
        \caption{The cusp of a front projection.}
        \label{figure:cusp}
\end{figure}

\begin{lem}
Let $\gamma\co (a,b)\rightarrow (\R^3,\xi )$ be a Legendrian
immersion. Then its front projection $\gamma_F(t)=(y(t),z(t))$ does
not have any vertical tangencies.
Away from the cusp points, $\gamma$ is recovered from its
front projection via
\[ x(t)=-\frac{\dot{z}(t)}{\dot{y}(t)}=-\frac{dz}{dy},\]
i.e.\ $x(t)$ is the negative slope of the front projection.
The curve $\gamma$ is embedded if and only
if $\gamma_F$ has only transverse self-intersections.

By a $C^{\infty}$--small perturbation of $\gamma$ we can obtain
a generic Legendrian curve $\tilde{\gamma}$
isotopic to~$\gamma$; by a $C^2$--small perturbation
we may achieve that the front
projection has only semi-cubical cusp singularities, i.e.\
around a cusp point at $t=0$ the curve $\tilde{\gamma}$
looks like
\[ \tilde{\gamma} (t)=(t+a,\lambda t^2+b, -\lambda (2t^3/3+at^2)+c)\]
with $\lambda\neq 0$, see Figure~\ref{figure:cusp}.

Any regular curve in the $(y,z)$--plane with semi-cubical cusps and
no vertical tangencies can be lifted to a unique
Legendrian curve in~$\R^3$.
\end{lem}

\begin{proof}
The Legendrian condition is $\dot{z}+x\dot{y}=0$. Hence $\dot{y}=0$
forces $\dot{z}=0$, so $\gamma_F$ cannot have any vertical tangencies.

Away from the cusp points, the Legendrian condition tells us how
to recover $x$ as the negative slope of the front projection. (By
continuity, the equation $x=\frac{dz}{dy}$ also makes sense at
generic cusp points.) In particular,
a self-intersecting front projection lifts to a non-intersecting curve
if and only if the slopes at the intersection point are different, i.e.\
if and only if the intersection is transverse.

That $\gamma$ can be approximated in the $C^{\infty}$--topology
by a generic immersion $\tilde{\gamma}$ follows from
the usual transversality theorem (in its most simple form,
viz., applied to the function $y(t)$; the function $x(t)$ may be
left unchanged, and the new $z(t)$ is then found by integrating
the new $-x\dot{y}$).

At a cusp point of $\tilde{\gamma}$ we have $\dot{y}=\dot{z}=0$.
Since $\tilde{\gamma}$ is an immersion, this forces
$\dot{x}\neq 0$, so $\tilde{\gamma}$ can be parametrised around a cusp
point by the $x$--coordinate, i.e.\ we may choose
the curve parameter $t$ such that the cusp lies at $t=0$
and $x(t)=t+a$. Since $\ddot{y}(0)\neq 0$ by the genericity
condition, we can write $y(t)=t^2g(t)+y(0)$ with a smooth function
$g(t)$ satisfying $g(0)\neq 0$ (This is proved like the
`Morse lemma' in Appendix~2 of~\cite{kosi93}). A $C^0$--approximation
of $g(t)$ by a function $h(t)$ with $h(t)\equiv g(0)$ for $t$ near
zero and $h(t)\equiv g(t)$ for $|t|$ greater than some small
$\varepsilon >0$ yields a $C^2$--approximation of $y(t)$ with
the desired form around the cusp point.
\end{proof}

\begin{figure}[ht]
        {\epsfysize=5cm \centerline{\epsfbox{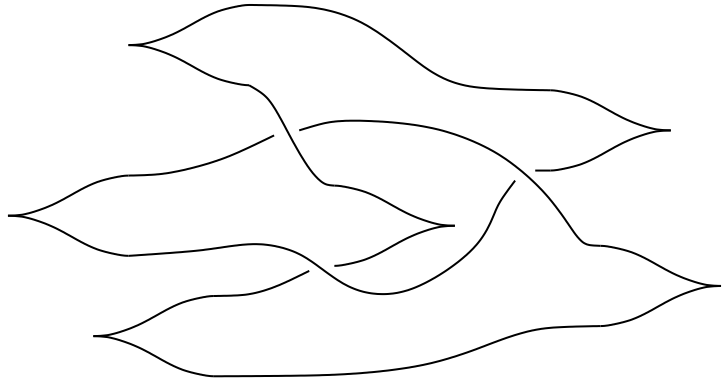}}}
        \caption{Front projection of a Legendrian trefoil knot.}
        \label{figure:trefoil}
\end{figure}

\begin{exam}
{\rm
Figure~\ref{figure:trefoil} shows the front projection of a Legendrian
trefoil knot.
}
\end{exam}

\begin{proof}[Proof of Theorem~\ref{thm:approximation} - Legendrian case]
First of all, we consider a curve $\gamma$ in standard~$\R^3$. In order to
find a $C^0$--close approximation of $\gamma$, we simply need to
choose a $C^0$--close approximation of its front projection $\gamma_F$
by a regular curve without vertical tangencies and with isolated
cusps (we call such a curve a {\it front})
in such a way, that the slope of the front at time~$t$
is close to $-x(t)$ (see Figure~\ref{figure:approxleg}).
Then the Legendrian lift of this front is the desired
$C^0$--approximation of~$\gamma$.

\begin{figure}[h]
\centerline{\relabelbox\small
\epsfxsize 7cm \epsfbox{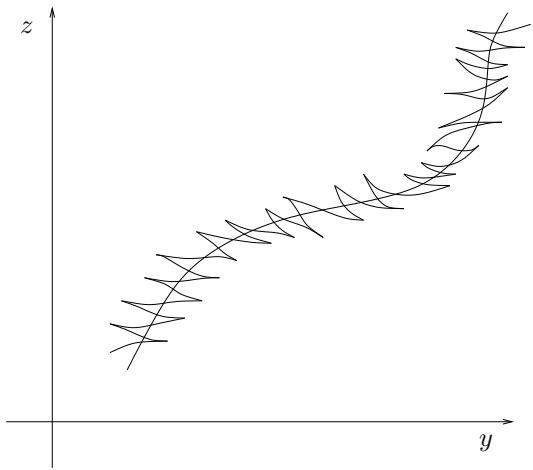}
\extralabel <-.7cm,.3cm> {$y$}
\extralabel <-6.8cm,5.8cm> {$z$}
\endrelabelbox}
\caption{Legendrian $C^0$--approximation via front projection.}
\label{figure:approxleg}
\end{figure}

If $\gamma$ is defined on an interval $(a,b)$ and is already Legendrian
near its endpoints, then the approximation of $\gamma_F$ may be
assumed to coincide with $\gamma_F$ near the endpoints, so that
the Legendrian lift coincides with $\gamma$ near the endpoints.

Hence, given a knot in an arbitrary contact $3$--manifold, we
can cut it (by the Lebesgue lemma) into little pieces that lie
in Darboux charts. There we can use the preceding recipe to
find a Legendrian approximation. Since, as just observed, one
can find such approximations on intervals with given boundary
condition, this procedure yields a Legendrian approximation of
the full knot.

Locally (i.e.\ in~$\R^3$) the described procedure does not
introduce any self-intersections in the approximating curve,
provided we approximate $\gamma_F$ by a front with only transverse
self-intersections. Since the original knot was embedded, the
same will then be true for its Legendrian $C^0$--approximation.
\end{proof}

The same result may be derived using the Lagrangian projection:

\begin{lem}
Let $\gamma\co (a,b)\rightarrow (\R^3,\xi )$ be a Legendrian immersion.
Then its Lagrangian projection $\gamma_L(t)=(x(t),y(t))$ is also
an immersed curve. The curve $\gamma$ is recovered from $\gamma_L$ via
\[ z(t_1)=z(t_0)-\int_{t_0}^{t_1}x\, dy.\]
A Legendrian immersion $\gamma\co S^1\rightarrow (\R^3,\xi )$ has
a Lagrangian projection that encloses zero area.
Moreover, $\gamma$ is embedded if and only if every loop in $\gamma_L$
(except, in the closed case, the full loop~$\gamma_L$)
encloses a non-zero oriented area.

Any immersed curve in the $(x,y)$--plane is the Lagrangian projection of
a Legendrian curve in~$\R^3$, unique up to translation in the
$z$--direction.
\end{lem}

\begin{proof}
The Legendrian condition $\dot{z}+x\dot{y}$ implies that if $\dot{y}=0$
then $\dot{z}=0$, and hence, since $\gamma$ is an immersion,
$\dot{x}\neq 0$. So $\gamma_L$ is an immersion.

The formula for $z$ follows by integrating the Legendrian condition.
For a closed curve $\gamma_L$ in the $(x,y)$--plane, the integral
$\oint_{\gamma_L}x\, dy$ computes the oriented area enclosed by~$\gamma_L$.
From that all the other statements follow.
\end{proof}

\begin{figure}[h]
\centerline{\relabelbox\small
\epsfxsize 10cm \epsfbox{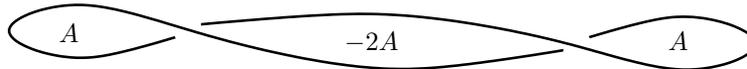}
\extralabel <-1.2cm,.3cm> {$A$}
\extralabel <-9.3cm,.4cm> {$A$}
\extralabel <-5.5cm,.3cm> {$-2A$}
\endrelabelbox}
\caption{Lagrangian projection of a Legendrian unknot.}
\label{figure:unknot}
\end{figure}

\begin{exam}
{\rm
Figure~\ref{figure:unknot} shows the Lagrangian projection of
a Legendrian unknot.
}
\end{exam}

\begin{figure}[h]
\centerline{\relabelbox\small
\epsfxsize 7cm \epsfbox{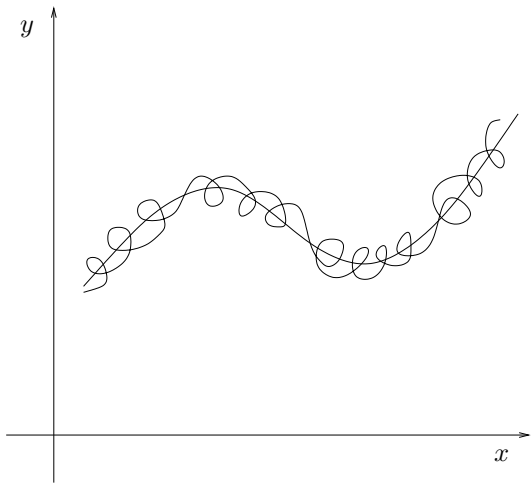}
\extralabel <-.5cm,.3cm> {$x$}
\extralabel <-6.8cm,6.0cm> {$y$}
\endrelabelbox}
\caption{Legendrian $C^0$--approximation via Lagrangian projection.}
\label{figure:approxleg2}
\end{figure}

\begin{proof}[Alternative proof of Theorem~\ref{thm:approximation} --
Legendrian case]
Again we consider a curve $\gamma$ in standard~$\R^3$ defined on
an interval. The generalisation to arbitrary contact manifolds
and closed curves is achieved as in the proof using front projections.

In order to find a $C^0$--approximation of $\gamma$ by a
Legendrian curve, one only has to approximate its Lagrangian projection
$\gamma_L$ by an immersed curve whose `area integral'
\[ z(t_0)-\int_{t_0}^t x\, dy\]
lies as close to the original
$z(t)$ as one wishes. This can be achieved by using
small loops oriented positively or negatively
(see Figure~\ref{figure:approxleg2}). If $\gamma_L$ has self-intersections,
this approximating curve can be chosen in such a way that along loops
properly contained in that curve the area integral is non-zero,
so that again we do not introduce any self-intersections in
the Legendrian approximation of~$\gamma$.
\end{proof}

\subsubsection{Transverse knots}
\label{section:transverse}
The quickest proof of the transverse case of
Theorem~\ref{thm:approximation} is via the Legendrian case. However, it
is perfectly feasible to give a direct proof along the lines
of the preceding discussion, i.e.\ using the front or the
Lagrangian projection. Since this picture is useful elsewhere, I
include a brief discussion, restricting attention to the front
projection.

Thus, let $\gamma (t)=(x(t),y(t),z(t))$ be a curve in~$\R^3$.
The condition for $\gamma$ to be positively transverse to the
standard contact structure $\xi =\ker (dz+x\, dy)$ is that
$\dot{z}+x\dot{y}>0$. Hence,
\[ \left\{\begin{array}{l}
\mbox{\rm if}\; \dot{y}=0,\;\mbox{\rm then}\; \dot{z}>0,\\
\mbox{\rm if}\; \dot{y}>0,\;\mbox{\rm then}\; x>-\dot{z}/\dot{y},\\
\mbox{\rm if}\; \dot{y}<0,\;\mbox{\rm then}\; x<-\dot{z}/\dot{y}.\\
\end{array}\right. \]

\begin{figure}[h]
\centerline{\relabelbox\small
\epsfxsize 6cm \epsfbox{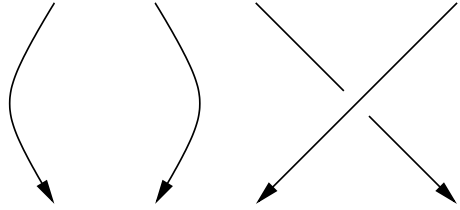}
\endrelabelbox}
\caption{Impossible front projections of positively transverse curve.}
\label{figure:transverse-not}
\end{figure}

The first statement says that there are no vertical tangencies
oriented downwards in the front projection. The second statement says
in particular that for $\dot{y}>0$ and $\dot{z}<0$ we have $x>0$;
the third, that for $\dot{y}<0$ and $\dot{z}<0$ we have $x<0$.
This implies that the situations shown in
Figure~\ref{figure:transverse-not} are not possible in the
front projection of a positively transverse curve. I leave it to the reader
to check that all other oriented crossings are possible in the front
projection of a positively transverse curve, and that any curve in
the $(y,z)$--plane without the forbidden crossing or downward
vertical tangencies admits a lift to a positively transverse
curve.

\begin{figure}[h]
\centerline{\relabelbox\small
\epsfxsize 6cm \epsfbox{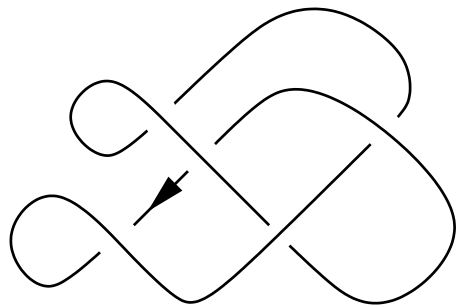}
\endrelabelbox}
\caption{Front projection of a positively transverse trefoil knot.}
\label{figure:trefoil2}
\end{figure}

\begin{exam}
{\rm
Figure~\ref{figure:trefoil2} shows the front projection of a positively
transverse trefoil knot.
}
\end{exam}

\begin{proof}[Proof of Theorem~\ref{thm:approximation} -- transverse case]
By the Legendrian case of this theorem, the given knot $\gamma$ can be
$C^0$--approximated by a Legendrian knot~$\gamma_1$. By
Example~\ref{exam:Legendrian}, a neighbourhood of $\gamma_1$ in
$(M,\xi )$ looks like a solid torus $S^1\times D^2$ with
contact structure $\cos\theta\, dx-\sin\theta\, dy=0$, where
$\gamma_1\equiv S^1\times\{ 0\}$. Then the curve
\[ \gamma_2(t)=(\theta =t,x=\delta\sin t,y=\delta\cos t),\;\; t\in [0,
2\pi ],\]
is a positively (resp.\ negatively)
transverse knot if $\delta >0$ (resp.~$<0$). By choosing $|\delta |$ small
we obtain as good a $C^0$--approximation of $\gamma_1$ and
hence of $\gamma$ as we wish.
\end{proof}


%% file: section3.tex
\section{Contact structures on $3$--manifolds}
\label{section:3-manifolds}

Here is the main theorem proved in this section:

\begin{thm}[Lutz-Martinet]
\label{thm:lutz-martinet}
Every closed, orientable $3$--manifold admits a contact structure in each
homotopy class of tangent $2$--plane fields.
\end{thm}

In Section~\ref{subsection3.1} I present what is essentially
J.~Martinet's~\cite{mart71}
proof of the existence of a contact structure on every $3$--manifold.
This construction 
is based on a surgery description of $3$--manifolds due to R.~Lickorish
and A.~Wallace. For the key step, showing how to extend over a solid torus
certain contact forms defined near the boundary of that torus (which then
makes it possible to perform Dehn surgeries), we use an approach due to
W.~Thurston and H.~Winkelnkemper;
this allows to simplify Martinet's proof slightly.

In Section~\ref{subsection3.2}
we show that every orientable $3$--manifold is parallelisable
and then build on this to classify (co-)oriented tangent $2$--plane
fields up to homotopy.

In Section~\ref{subsection3.3}
we study the so-called Lutz twist, a topologically trivial
Dehn surgery on a contact manifold $(M,{\xi})$ which yields a contact
structure ${\xi}'$ on $M$ that is not homotopic (as $2$--plane field)
to~$\xi$. We then complete the proof
of the main theorem stated above. These results are
contained in R.~Lutz's thesis~\cite{lutz71} (which, I have to admit,
I've never seen). Of Lutz's published work, \cite{lutz70} only deals
with the $3$--sphere (and is only an announcement);
\cite{lutz77} also deals with a more restricted problem.
I learned the key steps
of the construction from
an exposition given in V.~Ginzburg's 
thesis~\cite{ginz90}.
I have added proofs of many topological details that do not seem
to have appeared in a readily accessible source before.

In Section~\ref{subsection3.4}
I indicate two further proofs for the existence of contact 
structures on every $3$--manifold (and provide references
to others). The one by Thur\-ston and Winkelnkemper
uses a description of $3$--manifolds as open books due to J.~Alexander; the
crucial idea in their proof is the one we also use to
simplify Martinet's argument. Indeed, my discussion of the
Lutz twist in the present section owes more to the paper by
Thurston-Winkelnkemper than to any other reference. The second proof, by
J.~Gonzalo, is based on a branched cover description of $3$--manifolds
found by H.~Hilden, J.~Montesinos and T.~Thickstun.
This branched cover description
also yields a very simple geometric proof that every orientable
$3$--manifold
is parallelisable.

In Section~\ref{subsection3.5} we discuss the fundamental dichotomy
between tight and overtwisted contact structures, introduced by
Eliashberg, as well as the relation of these types of contact
structures with the concept of symplectic fillability.
The chapter concludes in
Section~\ref{subsection3.6} with a survey of classification results
for contact structures on $3$--manifolds.

But first we discuss, in
Section~\ref{subsection3.0}, an invariant of
transverse knots in $\R^3$ with its standard contact structure.
This invariant will be an ingredient in the proof of the Lutz-Martinet
theorem, but is also of independent interest.

I do not feel embarrassed to use quite a bit of machinery from
algebraic and differential topology in this chapter. However, I believe
that nothing that cannot be found in such standard texts as \cite{bred93},
\cite{kosi93} and \cite{mist74} is used without proof or an explicit reference.

Throughout this third section $M$ denotes a closed, orientable 3-manifold.
\subsection{An invariant of transverse knots}
\label{subsection3.0}
Although the invariant in question can be defined for transverse knots
in arbitrary contact manifolds (provided the knot is homologically
trivial),
for the sake of clarity I restrict attention to transverse
knots in $\R^3$ with its standard contact structure $\xi_0=\ker
(dz+x\, dy)$.
This will be sufficient for the proof of the Lutz-Martinet theorem.
In Section~\ref{subsection3.6} I say a few words about the
general situation and related invariants for Legendrian knots.

Thus, let $\gamma$ be a transverse knot in $(\R^3,\xi_0)$.
Push $\gamma$ a little in the direction of $\partial_x$ -- notice that
this is a nowhere zero vector field contained in~$\xi_0$, and in particular
transverse to~$\gamma$ -- to
obtain a knot~$\gamma '$. An orientation of $\gamma$ induces an
orientation of $\gamma'$. The orientation of $\R^3$ is given by
$dx\wedge dy\wedge dz$.

\begin{defn}
The {\bf self-linking number} $l(\gamma )$ of the transverse knot
$\gamma$ is the linking number of $\gamma$ and $\gamma'$.
\end{defn}

Notice that this definition is independent of the choice of
orientation of~$\gamma$ (but it changes sign if the orientation
of $\R^3$ is reversed). Furthermore,
in place of $\partial_x$ we could have chosen any
nowhere zero vector field $X$ in $\xi_0$ to define $l(\gamma )$:
The difference between the the self-linking number defined via
$\partial_x$ and that defined via $X$ is the degree of the map
$\gamma\rightarrow S^1$ given by associating to a point on $\gamma$
the angle between $\partial_x$ and $X$. This degree is computed with
the induced map $\Z\cong H_1(\gamma )\rightarrow H_1(S^1)\cong Z$.
But the map $\gamma\rightarrow S^1$ factors through $\R^3$, so the
induced homomorphism on homology is the zero homomorphism.

Observe that
$l(\gamma )$ is an invariant under isotopies of $\gamma$ within
the class of transverse knots.

\begin{figure}[h]
\centerline{\relabelbox\small
\epsfxsize 8cm \epsfbox{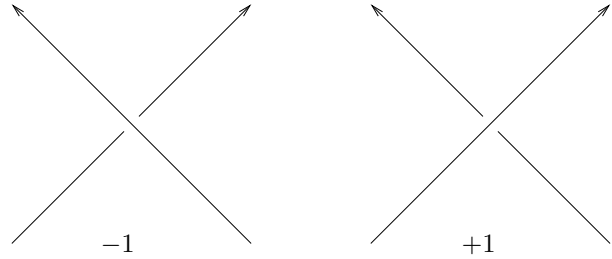}
\extralabel <-2cm,-.1cm> {$+1$}
\extralabel <-6.8cm,-.1cm> {$-1$}
\endrelabelbox}
\caption{Signs of crossings in a knot diagram.}
\label{figure:sign}
\end{figure}

We now want to compute $l(\gamma )$ from the front projection
of~$\gamma$. Recall that the {\bf writhe} of an oriented knot diagram
is the signed number of self-crossings of the diagram,
where the sign of the crossing is given in Figure~\ref{figure:sign}.

\begin{lem}
The self-linking number $l(\gamma )$ of a transverse knot is
equal to the writhe $w(\gamma )$ of its front projection.
\end{lem}

\begin{proof}
Let $\gamma'$ be the push-off of $\gamma$ as described. Observe that
each crossing of the front projection of $\gamma$ contributes
a crossing of $\gamma'$ underneath $\gamma$ of the corresponding sign.
Since the linking number of $\gamma$ and $\gamma'$ is equal
to the signed number of times that $\gamma'$ crosses
underneath~$\gamma$ (cf.~\cite[p.~37]{save99}),
we find that this linking number
is equal to the signed number of self-crossings of~$\gamma$,
that is, $l(\gamma )=w(\gamma )$.
\end{proof}

\begin{prop}
\label{prop:linking}
Every self-linking number is realised by a transverse link in standard $\R^3$.
\end{prop}

\begin{figure}[h]
\centerline{\relabelbox\small
\epsfxsize 10cm \epsfbox{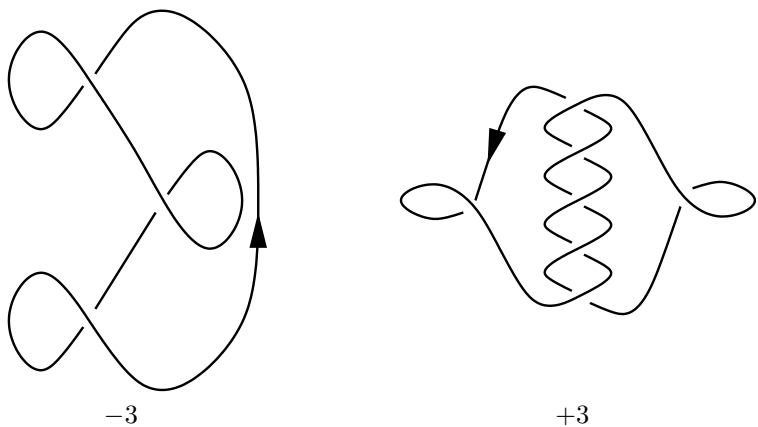}
\extralabel <-2.7cm,-.4cm> {$+3$}
\extralabel <-8.7cm,-.4cm> {$-3$}
\endrelabelbox}
\caption{Transverse knots with self-linking number $\pm 3$.}
\label{figure:linking}
\end{figure}

\begin{proof}
Figure~\ref{figure:linking} shows front projections of
positively transverse knots
(cf.\ Section~\ref{section:transverse}) with self-linking number~$\pm 3$.
From that
the construction principle for realising any odd integer
should be clear.
With a two component link any even integer can be realised.
\end{proof}

\begin{rem}
{\rm It is no accident that I do not give an example of a transverse
knot with {\it even} self-linking number. By a theorem of
Eliashberg~\cite[Prop.~2.3.1]{elia93a} that relates $l(\gamma )$
to the Euler characteristic of a Seifert surface $S$ for $\gamma$ and
the signed number of singular points of the characteristic
foliation~$S_{\xi}$, the self-linking number $l(\gamma )$ of
a knot can only take {\it odd} values.
}
\end{rem}

\subsection{Martinet's construction}
\label{subsection3.1}
According to Lickorish~\cite{lick62} and Wallace~\cite{walla60} $M$ can be
obtained from $S^3$ by Dehn surgery along a link of 1--spheres. This means
that we have to remove a disjoint set of embedded solid tori $S^1\times D^2$ 
from $S^3$ and glue back solid tori with suitable identification by a
diffeomorphism along the boundaries $S^1\times S^1$. The effect of such
a surgery (up to diffeomorphism of the resulting manifold)
is completely determined by the induced map in homology
\begin{eqnarray*}
H_1(S^1\times\partial D^2) & \longrightarrow & H_1(S^1\times\partial D^2)\\
{\Z}\oplus {\Z} & \longrightarrow & {\Z}\oplus {\Z},
\end{eqnarray*}
which is given by a unimodular matrix $\left( \begin{array}{cc}n&q\\m&p
\end{array}\right) \in \mbox{\rm GL}(2,{\Z})$.
Hence, denoting coordinates in $S^1\times S^1$ by $(\theta ,\varphi )$,
we may always assume the identification maps to be of the form
\[ \left( \begin{array}{c}\theta\\ \varphi\end{array}\right) \longmapsto
\left( \begin{array}{cc}n&q\\m&p\end{array}\right) \left( \begin{array}{c}
\theta\\ \varphi \end{array}\right) . \]
The curves $\mu$ and $\lambda$ on $\partial (S^1\times D^2)$ given
respectively by $\theta =0$ and $\varphi =0$ are called
{\it meridian} and {\it longitude}. We keep the same notation
$\mu$, $\lambda$ for the homology classes these curves
represent. It turns out that the
diffeomorphism type of the surgered manifold is completely
determined by the class $p\mu +q\lambda$, which is the class
of the curve that becomes homotopically trivial in the
surgered manifold (cf.~\cite[p.~28]{save99}). In fact, the
Dehn surgery is completely determined by the surgery coefficient
$p/q$, since the diffeomorphism of $\partial (S^1\times D^2)$ given
by $(\lambda ,\mu )\mapsto (\lambda ,-\mu )$ extends to a diffeomorphism
of the solid torus that we glue back.

Similarly, the diffeomorphism given by $(\lambda ,\mu )\mapsto
(\lambda +k\mu ,\mu )$ extends. By such a change of longitude in
$S^1\times D^2\subset M$, which amounts to choosing a different
trivialisation of the normal bundle (i.e.\ {\it framing})
of $S^1\times\{ 0\}\subset M$,
the gluing map is changed to
$\left(\begin{array}{cc}n&q\\m-kn&p-kq\end{array}\right)$. By a change
of longitude in the solid torus that we glue back, the gluing map
is changed to $\left(\begin{array}{cc}n+kq&q\\m+kp&p\end{array}
\right)$. Thus, a Dehn surgery is a so-called handle surgery (or
`honest surgery' or simply `surgery')
if and only if the surgery coefficient is an integer, that is, $q=\pm 1$.
For in exactly this case we may assume $\left( \begin{array}{cc}n&q\\m&p
\end{array}\right)
=\left( \begin{array}{cc}0&1\\1&0\end{array}\right)$, and the surgery is
given by cutting out $S^1\times D^2$ and gluing back $S^1\times D^2$
with the identity map
\[ \partial (D^2\times S^1)\longrightarrow \partial (S^1\times D^2). \]

The theorem of Lickorish and Wallace remains true if we only allow handle
surgeries. However, this
assumption does not entail any great simplification of the existence proof
for contact structures, so we shall work with general Dehn surgeries.

Our aim in this section is to use this topological description of
3--manifolds for a proof of the following theorem, first proved by
Martinet~\cite{mart71}. The proof presented here is in spirit the one
given by Martinet, but, as indicated in the introduction to
this third section, amalgamated with ideas of Thurston and 
Winkelnkemper~\cite{thwi75}, whose proof of the same theorem we shall
discuss later.

\begin{thm}[Martinet]
Every closed, orientable $3$--manifold $M$ admits a contact structure.
\end{thm}

In view of the theorem of Lickorish and Wallace and the fact that $S^3$
admits a contact structure, Martinet's theorem is a direct
consequence of the following result.

\begin{thm}
\label{thm:c-surgery}
Let ${\xi}_0$ be a contact structure on a $3$--manifold $M_0$. Let $M$ be
the manifold obtained from $M_0$ by a Dehn surgery along a knot~$K$. Then
$M$ admits a contact structure $\xi$ which coincides with ${\xi}_0$
outside the neighbourhood of $K$ where we perform surgery.
\end{thm}

\noindent {\it Proof.}
By Theorem~\ref{thm:approximation} we may assume that $K$ is 
positively transverse to~${\xi}_0$. Then, by
the contact neighbourhood theorem (Example~\ref{exam:c-nbhd}),
we can find a tubular
neighbourhood of $K$ diffeomorphic to $S^1\times D^2(\delta_0)$, where
$K$ is identified with $S^1\times\{ 0\}$ and $D^2(\delta_0)$ denotes a disc
of radius~$\delta_0$, such that the contact structure ${\xi}_0$ is
given as the kernel of $d\overline{\theta}+\overline{r}^2d\overline{\varphi}$,
with $\overline{\theta}$ denoting the $S^1$--coordinate and $(\overline{r},
\overline{\varphi})$ polar coordinates on~$D^2(\delta_0)$.

Now perform a Dehn surgery along $K$ defined by the unimodular matrix
$\left( \begin{array}{cc}n&q\\m&p\end{array} \right)$. This corresponds to
cutting out $S^1\times D^2(\delta_1)\subset S^1\times D^2(\delta_0)$ for 
some $\delta_1<\delta_0$ and gluing it back in the way described
above.

Write $(\theta ;r,\varphi )$ for the coordinates on the copy of $S^1\times
D^2(\delta_0)$ that we want to glue back. Then the contact form
$d\overline{\theta}+\overline{r}^2d\overline{\varphi}$ given on
$S^1\times D^2(\delta_0)$ pulls back (along $S^1\times\partial D^2(\delta_1)$)
to
\[ d(n\theta +q\varphi )+r^2d(m\theta +p\varphi ). \]
This form is defined on all of $S^1\times (D^2(\delta_1)-\{ 0\})$,
and to complete the proof it only remains to find a contact form on
$S^1\times D^2(\delta_1)$ that coincides with this form near $S^1\times
\partial D^2(\delta_1)$. It is at this point that we use an argument
inspired by the Thurston-Winkelnkemper proof (but which goes back to Lutz).

\begin{lem}
Given a unimodular matrix $\left( \begin{array}{cc}n&q\\m&p\end{array}
\right)$, there is a contact form on $S^1\times D^2(\delta )$ that coincides
with $(n+mr^2)\, d\theta +(q+pr^2)\, d\varphi$ near $r=\delta$ and with
$\pm d\theta +r^2d\varphi$ near $r=0$.
\end{lem}

\begin{proof}
We make the ansatz
\[ \alpha = h_1(r)\, d\theta +h_2(r)\, d\varphi \]
with smooth functions $h_1(r),h_2(r)$. Then
\[ d\alpha = h_1'\, dr\wedge d\theta +h_2'\, dr\wedge d\varphi \]
and 
\[ \alpha\wedge d\alpha = \left| \begin{array}{cc} h_1&h_2\\ h_1'&h_2'
\end{array}\right| \, d\theta\wedge dr\wedge d\varphi .\]
So to satisfy the contact condition $\alpha\wedge d\alpha \neq 0$ all we
have to do is to find a parametrised curve
\[ r\longmapsto (h_1(r),h_2(r)),\;\; 0\leq r\leq \delta, \]
in the plane satisfying the following conditions:
\begin{itemize}
\item[1.] $h_1(r)=\pm 1$ and $h_2(r)=r^2$ near $r=0$,
\item[2.] $h_1(r)=n+mr^2$ and $h_2(r)=q+pr^2$ near $r=\delta$,
\item[3.] $(h_1(r),h_2(r))$ is never parallel to $(h_1'(r),h_2'(r))$.
\end{itemize}
Since $np-mq=\pm 1$, the vector $(m,p)$ is not a multiple of $(n,q)$.
Figure~\ref{figure:Dehn}
shows possible solution curves for the two cases $np-mq=\pm 1$.

\begin{figure}[h]
\centerline{\relabelbox\small
\epsfxsize 12cm \epsfbox{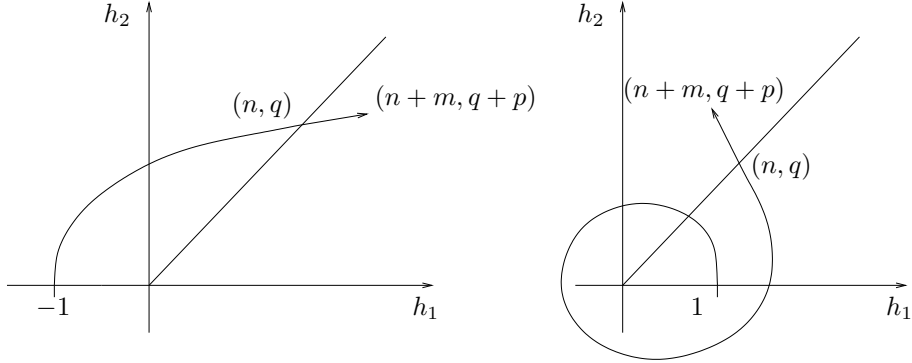}
\extralabel <-.3cm,.6cm> {$h_1$}
\extralabel <-6.6cm,.6cm> {$h_1$}
\extralabel <-2.9cm,.6cm> {$1$}
\extralabel <-11.6cm,.6cm> {$-1$}
\extralabel <-4.4cm,4.5cm> {$h_2$}
\extralabel <-10.7cm,4.5cm> {$h_2$}
\extralabel <-9cm,3.3cm> {$(n,q)$}
\extralabel <-7.1cm,3.4cm> {$(n+m,q+p)$}
\extralabel <-2.1cm,2.5cm> {$(n,q)$}
\extralabel <-3.8cm,3.5cm> {$(n+m,q+p)$}
\endrelabelbox}
\caption{Dehn surgery.}
\label{figure:Dehn}
\end{figure}

\noindent
This completes the proof of the lemma and hence that of
Theorem~\ref{thm:c-surgery}.
\end{proof}

\begin{rem}
{\rm
On $S^3$ we have the standard contact forms $\alpha_{\pm}=x\, dy-y\, dx
\pm (z\, dt-t\, dz)$ defining opposite orientations
(cf.~Remark~\ref{rem:posneg}). Performing the above
surgery construction either on $(S^3,\ker\alpha_+)$ or on $(S^3,\ker\alpha_-)$
we obtain both positive and negative contact structures on any given~$M$.
The same is true for the Lutz construction that we study in the next two
sections. Hence: {\it A closed oriented $3$--manifold admits both a positive
and a negative contact structure in each homotopy class of
tangent $2$--plane fields.}}
\end{rem}
\subsection{$2$--plane fields on $3$--manifolds}
\label{subsection3.2}
First we need the following well-known fact.

\begin{thm}
\label{thm:parallelisable}
Every closed, orientable $3$--manifold $M$ is parallelisable.
\end{thm}

\noindent {\em Remark.}
The most geometric proof of this theorem can be given based on a structure
theorem of Hilden, Montesinos and Thickstun. This will be discussed in
Section~\ref{section:branched}.
Another proof can be found in~\cite{kirb89}. Here we
present the classical algebraic proof.

\begin{proof}
The main point is to show the vanishing of the second Stiefel-Whitney class
$w_2(M)=w_2(TM)\in H^2(M;{\Z}_2)$. Recall the following
facts, which can be found in~\cite{bred93}; for the interpretation of
Stiefel-Whitney classes as obstruction classes see also~\cite{mist74}.

There are Wu classes $v_i\in H^i(M;{\Z}_2)$ defined by
\[ \langle \mbox{\rm Sq}^i(u),[M]\rangle =\langle v_i\cup u ,[M]\rangle \]
for all $u\in H^{3-i}(M;{\Z}_2)$, where $\mbox{\rm Sq}$ denotes the
Steenrod squaring operations. Since $\mbox{\rm Sq}^i(u)=0$ for
$i>3-i$, the only (potentially)
non-zero Wu classes are $v_0=1$ and $v_1$. The Wu classes
and the Stiefel-Whitney classes are related by $w_q=\sum_j
\mbox{\rm Sq}^{q-j}(v_j)$. Hence $v_1=\mbox{\rm Sq}^0(v_1)=w_1$, which
equals zero since $M$ is orientable. We conclude $w_2=0$.

Let $V_2({\R}^3)=\SO (3)/\SO (1)=\SO (3)$ be the Stiefel manifold of
oriented, orthonormal
2--frames in~${\R}^3$. This is connected, so there exists a section 
over the 1--skeleton of $M$ of
the 2--frame bundle $V_2(TM)$ associated with~$TM$ (with a choice of
Riemannian metric on $M$ understood\footnote{This is not
necessary, of course. One may also work with arbitrary $2$--frames
without reference to a metric. This does not affect the homotopical
data.}).
The obstruction to extending 
this section over the 2--skeleton is equal to~$w_2$, which vanishes as we
have just seen.
The obstruction to extending the section over all of $M$ lies in
$H^3(M;\pi_2(V_2({\R}^3)))$, which is the zero group because of
$\pi_2(\SO (3))=0$.

We conclude that $TM$ has a trivial 2--dimensional sub-bundle $\varepsilon^2$.
The complementary 1--dimensional bundle $\lambda =TM/\varepsilon^2$ is
orientable and hence trivial since $0=w_1(TM)=w_1(\varepsilon^2) +
w_1(\lambda )=w_1(\lambda )$. Thus $TM=\varepsilon^2\oplus \lambda$
is a trivial bundle.
\end{proof}

Fix an arbitrary Riemannian metric
on $M$ and a trivialisation of the unit tangent bundle $STM\cong 
M\times S^2$. This sets up a one-to-one correspondence between the following
sets, where all maps, homotopies etc.\ are understood to be smooth.
\begin{itemize}
\item Homotopy classes of unit vector fields $X$ on $M$,
\item Homotopy classes of (co-)oriented 2--plane distributions $\xi$
in~$TM$,
\item Homotopy classes of maps $f\co M\rightarrow S^2$. 
\end{itemize}
(I use the term `$2$--plane distribution' synomymously with
`$2$--dimensional sub-bundle of the tangent bundle'.)
Let ${\xi}_1,{\xi}_2$ be two arbitrary 2--plane distributions (always 
understood to be cooriented).
By elementary obstruction theory there is an obstruction
\[ d^2({\xi}_1,{\xi}_2)\in H^2(M;\pi_2(S^2))\cong H^2(M;{\Z}) \]
for ${\xi}_1$ to be homotopic to ${\xi}_2$ over the 2--skeleton
of $M$ and, if $d^2({\xi}_1,{\xi}_2)=0$ and after homotoping $\xi_1$ to
$\xi_2$ over the $2$--skeleton, an obstruction (which will depend, in
general, on that first homotopy)
\[ d^3({\xi}_1,{\xi}_2)\in H^3(M;\pi_3(S^2))\cong H^3(M;{\Z})
\cong {\Z} \]
for ${\xi}_1$ to be homotopic to ${\xi}_2$ over all of $M$.
(The identification of $H^3(M;{\Z})$ with $\Z$ is determined by the
orientation of~$M$.) However,
rather than relying on general obstruction theory, we shall interpret
$d^2$ and $d^3$ geometrically, which will later allow us to give a geometric 
proof that every homotopy class of 2--plane fields $\xi$ on $M$ contains
a contact structure.

The only fact that I want to quote here is that, by the Pontrjagin-Thom
construction, homotopy classes of maps $f\co M\rightarrow S^2$ are in
one-to-one correspondence with framed cobordism classes of framed 
(and oriented) links of 1--spheres in $M$. The general theory can be found
in~\cite{bred93} and~\cite{kosi93};
a beautiful and elementary account is given
in~\cite{miln65}. 

For given $f$, the correspondence is defined by choosing a regular value
$p\in S^2$ for $f$ and a positively oriented basis $\mathfrak b$ of $T_pS^2$,
and associating with it the oriented framed link $(f^{-1}(p),f^*{\mathfrak b})$,
where $f^*{\mathfrak b}$ is the pull-back of $\mathfrak b$ under the fibrewise
bijective map $Tf: T(f^{-1}(p))^{\perp}\rightarrow T_pS^2$.
The orientation of $f^{-1}(p)$ is the one which together with the frame
$f^*{\mathfrak b}$ gives the orientation of~$M$.

For a given framed link $L$ the corresponding $f$ is defined by projecting 
a (trivial) disc bundle neighbourhood $L\times D^2$ of $L$ in $M$ onto
the fibre $D^2\cong S^2-p^*$, where $0$ is identified with $p$ and $p^*$
denotes the antipode of~$p$, and sending $M-(L\times D^2)$ to $p^*$. Notice
that the orientations of $M$ and the components of $L$ determine that of
the fibre~$D^2$, and hence determine the map~$f$.

Before proceeding to define the obstruction classes $d^2$ and $d^3$ we make
a short digression and discuss some topological background material which is
fairly standard but not contained in our basic textbook
references \cite{bred93} and~\cite{kosi93}.
\subsubsection{Hopf's Umkehrhomomorphismus}
If $f\co M^m\rightarrow N^n$ is a continuous map between smooth manifolds,
one can define a homomorphism $\varphi :H_{n-p}(N)\rightarrow 
H_{m-p}(M)$ on homology classes 
represented by submanifolds as follows. Given a homology class $[L]_N\in
H_{n-p}(N)$ represented by a codimension $p$ submanifold~$L$, replace $f$ by
a smooth approximation transverse to $L$ and define $\varphi ([L]_N)=
[f^{-1}(L)]_M$. This is essentially Hopf's
{\em Umkehrhomomorphismus} \cite{hopf30}, except that he worked with 
combinatorial manifolds of equal dimension and made no assumptions on the
homology class. The following theorem, which in spirit is contained 
in~\cite{freu37}, shows that $\varphi$ is independent of choices (of
submanifold $L$ representing a class and smooth transverse approximation to~$f$)
and actually a homomorphism of intersection rings. This statement is not as
well-known as it should be, and I only know of a proof in the literature
for the special case where $L$ is a point~\cite{gott78}. In \cite{bred93}
this map is called {\em transfer map} (more general transfer maps
are discussed in~\cite{gott78}), but is only defined
indirectly via Poincar\'e duality (though implicitly the statement of the
following theorem is contained in~\cite{bred93}, see for instance p.~377).

\begin{thm}
Let $f\co M^m\rightarrow N^n$ be a smooth map between
closed, oriented mani\-folds and $L^{n-p}\subset N^n$ a closed,
oriented submanifold of codimension~$p$ such that
$f$ is transverse to~$L$. Write $u\in H^p(N)$ for the Poincar\'e dual of
$[L]_N$, that is, $u\cap [N]=[L]_N$. Then $[f^{-1}(L)]_M=
f^*u\cap [M]$. In other words: If $u$ is Poincar\'e dual to
$[L]_N$, then $f^*u\in H^p(M)$ is Poincar\'e dual to $[f^{-1}(L)]_M$.
\end{thm}

\begin{proof}
Since $f$ is transverse to $L$, the differential 
$Tf$ induces a fibrewise isomorphism between the normal bundles of
$f^{-1}(L)$ and~$L$, and we find (closed) tubular neighbourhoods 
$W\rightarrow L$ and $V=f^{-1}(W)\rightarrow f^{-1}(L)$ (considered as
disc bundles) such that $f\co V\rightarrow W$ is a fibrewise isomorphism.
Write $[V]_0$ and $[W]_0$ for the orientation classes in $H_m(V,V-f^{-1}(L))$
and $H_n(W,W-L)$, respectively. We can identify these
homology groups with $H_m(V,\partial V)$ and $H_n(W,\partial W)$, respectively.
Let $\tau_W\in H^p(W,\partial W)$ and $\tau_V\in H^p(V,\partial V)$ be the Thom
classes of these disc bundles defined by
\begin{eqnarray*}
\tau_W\cap [W]_0 & = & [L]_N, \\
\tau_V\cap [V]_0 & = & [f^{-1}(L)]_M.
\end{eqnarray*}
Notice that $f^*\tau_W=\tau_V$ since 
$f\co W\rightarrow V$ is fibrewise isomorphic and
the Thom class of an oriented disc bundle
is the unique class whose restriction to each fibre is a positive
generator of $H^p(D^p,\partial D^p)$.
Writing $i\co M\rightarrow (M,M-f^{-1}(L))$ and $j\co N\rightarrow(N,N-L)$
for the inclusion maps we have
\begin{eqnarray*}
[f^{-1}(L)]_M & = & \tau_V\cap [V]_0\\
  & = & f^*\tau_W\cap [V]_0 \\
  & = & f^*\tau_W\cap i_*[M],
\end{eqnarray*}
where we identify $H_m(M,M-f^{-1}(L))$ with $H_m(V,V-f^{-1}(L))$ under the
excision isomorphism. Then we have further
\begin{eqnarray*}
[f^{-1}(L)]_M & = & i^*f^*\tau_W\cap [M]\\
  & = & f^*j^*\tau_W\cap [M].
\end{eqnarray*}
So it remains to identify $j^*\tau_W$ as the Poincar\'e dual $u$
of~$[L]_N$. Indeed,
\begin{eqnarray*}
j^*\tau_W\cap [N] & = & \tau_W\cap j_*[N]\\
  & = & \tau_W\cap [W]_0\\
  & = & [L]_N,
\end{eqnarray*}
where we have used the excision isomorphism between the groups
$H_n(W,W-L)$ and $H_n(N,N-L)$.
\end{proof}
\subsubsection{Representing homology classes by submanifolds}
\label{section:rhcbs}
We now want to relate
elements in $H_1(M;{\Z})$ to cobordism classes of links in~$M$.

\begin{thm}
\label{thm:replink}
Let $M$ be a closed, oriented $3$--manifold. Any
$c\in H_1(M;{\Z})$ 
is represented by an embedded, oriented link (of $1$--spheres) $L_c$
in~$M$. Two links $L_0,L_1$ represent the same class
$[L_0]=[L_1]$ if and only if
they are cobordant in~$M$, that is, there is an embedded,
oriented surface $S$ in $M\times [0,1]$ with
\[ \partial S=L_1\sqcup (-L_0)\subset M\times\{ 1\} \sqcup M\times\{ 0\} ,\]
where $\sqcup$ denotes disjoint union.
\end{thm}

\begin{proof}
Given $c\in H_1(M;{\Z})$, set $u =
PD(c)\in H^2(M;{\Z})$, where $PD$ denotes the Poincar\'e duality
map from homology to cohomology. There is a well-known isomorphism
\[ H^2(M;{\Z})\cong [M,K({\Z},2)]=[M,{\C}P^{\infty}], \]
where brackets denote homotopy classes of maps (cf.~\cite[VII.12]{bred93}). 
So $u$ corresponds to a homotopy class of maps $[f]\co
M\rightarrow {\C}P^{\infty}$ such that
$f^*u_0 =u$, where $u_0$ is the positive generator of $H^2(
{\C}P^{\infty})$ (that is, the one that pulls back to the
Poincar\'e dual of $[{\C}P^{k-1}]_{{\C}P^k}$ under the natural
inclusion ${\C}P^k\subset {\C}P^{\infty}$). 
Since $\dim M=3$, any map $f\co M\rightarrow 
{\C}P^{\infty}$ is homotopic to a smooth map $f_1\co M\rightarrow {\C}P^1$. Let
$p$ be a regular value of~$f_1$. Then
\[ PD(c)=u=f_1^*u_0=f_1^*PD[p]=PD[f_1^{-1}(p)] \]
by our discussion above, and hence $c=[f_1^{-1}(p)]$. So $L_c=f_1^{-1}(p)$
is the desired link.

It is important to note that in spite of what we have just said it is not true
that $[M,{\C}P^{\infty}]=[M,{\C}P^1]$, since a map
$F\co M\times [0,1]\rightarrow {\C}P^{\infty}$ with $F(M\times\{ 0,1\})
\subset {\C}P^1$ is not, in general,
homotopic $\mbox{\rm rel}\, (M\times\{ 0,1\} )$ to
a map into ${\C}P^1$. However, we do have $[M,{\C}P^{\infty}]=
[M,{\C}P^2]$.

If two links $L_0,L_1$ are cobordant in~$M$, then clearly
\[ [L_0]=[L_1] \in H_1(M\times [0,1];{\Z})\cong H_1(M;{\Z}). \]
For the converse,
suppose we are given two links $L_0,L_1\subset M$ with $[L_0]=[L_1]$.
Choose arbitrary framings for these links and use this, as described at
the beginning of this section, to define smooth
maps $f_0,f_1\co M\rightarrow S^2$
with common regular value $p\in S^2$ such that $f_i^{-1}(p)=L_i,
i=0,1$. Now identify $S^2$ with the standardly embedded ${\C}P^1
\subset {\C}P^2$. Let $P\subset {\C}P^2$ be a second copy of
${\C}P^1$, embedded in such a way that $[P]_{{\C}P^2}=
[{\C}P^1]_{{\C}P^2}$ and $P$ intersects ${\C}P^1$ transversely
in $p$ only. This is possible since ${\C}P^1\subset {\C}P^2$ has
self-intersection one. Then the maps $f_0,f_1$, regarded as maps into
${\C}P^2$, are transverse to $P$ and we have $f_i^{-1}(P)=L_i, i=0,1$.
Hence
\begin{eqnarray*}
f_i^*u_0 & = & f_i^*(PD[P]_{{\C}P^2})= PD[f_i^{-1}(P)]_M\\
  & = & PD[L_i]_M
\end{eqnarray*}
is the same for $i=0$ or 1, and from the identification
\begin{eqnarray*}
[M,{\C}P^2] & \stackrel{\cong}{\longrightarrow} & H^2(M,{\Z})\\
\mbox{$[f]$} & \longmapsto & f^*u_0
\end{eqnarray*}
we conclude that $f_0$ and $f_1$ are homotopic as maps into~${\C}P^2$.

Let $F\co M\times [0,1]\rightarrow {\C}P^2$ be a homotopy between $f_0$
and $f_1$, which we may assume to be constant near $0$ and~$1$. This
$F$ can be smoothly approximated by a map $F'\co M\times [0,1]\rightarrow
{\C}P^2$ which is transverse to $P$ and coincides with $F$ near
$M\times 0$ and $M\times 1$ (since there the transversality condition was
already satisfied). In particular, $F'$ is still a homotopy between
$f_0$ and~$f_1$, and $S=(F')^{-1}(P)$ is a surface with the desired
property
$\partial S=L_1\sqcup (-L_0).$
\end{proof}

Notice that in the course of this proof we have observed that cobordism classes
of links in $M$ (equivalently, classes in $H_1(M;{\Z})$) correspond to
homotopy classes of maps $M\rightarrow {\C}P^2$, whereas framed 
cobordism classes of framed links correspond to homotopy classes of
maps $M\rightarrow {\C}P^1$.

By forming the connected sum of the components of a link
representing a certain class in $H_1(M;{\Z})$, one may actually
always represent such a class by a link with only one component, that is, 
a knot.
\subsubsection{Framed cobordisms}
We have seen that if $L_1,L_2\subset M$
are links with $[L_1]=[L_2]\in H_1(M;{\Z})$, then $L_1$ and $L_2$
are cobordant in $M$. In general, however, a given framing on $L_1$ and $L_2$
does not extend over the cobordism. The following observation will
be useful later on.

Write $(S^1,n)$ for a contractible loop in $M$ with framing $n\in{\Z}$
(by which we mean that $S^1$ and a second copy of $S^1$ obtained by
pushing it away in the direction of one of the vectors in
the frame have linking number~$n$). When
writing $L=L'\sqcup (S^1,n)$ it is understood that $(S^1,n)$ is not
linked with any component of~$L'$.

\begin{figure}[h]
\centerline{\relabelbox\small
\epsfxsize 8cm \epsfbox{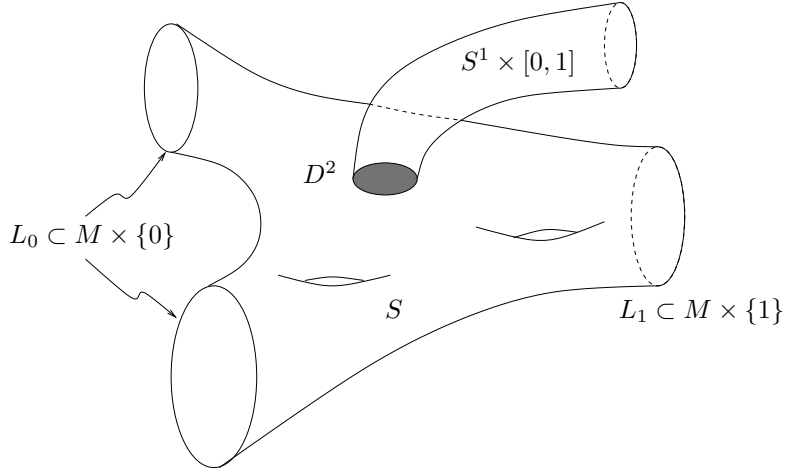}
\extralabel <-.9cm,2cm> {$L_1\subset M\times\{ 1\}$}
\extralabel <-9cm,3cm> {$L_0\subset M\times\{ 0\}$}
\extralabel <-5.1cm,3.8cm> {$D^2$}
\extralabel <-3cm,5.3cm> {$S^1\times [0,1]$}
\extralabel <-4cm,2cm> {$S$}
\endrelabelbox}
\caption{The framed cobordism between $L_0$ and $L_1\sqcup (S^1,n)$.}
\label{figure:cobord}
\end{figure}

Suppose we have two framed links $L_0,L_1\subset M$ with $[L_0]=[L_1]$.
Let $S\subset M\times [0,1]$ be an embedded surface with 
\[ \partial S=L_1\sqcup (-L_0)\subset M\times\{ 1\} \sqcup M\times\{ 0\} .\]
With $D^2$ a small disc embedded in $S$, the framing of $L_1$
and $L_2$ in $M$ extends to a framing of $S-D^2$ in $M\times [0,1]$ (since
$S-D^2$ deformation retracts to a $1$--dimensional complex containing
$L_0$ and~$L_1$, and over such a complex
an orientable $2$--plane bundle is trivial). Now
we embed a cylinder $S^1\times [0,1]$ in $M\times [0,1]$ such that
\[ S^1\times [0,1]\cap M\times \{ 0\} =\emptyset, \]
\[ S^1\times [0,1]\cap M\times\{ 1\} =S^1\times \{ 1\} ,\]
and
\[ S^1\times [0,1]\cap (S-D^2)=S^1\times\{ 0\} =\partial D^2, \]
see Figure~\ref{figure:cobord}.
This shows that $L_0$ is framed cobordant in $M$ to $L_1\sqcup (S^1,n)$
for suitable $n\in {\Z}$. 

\subsubsection{Definition of the obstruction classes}
We are now in a position to define the obstruction classes $d^2$ and~$d^3$.
With a choice of Riemannian metric on $M$ and a trivialisation
of $STM$ understood,
a 2--plane distribution ${\xi}$ on $M$ defines a map $f_{\xi}\co M\rightarrow
S^2$ and hence an oriented framed link $L_{\xi}$ as described above. Let
$[L_{\xi}]\in H_1(M;{\Z})$ be the homology class represented 
by~$L_{\xi}$. This only depends on the homotopy class of~$\xi$, since
under homotopies of $\xi$ or choice of different regular
values of $f_{\xi}$ the cobordism class of $L_{\xi}$ remains
invariant. We define
\[ d^2(\xi_1,\xi_2)=PD[L_{\xi_1}]-PD[L_{\xi_2}]. \]
With this definition $d^2$ is clearly additive, that is,
\[ d^2({\xi}_1,{\xi}_2)+d^2({\xi}_2,{\xi}_3)=
d^2({\xi}_1,{\xi}_3). \]
The following lemma shows that $d^2$ is indeed the desired obstruction class.

\begin{lem}
The $2$--plane distributions ${\xi}_1$ and ${\xi}_2$ are
homotopic over the $2$--skeleton $M^{(2)}$ of~$M$ if and only if 
$d^2(\xi_1,\xi_2)=0$.
\end{lem}

\begin{proof}
Suppose $d^2(\xi_1,\xi_2)=0$, that is, $[L_{{\xi}_1}]=
[L_{{\xi}_2}]$. By Theorem~\ref{thm:replink}
we find a surface $S$ in $M\times [0,1]$ with
\[ \partial S=L_{{\xi}_2}\sqcup (-L_{{\xi}_1})\subset
M\times\{ 1\} \sqcup M\times\{ 0\} .\]
From the discussion on framed cobordism above we know that
for suitable $n\in {\Z}$ we find
a {\em framed} surface $S'$ in $M\times [0,1]$ such that 
\[ \partial S'=\left( L_{{\xi}_2}\sqcup (S^1,n)\right) \sqcup 
(-L_{{\xi}_1})\subset M\times\{ 1\} \sqcup M\times\{ 0\} \]
as framed manifolds.

Hence ${\xi}_1$ is homotopic to a 2--plane distribution ${\xi}_1'$
such that $L_{{\xi}_1'}$ and $L_{{\xi}_2}$ differ only by one
contractible framed loop (not linked with any other component). Then the
corresponding maps $f_1',f_2$ differ only in a neighbourhood of this loop,
which is contained in a 3--ball, so $f_1'$ and $f_2$ (and hence ${\xi}_1'$
and ${\xi}_2$) agree over the 2--skeleton.

Conversely, if $\xi_1$ and $\xi_2$ are homotopic over~$M^{(2)}$, we may
assume $\xi_1=\xi_2$ on $M-D^3$ for some embedded 3--disc $D^3\subset M$
without changing $[L_{\xi_1}]$ and~$[L_{\xi_2}]$. Now $[L_{\xi_1}]=
[L_{\xi_2}]$ follows from $H_1(D^3,S^2)=0$.
\end{proof}

\begin{rem}
{\rm
By \cite[\S~37]{stee51} the obstruction to homotopy
between $\xi$ and ${\xi}_0$ (corresponding to the constant map
$f_{\xi_0}\co M\rightarrow S^2$) over the 2--skeleton of $M$ is given by
$f_{\xi}^*u_0$, where $u_0$ is the positive generator of $H^2(S^2;{\Z})$.
So $u_0 =PD[p]$ for any $p\in S^2$, and taking $p$ to be
a regular value of $f_{\xi}$ we have
\begin{eqnarray*}
f_{\xi}^*u_0 & = & f_{\xi}^*PD[p]\, =\, PD[f_{\xi}^{-1}(p)]\\
  & = & PD[L_{\xi}]\, =\, d^2({\xi},{\xi}_0).
\end{eqnarray*}
This gives an alternative way to see that our geometric definition of
$d^2$ does indeed coincide with the class defined by classical obstruction
theory.
}
\end{rem}

Now suppose $d^2({\xi}_1,{\xi}_2)=0$. We may then assume that ${\xi}_1
={\xi}_2$ on $M-\mbox{\rm int}(D^3)$, and we define
$d^3({\xi}_1,{\xi}_2)$ to be the Hopf invariant $H(f)$ of the map
$f\co S^3\rightarrow S^2$ defined as $f_1\circ\pi_+$ on the upper
hemisphere and $f_2\circ\pi_-$ on the lower hemisphere, where $\pi_+,\pi_-$ are
the orthogonal projections of the upper resp.\ lower hemisphere onto
the equatorial disc, which we identify with $D^3\subset M$.
Here, given an orientation of~$M$,
we orient $S^3$ in such a way that $\pi_+$ is orientation-preserving
and $\pi_-$ orientation-reversing; the orientation of $S^2$ is
inessential for the computation of~$H(f)$.
Recall that $H(f)$ is defined as the linking number of the preimages of
two distinct regular values of a smooth map homotopic to~$f$. Since
the Hopf invariant classifies homotopy classes of maps $S^3\rightarrow S^2$
(it is in fact an isomorphism $\pi_3(S^2)\rightarrow {\mathbb Z}$),
this is a suitable definition for the obstruction class~$d^3$.
Moreover, the homomorphism property of $H(f)$ and the way addition in
$\pi_3(S^2)$ is defined entail the additivity of $d^3$ analogous
to that of~$d^2$.

For $M=S^3$ there is another way to interpret $d^3$. Oriented
$2$--plane distributions
on $M$ correspond to sections of the bundle associated to $TM$
with fibre $\SO (3)/\U (1)$, hence to maps $M\rightarrow
\SO (3)/\U (1)\cong S^2$ since
$TM$ is trivial. Similarly, almost complex structures on $D^4$ correspond
to maps $D^4\rightarrow\SO (4)/\U (2)\cong\SO (3)/\U (1)$ (cf.~\cite{gray59}
for this isomorphism). A cooriented $2$--plane distribution on $M$
can be interpreted as a triple $(X,\xi ,J)$, where $X$ is a vector field
transverse to $\xi$ defining the coorientation, and $J$
a complex structure on $\xi$ defining the orientation. Such a triple
is called an {\bf almost contact structure}. (This notion generalises
to higher (odd) dimensions, and by Remark~\ref{rem:symplectic} every
{\it cooriented} contact structure induces an almost contact
structure, and in fact a unique one up to homotopy as follows
from the result cited in that remark.) Given an almost
contact structure $(X,\xi ,J)$ on $S^3$, one defines an almost complex
structure $\widetilde{J}$ on $TD^4|S^3$ by $\widetilde{J}|\xi =J$ and
$\widetilde{J}N=X$, where $N$ denotes the outward normal vector field.
So there is a canonical way to identify homotopy classes of almost contact
structures on $S^3$ with elements of $\pi_3(\SO (3)/\U (1))\cong {\Z}$
such that the value zero corresponds to the almost contact structure
that extends as almost complex structure over~$D^4$.
\subsection{Let's Twist Again}
\label{subsection3.3}
Consider a 3--manifold $M$ with cooriented contact structure~$\xi$ and an
oriented 1--sphere $K\subset M$ embedded transversely to $\xi$ such that
the positive orientation of $K$ coincides with the positive coorientation
of~$\xi$. Then in suitable local coordinates we can identify $K$
with $S^1\times\{ 0\}\subset S^1\times D^2$ such that ${\xi}=\ker (d\theta
+r^2d\varphi )$ and $\partial_{\theta}$ corresponds to the positive
orientation of~$K$ (see Example~\ref{exam:c-nbhd}).
Strictly speaking, if, as we shall always assume,
$S^1$ is parametrised by $0\leq \theta\leq 2\pi$, then this formula
for $\xi$ holds on $S^1\times D^2(\delta )$ for some, possibly small,
$\delta >0$. However, to simplify notation we usually work with $S^1\times
D^2$ as local model.

\begin{figure}[h]
\centerline{\relabelbox\small
\epsfxsize 8cm \epsfbox{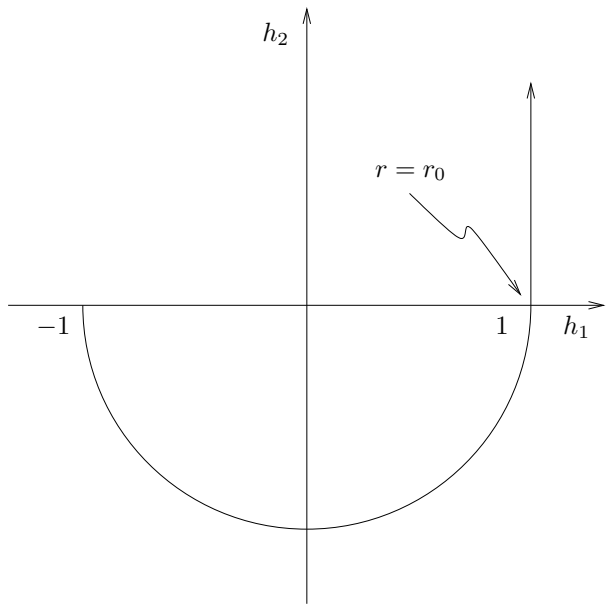}
\extralabel <-.6cm,3.6cm> {$h_1$}
\extralabel <-7.6cm,3.6cm> {$-1$}
\extralabel <-1.5cm,3.6cm> {$1$}
\extralabel <-4.6cm,7.5cm> {$h_2$}
\extralabel <-3.1cm,5.7cm> {$r=r_0$}
\endrelabelbox}
\caption{Lutz twist.}
\label{figure:lutz}
\end{figure}

We say that ${\xi}'$ is obtained from $\xi$ by a {\bf Lutz twist}
along $K$ and write ${\xi}'={\xi}^K$ if on $S^1\times D^2$ the new
contact structure ${\xi}'$ is defined by
\[ {\xi}'=\ker (h_1(r)\, d\theta +h_2(r)\, d\varphi ) \]
with $(h_1(r),h_2(r))$ as in Figure~\ref{figure:lutz},
and ${\xi}'$ coincides with $\xi$
outside $S^1\times D^2$.

More precisely, $(h_1(r),h_2(r))$ is required to satisfy the conditions
\begin{itemize}
\item[1.] $h_1(r)=-1$ and $h_2(r)=-r^2$ near $r=0$,
\item[2.] $h_1(r)=1$ and $h_2(r)=r^2$ near $r=1$,
\item[3.] $(h_1(r),h_2(r))$ is never parallel to $(h_1'(r),h_2'(r))$.
\end{itemize}
This is the same as applying the construction of Section~\ref{subsection3.1}
to the topologically trivial Dehn surgery given by the matrix
$\left( \begin{array}{cc}-1&0\\0&-1\end{array}\right)$.

\vspace{2mm}

It will be useful later on to understand more precisely the behaviour
of the map $f_{\xi'}\co S^3\rightarrow S^2$. For the definition of this
map we assume -- this assumption will be justified below --
that on $S^1\times D^2$ the map $f_{\xi}$ was defined in
terms of the standard metric $d\theta^2+du^2+dv^2$ (with $u,v$ cartesian
coordinates on~$D^2$ corresponding to the polar coordinates
$r,\varphi$) and the trivialisation $\partial_{\theta},
\partial_u,\partial_v$ of $T(S^1\times D^2)$.
Since $\xi '$ is spanned by $\partial_r$ and $h_2(r)\partial_\theta -
h_1(r)\partial_{\varphi}$ (resp.\ $\partial_u,\partial_v$ for $r=0$), a
vector positively orthogonal to $\xi '$ is given by 
\[ h_1(r)\partial_{\theta}+h_2(r)\partial_{\varphi}, \]
which makes sense even for $r=0$. Observe that the ratio $h_1(r)/h_2(r)$ (for
$h_2(r)\neq 0$) is a strictly monotone decreasing function since by
the third condition above we have
\[ (h_1/h_2)'=(h_1'h_2-h_1h_2')/h_2^2<0. \]
This implies that any value on $S^2$ other than $(1,0,0)$ (corresponding
to $\partial_{\theta}$) is regular for the map~$f_{\xi'}$; the value $(1,0,0)$
is attained along the torus $\{ r=r_0\}$, with $r_0>0$ determined by
$h_2(r_0)=0$, and hence not regular.

If $S^1\times D^2$ is endowed with the orientation defined by the volume
form $d\theta\wedge r\, dr\wedge d\varphi =d\theta \wedge du\wedge dv$
(so that $\xi$ and $\xi'$ are positive contact structures) and
$S^2\subset {\mathbb R}^3$ is given its `usual' orientation defined by
the volume form $x\, dy\wedge dz+y\, dz\wedge dx+z\, dx\wedge dy$, then
\[ f_{\xi'}^{-1}(-1,0,0)=S^1\times\{ 0\} \]
with orientation given by $-\partial_{\theta}$, since $f_{\xi'}$ maps the
slices $\{\theta\}\times D^2(r_0)$ orientation-reversingly onto~$S^2$.

More generally, for any
$p\in S^2-\{ (1,0,0)\}$ the preimage $f_{\xi'}^{-1}(p)$
(inside the domain $\{ (\theta ,r,\varphi )\co h_2(r)<0\}
=\{ r=r_0\}$) is a
circle $S^1\times\{ {\bf u}\}$, ${\bf u}\in D^2$, with orientation
given by~$-\partial_{\theta}$.

\vspace{2mm}

We are now ready to show how to construct a 
contact structure on $M$ in any given
homotopy class of $2$--plane distributions by starting with an arbitrary
contact structure and performing suitable Lutz twists. First we deal
with homotopy over the 2--skeleton. One way to proceed
would be to prove directly, with notation as above, that $d^2(\xi^K,\xi )=
-PD[K]$. However, it is somewhat easier to compute $d^2(\xi^K,\xi )$ in the 
case where $\xi$ is a trivial 2--plane bundle and the trivialisation
of $STM$ is adapted to~$\xi$. Since I would anyway like
to present an alternative argument for computing the effect of a Lutz twist
on the Euler class of the contact structure, and thus relate $d^2(\xi_1,\xi_2)$
with the Euler classes of $\xi_1$ and~$\xi_2$, it seems opportune to do this
first and use it to show the existence of a contact structure with Euler
class zero. In the next section we shall actually discuss a direct
geometric proof, due to Gonzalo, of the existence of a contact structure 
with Euler class zero.

Recall that the Euler class $e(\xi )\in H^2(B;{\Z})$ of a 2--plane
bundle over a complex $B$ (of arbitrary dimension) is the obstruction to
finding a nowhere zero section of $\xi$ over the 2--skeleton of~$B$. Since
$\pi_i(S^1)=0$ for $i\geq 2$, all higher obstruction groups $H^{i+1}(B;
\pi_i(S^1))$ are trivial, so a 2--dimensional orientable bundle $\xi$ is
trivial if and only if $e(\xi )=0$, no matter what the dimension of~$B$.

Now let $\xi$ be an arbitrary cooriented 2--plane distribution on an 
oriented 3--manifold~$M$. Then $TM\cong \xi\oplus\varepsilon^1$, where
$\varepsilon^1$ denotes a trivial line bundle. Hence $w_2(\xi )=w_2(
\xi\oplus\varepsilon^1)=w_2(TM)=0$, and since $w_2(\xi )$ is the
mod~2 reduction of $e(\xi )$ we infer that $e(\xi )$ has to be even.

\begin{prop}
For any even element $e\in H^2(M;{\Z})$ there is a contact structure $\xi$
on $M$ with $e(\xi )=e$.
\end{prop}

\begin{figure}[h]
\centerline{\relabelbox\small
\epsfxsize 10cm \epsfbox{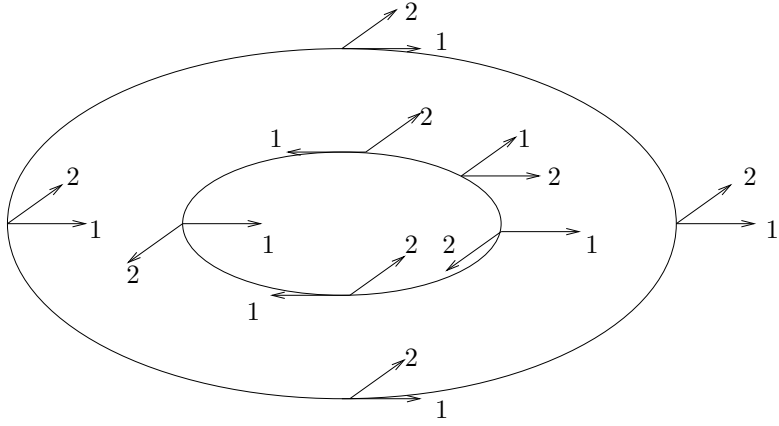}
\extralabel <-4.3cm,-.2cm> {$1$}
\extralabel <.1cm,2.2cm> {$1$}
\extralabel <-4.3cm,4.7cm> {$1$}
\extralabel <-8.9cm,2.2cm> {$1$}
\extralabel <-6.8cm,1.1cm> {$1$}
\extralabel <-2.3cm,2cm> {$1$}
\extralabel <-3.2cm,3.4cm> {$1$}
\extralabel <-6.5cm,3.4cm> {$1$}
\extralabel <-6.6cm,2cm> {$1$}
\extralabel <-4.7cm,.5cm> {$2$}
\extralabel <-.2cm,2.9cm> {$2$}
\extralabel <-4.7cm,5.1cm> {$2$}
\extralabel <-9.2cm,2.9cm> {$2$}
\extralabel <-4.7cm,2cm> {$2$}
\extralabel <-4.2cm,2cm> {$2$}
\extralabel <-2.8cm,2.9cm> {$2$}
\extralabel <-4.5cm,3.7cm> {$2$}
\extralabel <-8.4cm,1.6cm> {$2$}
\endrelabelbox}
\caption{Effect of Lutz twist on Euler class.}
\label{figure:euler}
\end{figure}

\begin{proof}
Start with an arbitrary contact structure $\xi_0$ on $M$ with $e(\xi_0)=e_0$
(which we know to be even). Given any even $e_1\in H^2(M;{\Z})$, represent
the Poincar\'e dual of $(e_0-e_1)/2$ by a collection of embedded
oriented circles
positively transverse to~$\xi_0$. (Here by $(e_0-e_1)/2$ I mean some
class whose double equals $e_0-e_1$; in the presence of $2$--torsion there
is of course a choice of such classes.)
Choose a section of $\xi_0$ transverse
to the zero section of~$\xi_0$, that is, a vector field in $\xi_0$ with 
generic zeros. We may assume that there are no zeros on the curves
representing $PD^{-1}(e_0-e_1)/2$. Now perform a Lutz twist as described
above along these curves and call $\xi_1$ the resulting contact
structure. It is easy to see that in the local model
for the Lutz twist a constant vector field tangent to
$\xi_0$ along $\partial (S^1\times D^2(r_0))$
extends to a vector field tangent to $\xi_1$ over
$S^1\times D^2(r_0)$ with zeros of index $+2$ along $S^1\times\{ 0\}$
(Figure~\ref{figure:euler}).
So the vector field in $\xi_0$ extends to a vector field in $\xi_1$
with new zeros of index $+2$ along the curves representing
$PD^{-1}(e_1-e_0)/2$ (notice that a Lutz twist along a positively
transverse knot $K$ turns $K$ into a negatively transverse knot).
Since the self-intersection class of $M$ in the total space of a vector
bundle is Poincar\'e dual to the Euler class of that bundle, this
proves $e(\xi_1)=e(\xi_0)+e_1-e_0=e_1.$
\end{proof}

We now fix a contact structure $\xi_0$ on $M$ with $e(\xi_0)=0$ and give
$M$ the orientation induced by~$\xi_0$ (i.e.\ the one for which
$\xi_0$ is a positive contact structure). Moreover,
we fix a Riemannian metric on $M$ and define $X_0$ as the vector field
positively orthonormal to~$\xi_0$. Since $\xi_0$ is a trivial
plane bundle, $X_0$ extends to an orthonormal frame $X_0,X_1,X_2$,
hence a trivialisation of $STM$,
with $X_1,X_2$ tangent to $\xi_0$ and defining the orientation of~$\xi_0$.
With these choices, $\xi_0$ corresponds to the constant map
$f_{\xi_0}\co M\rightarrow (1,0,0)\in S^2$.

\begin{prop}
Let $K\subset M$ be
an embedded, oriented circle positively transverse to $\xi_0$. Then
$d^2(\xi_0^K,\xi_0)=-PD[K]$.
\end{prop}

\begin{proof}
Identify a tubular neighbourhood of $K\subset M$ with $S^1\times D^2$
with framing defined by~$X_1$,
and $\xi_0$ given in this neighbourhood as the kernel of
$d\theta +r^2d\varphi =d\theta +u\, dv-v\, du$. We may then change
the trivialisation $X_0,X_1,X_2$ by a homotopy, fixed outside
$S^1\times D^2$, such that $X_0=\partial_{\theta}$, $X_1=\partial_u$
and $X_2=\partial_v$ near~$K$; this does not change the homotopical data
of $2$--plane distributions computed via the Pontrjagin-Thom
construction. Then $f_{\xi_0}$ is no longer constant, but its
image still does not contain the point $(-1,0,0)$.

Now perform a Lutz twist along $K\times\{ 0\}$. Our discussion at the
beginning of this section shows that $(-1,0,0)$ is a regular value of
the map $f_{\xi}\co M\rightarrow S^2$ associated with $\xi =\xi_0^K$
and $f_{\xi}^{-1}(-1,0,0)=-K$. Hence, by definition of the obstruction
class $d^2$ we have $d^2(\xi_0^K,\xi_0)=-PD[K]$.
\end{proof}

\begin{proof}[Proof of Theorem~\ref{thm:lutz-martinet}]
Let $\eta$ be a $2$--plane distribution on $M$ and $\xi_0$ the
contact structure on $M$ with $e(\xi_0)=0$ that we fixed earlier on.
According to our discussion in Section~\ref{section:rhcbs} and
Theorem~\ref{thm:approximation}, we can find an oriented knot $K$
positively transverse to $\xi_0$ with $-PD[K]=d^2(\eta ,\xi_0)$.
Then $d^2(\eta ,\xi_0)=d^2(\xi_0^K,\xi_0)$ by the preceding proposition,
and therefore $d^2(\xi_0^K,\eta )=0$.

We may then assume that $\eta =\xi_0^K$ on $M-D^3$, where we choose
$D^3$ so small that $\xi_0^K$ is in Darboux normal form there (and
identical with $\xi_0$).
By Proposition~\ref{prop:linking} we can find a link
$K'$ in $D^3$ transverse to $\xi_0^K$ with self-linking number
$l(K')$ equal to $d^3(\eta ,\xi_0^K)$.

Now perform a Lutz twist of $\xi_0^K$ along each component of $K'$ and let
$\xi$ be the resulting contact structure. Since this does not change
$\xi_0^K$ over the $2$--skeleton of~$M$, we still have $d^2(\xi ,\eta )=0$.

Observe that $f_{\xi_0^K}|_{D^3}$
does not contain the point $(-1,0,0)\in S^2$, and -- since
$f_{\xi_0^K}(D^3)$ is compact -- there is a whole neighbourhood $U\subset S^2$
of $(-1,0,0)$ not contained in $f_{\xi_0^K}(D^3)$.
Let $f\co S^3\rightarrow S^2$ be the map used to compute $d^3(\xi ,\xi_0^K)$,
that is, $f$ coincides on the upper hemisphere with
$f_{\xi}|_{D^3}$ and on the lower hemisphere with $f_{\xi_0^K}|_{D^3}$.
By the discussion in Section~\ref{subsection3.3}, the preimage
$f^{-1}(u)$ of any $u\in U-\{ (-1,0,0)\}$
will be a push-off of $-K'$ determined
by the trivialisation of $\xi_0^K|_{D^3}=\xi_0|_{D^3}$. So the linking number
of $f^{-1}(u)$ with $f^{-1}(-1,0,0)$, which is by definition the Hopf invariant
$H(f)=d^3(\xi ,\xi_0^K)$, will be equal to~$l(K')$. By our choice of
$K'$ and the additivity of $d^3$ this implies $d^3(\xi ,\eta )=0$.
So $\xi$ is a contact structure that is homotopic to $\eta$ as a
$2$--plane distribution.
\end{proof}
\subsection{Other existence proofs}
\label{subsection3.4}
Here I briefly summarise the other known existence proofs for
contact structures on $3$--manifolds, mostly by pointing to
the relevant literature. In spirit, most of these proofs
are similar to the one by Lutz-Martinet: start with a structure theorem
for $3$--manifolds and show that the topological construction can be
performed compatibly with a contact structure.
\subsubsection{Open books}
\begin{figure}[h]
\centerline{\relabelbox\small
\epsfxsize 5cm \epsfbox{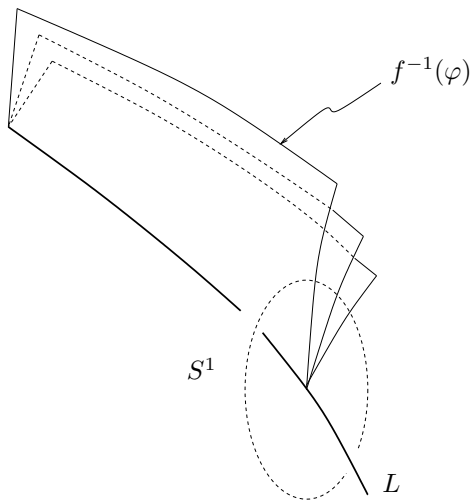}
\extralabel <0cm,.1cm> {$L$}
\extralabel <-2.6cm,1.6cm> {$S^1$}
\extralabel <.1cm,5.6cm> {$f^{-1}(\varphi )$}
\endrelabelbox}
\caption{An open book near the binding.}
\label{figure:book}
\end{figure}

According to a theorem of Alexander~\cite{alex23}, cf.~\cite{rolf76},
every closed, orientable
$3$--manifold $M$ admits an {\bf open book decomposition}.
This means that there
is a link $L\subset M$, called the {\bf binding}, and a fibration
$f\co M-L\rightarrow S^1$, whose fibres are called the
{\bf pages}, see Figure~\ref{figure:book}. It may be assumed that
$L$ has a tubular neighbourhood $L\times D^2$ such that $f$ restricted
to $L\times (D^2-\{ 0\})$ is given by $f(\theta ,r,\varphi )=\varphi$,
where $\theta$ is the coordinate along $L$ and $(r,\varphi )$ are polar
coordinates on~$D^2$.

At the cost of raising the genus of the pages, one may decrease the number
of components of~$L$, and in particular one may always assume $L$ to
be a knot. Another way to think of such an open book is as follows.
Start with a surface $\Sigma$ with boundary $\partial\Sigma =K\cong S^1$ and
a self-diffeomorphism $h$ of~$\Sigma$ with $h=\mbox{\rm id}$
near~$K$. Form the mapping torus
$T_h=\Sigma_h=\Sigma\times [0,2\pi ]/\!\!\sim$, where `$\sim$' denotes the
identification $(p,2\pi )\sim (h(p),0)$. Define a $3$--manifold $M$ by
\[ M=T_h\cup_{K\times S^1} (K\times D^2).\]
This $M$ carries by construction the structure of an open book
with binding $K$ and pages diffeomorphic to~$\Sigma$.

Here is a slight variation on a simple argument of Thurston and
Winkelnkemper~\cite{thwi75} for producing a contact structure on
such an open book (and hence on any closed, orientable $3$--manifold):

Start with a $1$--form $\beta_0$ on $\Sigma$ with $\beta_0=e^td\theta$
near $\partial\Sigma =K$, where $\theta$
denotes the coordinate along $K$
and $t$ is a collar parameter with $K=\{ t=0\}$ and $t<0$ in the interior
of~$\Sigma$. Then $d\beta_0$ integrates to $2\pi$ over $\Sigma$ by
Stokes's theorem. Given any area form $\omega$ on $\Sigma$ (with total
area equal to~$2\pi$) satisfying $\omega =e^tdt\wedge d\theta$
near~$K$, the $2$--form $\omega -d\beta_0$ is, by de Rham's
theorem, an exact $1$--form, say $d\beta_1$, where we may assume $\beta_1
\equiv 0$ near~$K$.

Set $\beta =\beta_0+\beta_1$. Then $d\beta =\omega$ is an area form
(of total area~$2\pi$) on $\Sigma$
and $\beta =e^td\theta$ near~$K$. The set of $1$--forms satisfying these two
properties is a convex set, so we can find a $1$--form (still denoted~$\beta$)
on $T_h$ which has these properties when restricted to the fibre
over any $\varphi\in S^1=[0,2\pi ]/_{0\sim 2\pi}$. We may (and shall)
require that $\beta =e^td\theta$ near~$\partial T_h$.

Now a contact form $\alpha$ on $T_h$ is found by setting
$\alpha =\beta +C\, d\varphi$ for a sufficiently large constant
$C\in{\mathbb R}^+$, so that in
\[ \alpha\wedge d\alpha = (\beta +C\, d\varphi )\wedge d\beta\]
the non-zero term $d\varphi\wedge d\beta =d\varphi\wedge\omega$ dominates.
This contact form can be extended to all of $M$ by
making the ansatz
$\alpha =h_1(r)d\theta +h_2(r)d\varphi$ on $K\times D^2$, as described
in our discussion of the Lutz twist. The boundary conditions in the
present situation are, say,
\begin{itemize}
\item[1.] $h_1(r)=2$ and $h_2(r)=r^2$ near $r=0$,
\item[2.] $h_1(r)=e^{1-r}$ and $h_2(r)=C$ near $r=1$.
\end{itemize}

Observe that subject to these boundary conditions
a curve $(h_1(r),h_2(r))$
can be found that does not pass the
$h_2$--axis (i.e.\ with $h_1(r)$ never being equal to zero).
In the $3$--dimensional setting this is not essential (and the
Thurston-Winkelnkemper ansatz lacked that feature), but it is crucial
when one tries to generalise this construction to higher dimensions.
This has recently been worked out by Giroux and J.-P.~Mohsen~\cite{gimo}.
This, however, is only the easy part of their work. Rather strikingly,
they have shown that a converse of this result holds: Given a compact contact
manifold of arbitrary dimension, it admits an open book decomposition
that is adapted to the contact structure in the way described above.
Full details have not been published at the time of writing, but see
Giroux's talk~\cite{giro02} at the ICM 2002.
\subsubsection{Branched covers}
\label{section:branched}
A theorem of Hilden, Montesinos and Thickstun~\cite{hmt76} states
that every closed, orientable $3$--manifold $M$ admits a branched
covering $\pi\co M\rightarrow S^3$ such that the upstairs branch
set is a simple closed curve that bounds an embedded disc.
(Moreover, the cover can be chosen $3$--fold and simple, i.e.\
the monodromy representation of $\pi_1(S^3-K)$, where $K$ is the
branching set downstairs (a knot in~$S^3$), represents the meridian
of $K$ by a transposition in the symmetric group~$S_3$. This, however,
is not relevant for our discussion.)

It follows immediately, as announced in Section~\ref{subsection3.2},
that every closed, orientable $3$--manifold is parallelisable: First
of all, $S^3$ is parallelisable since it carries a Lie group
structure (as the unit quaternions, for instance). Given $M$ and a
branched covering $\pi\co M\rightarrow S^3$ as above,
there is a $3$--ball $D^3\subset M$ containing the upstairs branch
set. Outside of $D^3$, the covering $\pi$ is unbranched, so the
$3$--frame on $S^3$ can be lifted to a frame on $M-D^3$.
The bundle $TM|_{D^3}$ is trivial, so the frame
defined along
$\partial D^3$ defines an element of~$\SO (3)$ (cf.\ the footnote in the proof
of Theorem~\ref{thm:parallelisable}).
Since $\pi_2(\SO (3))=0$,
this frame extends over~$D^3$.

In~\cite{gonz87}, Gonzalo uses this theorem to construct a
contact structure on every closed, orientable $3$--manifold $M$,
in fact one with zero Euler class:
Away from the branching set one can lift the standard
contact structure from~$S^3$ (which has Euler class zero: a
trivialisation is given by two of the three (quaternionic)
Hopf vector fields). A careful analysis of the
branched covering map near the branching set then shows how
to extend this contact structure over~$M$ (while keeping it trivial
as $2$--plane bundle).

A branched covering construction for higher-dimensional contact
manifolds is discussed in~\cite{geig97}.
\subsubsection{. . . and more}
The work of Giroux~\cite{giro91}, in which he initiated the study
of convex surfaces in contact $3$--manifolds, also contains a proof
of Martinet's theorem.

An entirely different proof, due to S. Altschuler~\cite{alts95},
finds contact structures from solutions to a certain parabolic
differential equation for $1$--forms on $3$--manifolds. Some of these
ideas have entered into the more far-reaching work of
Eliashberg and Thurston on so-called `confoliations'~\cite{elth98}, that is,
$1$--forms satisfying $\alpha\wedge d\alpha\geq 0$.

\subsection{Tight and overtwisted}
\label{subsection3.5}
The title of this section describes the
fundamental dichotomy of contact structures
in dimension~$3$ that has proved seminal for the development
of the field.

\begin{figure}[h]
\centerline{\relabelbox\small
\epsfxsize 8cm \epsfbox{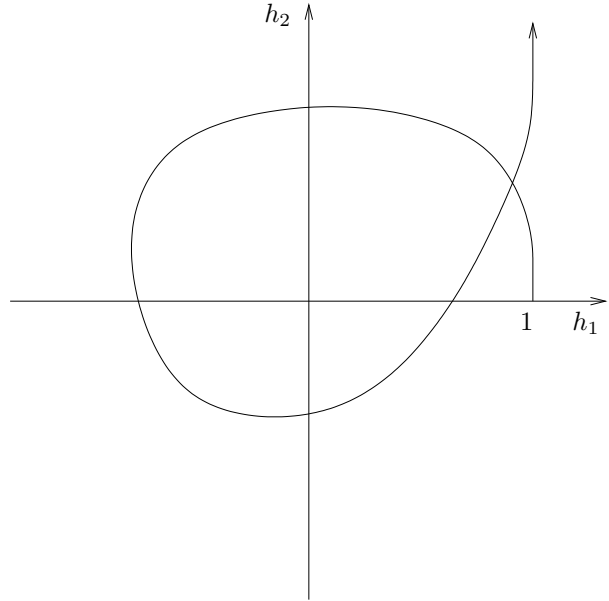}
\extralabel <-.5cm,3.6cm> {$h_1$}
\extralabel <-1.2cm,3.6cm> {$1$}
\extralabel <-4.6cm,7.7cm> {$h_2$}
\endrelabelbox}
\caption{A full Lutz twist.}
\label{figure:full-lutz}
\end{figure}

In order to motivate the notion of an overtwisted contact structure,
as introduced by Eliashberg~\cite{elia89}, we consider a `full'
Lutz twist as follows.
Let $(M,\xi )$ be a 
contact $3$--manifold and $K\subset M$ a knot transverse to~$\xi$. As
before, identify $K$ with $S^1\times\{ 0\}\subset S^1\times D^2\subset M$
such that $\xi =\ker (d\theta +r^2d\varphi )$ on $S^1\times D^2$. Now
define a new contact structure $\xi'$ as in Section~\ref{subsection3.3},
with $(h_1(r),h_2(r))$
now as in Figure~\ref{figure:full-lutz},
that is, the boundary conditions are now
\begin{itemize}
\item[] $h_1(r)=1$ and $h_2(r)=r^2$ for $r\in [0,\varepsilon ]\cup
[1-\varepsilon ,1]$
\end{itemize}
for some small $\varepsilon >0$.

\begin{lem}
A full Lutz twist does not change the homotopy class of $\xi$ as a $2$--plane
field.
\end{lem}

\begin{proof}
Let $(h_1^t(r),h_2^t(r))$, $r,t\in [0,1]$, be a homotopy of paths such
that
\begin{itemize}
\item[1.] $h_1^0\equiv 1, h_2^0(r)=r^2$,
\item[2.] $h_i^1\equiv h_i$, $i=1,2$,
\item[3.] $h_i^t(r)=h_i(r)$ for $r\in [0,\varepsilon ]\cup [1-\varepsilon ,1]$.
\end{itemize}
Let $\chi :[0,1]\rightarrow {\bf R}$ be a smooth function which is identically
zero near $r=0$ and $r=1$ and $\chi (r)>0$ for $r\in [\varepsilon ,1-
\varepsilon ]$. Then 
\[ \alpha_t=t(1-t)\chi (r)\, dr+h_1^t(r)\, d\theta +h_2^t(r)\, d\varphi \]
is a homotopy from $\alpha_0= d\theta+r^2d\varphi$ to $\alpha_1=
h_1(r)\, d\theta +h_2(r)\, d\varphi$ through non-zero $1$--forms. This homotopy
stays fixed near $r=1$, and so it defines a homotopy between $\xi$ and
$\xi'$ as $2$--plane fields.
\end{proof}

\begin{figure}[h]
\centerline{\relabelbox\small
\epsfxsize 6cm \epsfbox{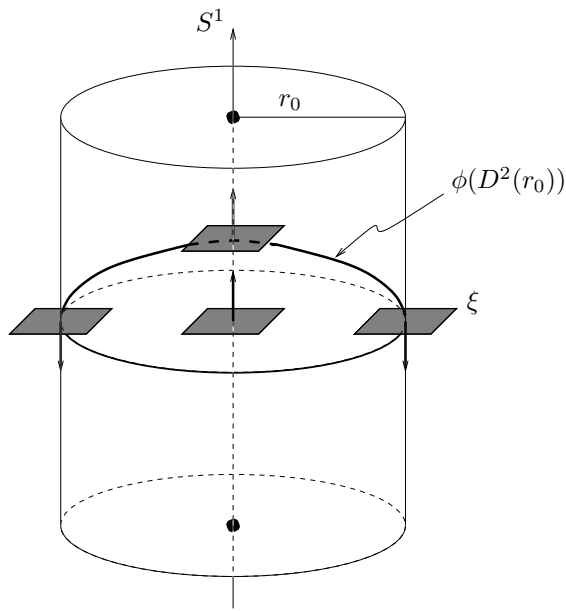}
\extralabel <.1cm,4cm> {$\xi$}
\extralabel <-0.1cm,5.6cm> {$\phi(D^2(r_0))$}
\extralabel <-3.5cm,7.7cm> {$S^1$}
\extralabel <-2.4cm,6.7cm> {$r_0$}
\endrelabelbox}
\caption{An overtwisted disc.}
\label{figure:overtwisted2}
\end{figure}

\begin{figure}[h]
\centerline{\relabelbox\small
\epsfxsize 5cm \epsfbox{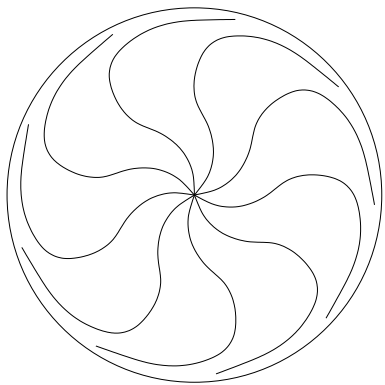}
\endrelabelbox}
\caption{Characteristic foliation on an overtwisted disc.}
\label{figure:overtwisted}
\end{figure}

Let $r_0$ be the smaller of the two positive radii with $h_2(r_0)=0$
and consider the embedding
\[ \begin{array}{rccc}
\phi : & D^2(r_0) & \longrightarrow & S^1\times D^2\\
 & (r,\varphi ) & \longmapsto & (\theta (r), r,\varphi ),
\end{array} \]
where $\theta (r)$ is a smooth function with $\theta (r_0)=0$, $\theta (r)>0$ 
for $0\leq r<r_0$, and $\theta'(r)=0$ only for $r=0$. We may require in
addition that $\theta (r)=\theta (0)-r^2$ near $r=0$. Then
\[ \phi^*(h_1(r)\, d\theta +h_2(r)\, d\varphi )=h_1(r)\theta'(r)\, dr+
h_2(r)\, d\varphi \]
is a differential
$1$--form on $D^2(r_0)$ which vanishes only for $r=0$, and along
$\partial D^2(r_0)$ the vector field $\partial_{\varphi}$ tangent to
the boundary lies in the kernel of this $1$--form,
see Figure~\ref{figure:overtwisted2}. In other words, the contact
planes $\ker (h_1(r)\, d\theta +h_2(r)\, d\varphi )$ intersected with
the tangent planes to the embedded disc $\phi (D^2(r_0))$ induce a singular
$1$--dimensional foliation on this disc with the boundary of this disc as
closed leaf and precisely one singular point in the interior of the disc
(Figure~\ref{figure:overtwisted};
notice that the leaves of this foliation are the integral curves of 
the vector field $h_1(r)\theta'(r)\, \partial_{\varphi}-h_2(r)\,
\partial_r$). Such a disc is called an {\bf overtwisted disc}.

A contact structure $\xi$ on a $3$--manifold $M$ is called {\bf overtwisted}
if $(M,\xi )$ contains an embedded overtwisted disc.
In order to justify this terminology,
observe that in the radially symmetric standard contact structure
of Example~\ref{exam:standard2}, the angle by which
the contact planes turn approaches $\pi /2$ asymptotically as
$r$ goes to infinity. By contrast, any contact manifold which
has been constructed using at least one (simple) Lutz twist
contains a similar cylindrical region where the contact planes
twist by more than $\pi$ in radial direction (at the smallest positive radius
$r_0$ with $h_2(r_0)=0$ the twisting angle has reached~$\pi$).

We have shown the following

\begin{prop}
Let $\xi$ be a contact structure on $M$. By a full Lutz twist along
any transversely embedded circle one obtains an overtwisted contact
structure $\xi'$ that is homotopic to $\xi$ as a $2$--plane distribution.
\hfill $\Box$
\end{prop}

Together with the theorem of Lutz and Martinet we find that $M$ contains 
an {\em overtwisted} contact structure in every homotopy class
of $2$--plane distributions. In fact, Eliashberg~\cite{elia89} has proved the
following much stronger theorem.

\begin{thm}[Eliashberg]
\label{thm:class-ot}
On a closed, orientable $3$--manifold there is a one-to-one correspondence
between homotopy classes of overtwisted contact structures and homotopy
classes of $2$--plane distributions.
\end{thm}

This means that two overtwisted contact structures which are homotopic
as $2$--plane fields are actually homotopic as contact structures
and hence isotopic by Gray's stability theorem.

\vspace{2mm}

Thus, it `only' remains to classify contact structures that are not
overtwisted. In~\cite{elia92} Eliashberg defined {\bf tight} contact structures
on a $3$--manifold $M$ as contact structures $\xi$ for which there is
no embedded disc $D\subset M$ such that $D_{\xi}$ contains a limit cycle.
So, by definition, overtwisted contact structures are not tight.
In that same paper, as mentioned above in
Section~\ref{section:nbhd-applications}, Eliashberg goes on to show the
converse with the help of the Elimination Lemma: non-overtwisted contact
structures are tight.

There are various ways to detect whether a contact structure is tight.
Historically the first proof that a certain contact structure is tight is
due to D.~Bennequin~\cite[Cor.~2, p.~150]{benn83}:

\begin{thm}[Bennequin]
The standard contact structure $\xi_0$ on $S^3$ is tight.
\end{thm}

The steps of the proof are as follows: 
(i) First, Bennequin shows that if
$\gamma_0$ is a transverse knot in $(S^3,\xi_0)$ with Seifert
surface $\Sigma$, then the self-linking number of $\gamma$ satisfies
the inequality
\[ l(\gamma_0 )\leq -\chi (\Sigma ).\]

(ii) Second, he introduces an invariant for Legendrian knots; nowadays this is
called the {\bf Thurston-Bennequin invariant}: Let $\gamma$ be a Legendrian
knot in $(S^3,\xi_0)$. Take a vector field $X$ along $\gamma$
transverse to~$\xi_0$, and let $\gamma'$ be the push-off of $\gamma$
in the direction of~$X$. Then the Thurston-Bennequin invariant $tb(\gamma)$
is defined to be the linking number of $\gamma$ and $\gamma'$.
(This invariant has an extension to homologically trivial Legendrian knots
in arbitrary contact $3$--manifolds.)

(iii) By pushing $\gamma$ in the direction of $\pm X$, one obtains
transverse curves $\gamma^{\pm}$ (either of which is a candidate for
$\gamma'$ in~(ii)). One of these curves is positively transverse,
the other negatively transverse to~$\xi_0$. The self-linking number
of $\gamma^{\pm}$ is related to the Thurston-Bennequin invariant
and a further invariant (the rotation number) of~$\gamma$. The equation
relating these three invariants implies $tb (\gamma )\leq -\chi (\Sigma )$,
where $\Sigma$ again denotes a Seifert surface for~$\gamma$. In particular,
a Legendrian unknot $\gamma$ satisfies $tb (\gamma )<0$. This inequality
would be violated by the vanishing cycle of an overtwisted disc
(which has $tb=0$), which
proves that $(S^3,\xi_0)$ is tight.

\begin{rem}
{\rm (1) Eliashberg~\cite{elia93} generalised the Bennequin inequality
$l(\gamma_0 )\leq-\chi (\Sigma )$ for transverse knots
(and the corresponding inequality for the Thurston-Bennequin invariant
of Legendrian knots) to arbitrary tight contact $3$--manifolds. Thus, whereas
Bennequin established the tightness (without that name) of the standard
contact structure on $S^3$
by proving the inequality that bears his name, that inequality is
now seen, conversely, as a consequence of tightness.

(2) In~\cite{benn83} Bennequin denotes the positively (resp.\
negatively) transverse push-off of the Legendrian knot $\gamma$ by
$\gamma^-$ (resp.~$\gamma^+$). This has led to some sign errors
in the literature. Notably, the $\pm$ in Proposition~2.2.1
of~\cite{elia93}, relating the described invariants of $\gamma$
and $\gamma^{\pm}$, needs to be reversed.
}
\end{rem}

\begin{cor}
The standard contact structure on $\R^3$ is tight.
\end{cor}

\begin{proof}
This is immediate from Proposition~\ref{prop:standard-S-R}.
\end{proof}

Here are further tests for tightness:

\vspace{2mm}

{\bf 1.} A closed contact $3$--manifold $(M,\xi )$ is called
{\bf symplectically
fillable} if there exists a compact symplectic manifold $(W,\omega )$
bounded by $M$ such that
\begin{itemize}
\item the restriction of $\omega$ to $\xi$ does not vanish anywhere,
\item the orientation of $M$ defined by $\xi$ (i.e.\ the one for which
$\xi$ is positive) coincides with the orientation of $M$ as boundary
of the symplectic manifold $(W,\omega )$ (which is oriented by~$\omega^2$).
\end{itemize}

We then have the following result of Eliashberg~\cite[Thm.~3.2.1]{elia88},
\cite{elia90} and
Gromov~\cite[2.4.$D_2'$(b)]{grom85}, cf.~\cite{benn90}:

\begin{thm}[Eliashberg-Gromov]
A symplectically fillable contact structure is tight.
\end{thm}

\begin{exam}
{\rm The $4$--ball $D^4\subset\R^4$ with symplectic form $\omega =
dx_1\wedge dy_1+dx_2\wedge dy_2$ is a symplectic filling of
$S^3$ with its standard contact structure~$\xi_0$. This gives an
alternative proof of Bennequin's theorem.
}
\end{exam}

\vspace{2mm}

{\bf 2.} Let $(\widetilde{M},\tilde{\xi})\rightarrow (M,\xi )$ be a
covering map and contactomorphism. If $(\widetilde{M},\tilde{\xi})$ is
tight, then so is $(M,\xi )$, for any overtwisted disc in $(M,\xi )$
would lift to an overtwisted disc in $(\widetilde{M},\tilde{\xi})$.

\begin{exam}
{\rm The contact structures $\xi_n$, $n\in{\mathbb N}$,
on the $3$--torus $T^3$ defined by
\[ \alpha_n=\cos (n\theta_1 )\, d\theta_2+\sin (n\theta_1)\, d\theta_3=0\]
are tight: Lift the contact structure $\xi_n$ to the universal cover
$\R^3$ of $T^3$; there the contact structure is defined by the same equation
$\alpha_n=0$,
but now $\theta_i\in\R$ instead of $\theta_i\in\R /2\pi\Z \cong S^1$.
Define a diffeomorphism $f$ of $\R^3$ by
\[ f(x,y,z)=(y/n,z\,\cos y+x\,\sin y,z\,\sin y-x\,\cos y)=:
(\theta_1,\theta_2,\theta_3).\]
Then $f^*\alpha_n=dz+x\, dy$, so the lift of $\xi_n$ to $\R^3$ is
contactomorphic to the tight standard contact structure on~$\R^3$.
}
\end{exam}

Notice that it is possible for a tight contact structure to be finitely
covered by an overtwisted contact structure. The first such examples
were due to S.~Makar-Limanov~\cite{maka98}. Other examples of this kind
have been found by V.~Colin~\cite{coli99} and R.~Gompf~\cite{gomp98}.

\vspace{2mm}

{\bf 3.} The following theorem of H.~Hofer~\cite{hofe93}
relates the dynamics of
the Reeb vector field to overtwistedness.

\begin{thm}[Hofer]
Let $\alpha$ be a contact form on a closed $3$--manifold such
that the contact structure $\ker\alpha$ is overtwisted. Then the
Reeb vector field of $\alpha$ has at least one contractible periodic orbit.
\end{thm}

\begin{exam}
{\rm The Reeb vector field $R_n$
of the contact form $\alpha_n$ of the preceding example is
\[ R_n =\cos (n\theta_1)\,\partial_{\theta_2}+\sin (n\theta_1)\,
\partial_{\theta_3}.\]
Thus, the orbits of $R_n$ define constant slope foliations of the $2$--tori
$\{ \theta_1=\mbox{\rm const.}\}$; in particular, the periodic orbits
of $R_n$ are even homologically non-trivial. It follows, again,
that the $\xi_n$ are tight contact structures on~$T^3$.
(This, admittedly, amounts to attacking starlings with rice puddings
fired from catapults\footnote{This turn of phrase
originates from~\cite{mill54}.}.)
}
\end{exam}

\subsection{Classification results}
\label{subsection3.6}
In this section I summarise some of the known classification
results for contact structures on $3$--manifolds. By Eliashberg's
Theorem~\ref{thm:class-ot} it suffices to list the tight contact
structures, up to isotopy or diffeomorphism, on a given
closed $3$--manifold.

\begin{thm}[Eliashberg~\cite{elia92}]
Any tight contact structure on $S^3$ is isotopic to the
standard contact structure~$\xi_0$.
\end{thm}

This theorem has a remarkable application in differential topology,
viz., it leads to a new proof of Cerf's theorem~\cite{cerf68} that any
diffeomorphism of $S^3$ extends to a diffeomorphism
of the $4$--ball~$D^4$. The idea is that the above theorem implies that
any diffeomorphism of $S^3$ is isotopic to a contactomorphism
of~$\xi_0$. Eliashberg's technique~\cite{elia90}
of filling by holomorphic discs
can then be used to show that such a contactomorphism extends
to a diffeomorphism of~$D^4$.

As remarked earlier (Remark~\ref{rem:gray}), Eliashberg has also
classified contact structures on~$\R^3$. Recall that homotopy classes
of $2$--plane distributions on $S^3$ are classified by $\pi_3(S^2)\cong
\Z$. By Theorem~\ref{thm:class-ot}, each of these classes contains
a unique (up to isotopy) overtwisted contact structure. When a point
of $S^3$ is removed, each of these contact structures induces one on~$\R^3$,
and Eliashberg~\cite{elia93} shows that they remain non-diffeomorphic
there. Eliashberg shows further that,
apart from this integer family of overtwisted contact
structures, there is a unique tight contact structure on~$\R^3$
(the standard one), and a single overtwisted one that is `overtwisted
at infinity' and cannot be compactified to a contact structure
on~$S^3$.

\vspace{2mm}

In general, the classification of contact structures up to
diffeomorphism will differ from the classification up to isotopy.
For instance, on the $3$--torus $T^3$ we have the following
diffeomorphism classification due to Y.~Kanda~\cite{kand97}:

\begin{thm}[Kanda]
Every (positive) tight contact structure on $T^3$ is contactomorphic to one
of the $\xi_n$, $n\in{\mathbb N}$, described above. For $n\neq m$,
the contact structures $\xi_n$ and $\xi_m$ are not contactomorphic.
\end{thm}

Giroux~\cite{giro94} had proved earlier that $\xi_n$ for $n\geq 2$ is not
contactomorphic to~$\xi_1$.

On the other hand, all the $\xi_n$ are homotopic as $2$--plane
fields to $\{ d\theta_1=0\}$.
This shows one way how Eliashberg's classification Theorem~\ref{thm:class-ot}
for overtwisted contact structures can fail for tight contact structures:

\begin{itemize}
\item There are tight contact structures on $T^3$ that are homotopic as plane
fields but not contactomorphic.
\end{itemize}

P.~Lisca and G.~Mati\'c~\cite{lima97} have found examples of the same kind on
homology spheres by applying Seiberg-Witten theory to Stein fillings
of contact manifolds, cf.\ also~\cite{krmr97}.

Eliashberg and L.~Polterovich~\cite{elpo94} have determined the
isotopy classes of diffeomorphisms of $T^3$ that contain a contactomorphism
of~$\xi_1$: they correspond to exactly those elements
of $\SL (3,\Z )=\pi_0(\mbox{\rm Diff}(T^3))$ that stabilise the
subspace $0\oplus\Z^2$ corresponding to the coordinates
$(\theta_2,\theta_3)$. In combination with Kanda's result, this allows
to give an isotopy classification of tight contact structures on~$T^3$.
One particular consequence of the result of Eliashberg and Polterovich
is the following:

\begin{itemize}
\item There are tight contact structures on $T^3$
that are contactomorphic and
homotopic as plane fields, but not isotopic (i.e.\ not
homotopic through contact structures).
\end{itemize}

Again, such examples also exist on homology spheres, as
S.~Akbulut and R.~Matveyev~\cite{akma98} have shown.

Another aspect of Eliashberg's classification of overtwisted contact
structures that fails to hold for tight structures is of course
the existence of such a structure in every homotopy class of $2$--plane
fields, as is already demonstrated by the classification of
contact structures on~$S^3$. Etnyre and K.~Honda~\cite{etho01} have
recently even found an example of a manifold -- the connected some
of two copies of the Poincar\'e sphere with opposite orientations --
that does not admit any tight contact structure at all.

For the classification of tight contact structures on lens spaces and
$T^2$--bundles over $S^1$ see \cite{giro00}, \cite{hond00}
and~\cite{hond00a}. A partial classification of tight contact structures
on lens spaces had been obtained earlier in~\cite{etny00}.

As further reading on $3$--dimensional contact geometry I can recommend
the lucid Bourbaki talk by Giroux~\cite{giro93}. This covers the
ground up to Eliashberg's classification of overtwisted contact structures
and the uniqueness of the tight contact structure on~$S^3$.

%% file: section4.tex
\section{A guide to the literature}
In this concluding section I give some recommendations for further reading,
concentrating on those aspects of contact geometry that have not
(or only briefly) been touched upon in earlier sections.

Two general surveys that emphasise historical matters and describe
the development of contact geometry from some of its earliest
origins are the one by Lutz~\cite{lutz88} and one by
the present author~\cite{geig01}.

One aspect of contact geometry that I have neglected in these
notes is the Riemannian geometry of contact manifolds (leading, for instance,
to Sasakian geometry). The survey by Lutz
has some material on that; D.~Blair~\cite{blai02}
has recently published a monograph on the topic.

There have also been various ideas for defining
interesting families of contact structures. Again, the survey by Lutz
has something to say on that; one such concept that has exhibited
very intriguing ramifications -- if this commercial break be 
permitted -- was introduced in~\cite{gego95}.
\subsection{Dimension $3$}
As mentioned earlier, Chapter~8 in~\cite{abklr94} is in parts complementary
to the present notes and has some material on surfaces in contact
$3$--manifolds. Other surveys and introductory texts on $3$--dimensional
contact geometry are the introductory lectures by
Etnyre~\cite{etny-lectures} and the Bourbaki talk
by Giroux~\cite{giro91}. Good places to
start further reading are two papers by Eliashberg:
\cite{elia92} for the classification of tight
contact structures and \cite{elia93a} for knots in contact $3$--manifolds.
Concerning the latter, there is also a chapter in preparation by
Etnyre~\cite{etny-HB} for a companion {\it Handbook} and an article by Etnyre
and Honda~\cite{etho01a} with an extensive introduction to that subject.

The surveys~\cite{elia88} and~\cite{elia98} by
Eliashberg are more general in scope,
but also contain material about contact $3$--manifolds.

$3$--dimensional contact topology has now become a fairly wide arena;
apart from the work of Eliashberg, Giroux, Etnyre-Honda and others described
earlier, I should also mention the results of Colin, who has, for instance,
shown that surgery of index one (in particular: taking the
connected sum) on a tight contact $3$--manifold leads again to
a tight contact structure~\cite{coli97}.

Finally, Etnyre and L.~Ng~\cite{etng-problems} have compiled a useful list
of problems in $3$--dimensional contact topology.
\subsection{Higher dimensions}
The article~\cite{geig01a} by the present author contains a survey
of what was known at the time of writing about the existence of contact
structures on higher-dimensional manifolds. One of the most important
techniques for constructing contact manifolds in higher dimensions
is the so-called contact surgery along isotropic spheres developed by
Eliashberg~\cite{elia90a} and A.~Weinstein~\cite{wein91}. The latter
is a very readable paper. For a recent application of this
technique see~\cite{geth01}. Other constructions of contact manifolds
(branched covers, gluing along codimension~$2$ contact
submanifolds) are described in my paper~\cite{geig97}.

Odd-dimensional tori are of course amongst the manifolds with
the simplest global description, but they do not easily lend
themselves to the construction of contact structures.
In~\cite{lutz79} Lutz found a contact structure on~$T^5$; since then
it has been one of the prize questions in contact geometry to find
a contact structure on higher-dimensional tori. That prize,
as it were, recently went to F.~Bourgeois~\cite{bour02}, who showed
that indeed all odd-dimensional tori do admit a contact structure.
His construction uses the result of Giroux and
Mohsen~\cite{giro02},~\cite{gimo}
about open book decompositions adapted to contact structures
in conjunction with the original proof of Lutz. With the help of the
branched cover theorem described in~\cite{geig97} one can conclude
further that every manifold of the form $M\times {\Sigma}$ with
$M$ a contact manifold and $\Sigma$ a surface of genus at least~$1$
admits a contact structure.

Concerning the classification of contact structures in higher
dimensions, the first steps have been taken by Eliashberg~\cite{elia98a}
with the development of contact homology, which has been taken
further in~\cite{egh00}. This has been used by I.~Ustilovsky~\cite{usti99}
to show that on $S^{4n+1}$ there exist infinitely many non-iso\-mor\-phic
contact structures that are homotopically equivalent (in the
sense that they induce the same almost contact structure, i.e.\
reduction of the structure group
of $TS^{4n+1}$ to $1\times \mbox{\rm U} (2n)$). Earlier results
in this direction can be found in~\cite{geig97a} in the context of
various applications of contact surgery.
\subsection{Symplectic fillings}
A survey on the various types of symplectic fillings of contact
manifolds is given by Etnyre~\cite{etny98}, cf.\ also the
survey by Bennequin~\cite{benn90}. Etnyre and Honda~\cite{etho02}
have recently shown that certain Seifert fibred $3$--manifolds $M$
admit tight contact structures $\xi$ that are not symplectically
semi-fillable, i.e.\ there is no symplectic filling $W$
of $(M,\xi )$ even if $W$ is allowed to have other contact boundary
components. That paper also contains a good update on the general
question of symplectic fillability.

A related question is whether symplectic manifolds can have disconnected
boundary of contact type (this corresponds to a stronger notion of
symplectic filling defined via a Liouville vector field transverse
to the boundary and pointing outwards).
For (boundary) dimension~$3$ this is discussed
by D.~McDuff~\cite{mcdu91}; higher-dimensional symplectic manifolds with
disconnected boundary of contact type have been constructed in~\cite{geig94}.
\subsection{Dynamics of the Reeb vector field}
In a seminal paper, Hofer~\cite{hofe93} applied the method
of pseudo-holomorphic curves, which had been introduced to symplectic
geometry by Gromov~\cite{grom85}, to solve (for large
classes of contact $3$--manifolds) the so-called
Weinstein conjecture~\cite{wein79} concerning the existence
of periodic orbits of the Reeb vector field of a given contact form.
(In fact, one of the remarkable aspects of Hofer's work is that in
many instances it shows the existence of a periodic orbit of the
Reeb vector field of any contact form defining a given contact
structure.) A Bourbaki talk on the state of the art around the time
when Weinstein formulated the conjecture that bears his name
was given by N.~Desolneux-Moulis~\cite{deso81}; another Bourbaki talk
by F.~Laudenbach describes Hofer's contribution to the problem.

The textbook by Hofer and E.~Zehnder~\cite{hoze94} addresses these issues,
although its main emphasis, as is clear from the title,
lies more in the direction of symplectic geometry and Hamiltonian dynamics.
Two surveys by Hofer~\cite{hofe99}, \cite{hofe00}, and one by
Hofer and M.~Kriener~\cite{hokr99}, are more directly concerned
with contact geometry. Of the three, \cite{hofe99} may be the
best place to start, since it derives from a course of five
lectures. In collaboration with K.~Wysocki and Zehnder, Hofer
has expanded his initial ideas into a far-reaching project
on the characterisation of contact manifolds via the dynamics of the
Reeb vector field, see e.g.~\cite{hwz95}.